\documentclass[11pt,a4paper]{amsart}

\usepackage{amsmath, amssymb, amsthm}
\usepackage{stmaryrd} 
\usepackage{mathrsfs}
\usepackage[unicode]{hyperref}
\usepackage{geometry}
\usepackage{graphicx}
\usepackage{enumitem}
\usepackage{xcolor}
\usepackage{tikz}
\usetikzlibrary{calc, angles, quotes, arrows.meta,decorations.pathreplacing}
\usepackage{tikzsymbols}

\usepackage{xifthen}

\definecolor{slice}{RGB}{80,65,130}
\definecolor{velocity}{RGB}{185,48,51}
\definecolor{tangent}{RGB}{0,115,125}

\newtheorem{theorem}{Theorem}[section]
\newtheorem{lemma}[theorem]{Lemma}
\newtheorem{corollary}[theorem]{Corollary}
\newtheorem{proposition}[theorem]{Proposition}

\theoremstyle{definition}
\newtheorem{definition}[theorem]{Definition}
\newtheorem{example}[theorem]{Example}
\newcommand{\continuedexamplename}{}
\newtheorem*{continuedexampleinner}{\continuedexamplename}
\newenvironment{examplecontinued}[2][continued]
  {%
    \renewcommand{\continuedexamplename}{Example~\ref{#2}}%
    \begin{continuedexampleinner}[#1]%
  }
  {\end{continuedexampleinner}}

\theoremstyle{remark}
\newtheorem{remark}[theorem]{Remark}

\numberwithin{equation}{section}

\title[Rate-dependent gradient flows of integral currents]{Rate-dependent gradient flows of\\integral currents via varicurrents}

\date{\today}

\newcommand{\orient}{\mathcal{O}}
\newcommand{\orientSpatial}{\mathcal{O}_{\operatorname{sp}}}
\newcommand{\intcur}[1]{\mathrm I_{#1}}
\newcommand{\varicur}[1]{\mathrm {VI}_{#1}}
\newcommand{\young}{\mathcal Y}

\newcommand{\mass}[1]{\mathbf{M}\mathopen{}\left(#1\right)\mathclose{}}
\newcommand{\variation}[1]{\operatorname{Var}\mathopen{}\left(#1\right)\mathclose{}}
\newcommand{\AC}{\mathrm{AC}}

\newcommand{\eym}[1]{\delta\mathopen{}\left[#1\right]\mathclose{}}
\newcommand{\stcurr}[1]{S\mathopen{}\left[#1\right]\mathclose{}}
\newcommand{\vstcurr}[1]{\overrightarrow{S\mathopen{}\left[#1\right]\mathclose{}}}
\newcommand{\timepush}[2]{{#1}_*{#2}}
\newcommand{\loc}{\operatorname{loc}}
\newcommand{\slice}[2]{#1\mathopen{}\left(#2\right)\mathclose{}}
\newcommand{\slicenov}[2]{\widehat{#1}\mathopen{}\left(#2\right)\mathclose{}}
\newcommand{\LinftySlices}{L^\infty_{\mathrm{sl}}}
\newcommand{\weakstarto}{\overset\ast\rightharpoonup}
\newcommand{\interior}[1]{\operatorname{int}#1}
\newcommand{\statespace}{X_\gamma}

\DeclareMathOperator{\Wedge}{{\textstyle\bigwedge}}

\newcommand{\funcarg}[1]{\ifthenelse{\isempty{#1}}{}{\mathopen{}\left({#1} \right)\mathclose{}}}
\newcommand{\dualR}[1]{\mathcal P^+\funcarg{#1}}
\newcommand{\dualRinf}[1]{\mathcal P^-\funcarg{#1}}

\newcommand{\xproj}[1]{\mathbf p #1}
\newcommand{\tproj}[1]{\mathbf t #1}
\newcommand{\norm}[1]{\mathopen{}\left\| #1\right\|\mathclose{}}
\newcommand{\abs}[1]{\left| #1\right|}
\renewcommand{\d}{\,\mathrm{d}}
\newcommand{\mres}{\mathrel{\tikz[baseline=0ex]{\draw[line width=.5pt,line cap=round] (0ex,1ex) -- (0ex,0ex) -- (1ex,0ex);}}}

\newcommand{\R}{\mathbb R}
\newcommand{\N}{\mathbb N}

\newcommand{\zigzag}{\operatorname{zig}}

\begin{document}

\begin{abstract}
Motivated by the analysis of recent models of dislocation motion in crystalline materials, this work introduces space-time varicurrents, a form of Young measures over integral currents, to describe their evolutions with rate-dependent dissipation potentials. This addresses the problem that rate-dependent versions of the (rate-independent) space-time variation suffer from a lack of lower semicontinuity with respect to the natural weak topology. This causes a dynamical relaxation to occur by means of infinitely fine oscillations. Varicurrents retain a description of this ``microstructure'' (in space-time) and, based on this theory, we are able to establish the existence of solutions to a rate-dependent gradient flow. This is in contrast to the classical formulation using integral currents, which does not have a solution in general.
\end{abstract}

\author{Anton\'{i}n \v{C}e\v{s}\'{i}k}
\address{Mathematics Institute, University of Warwick, Coventry, CV4 7AL, United Kingdom}
\email{Antonin.Cesik@warwick.ac.uk}
	
\author{Malte Kampschulte}
\address{Faculty of Mathematics and Physics, Charles University, Sokolovsk\'{a} 83, 186~75, Prague, Czech Republic}
\email{kampschulte@karlin.mff.cuni.cz} 
	
\author{Filip Rindler}
\address{Mathematics Institute, University of Warwick, Coventry, CV4 7AL, United Kingdom}
\email{F.Rindler@warwick.ac.uk}

\maketitle

\setcounter{tocdepth}{2}

\section{Introduction}

In crystalline materials, plastic deformation is predominantly driven by the motion of dislocation lines, which are one-dimensional topological defects within the crystal lattice~\cite{HullBacon11book,AndersonHirthLothe17book}. Recently, the work \cite{hudsonElastoplasticEvolutionSingle2022} formalized this within a large-strain mesoscopic model for dislocation-driven plasticity, starting from first principles. A central idea in that work is to geometrically express dislocation motion through so-called ``slip trajectories'', that is, the trajectory of a dislocation line traced out in space-time. To express this in a mathematically rigorous way, \cite{rindlerSpacetimeIntegralCurrents2023} introduced the theory of space-time integral currents (we remark that some of these ideas were also developed independently in \cite{Kampschulte17PhD}). Using this framework, it is possible to establish the existence of energetic solutions for elasto-plastic evolutions driven by (discrete) dislocation motion in the rate-independent regime, as carried out in \cite{rindlerEnergeticSolutionsRateindependent2024}.

The works \cite{rindlerSpacetimeIntegralCurrents2023,rindlerEnergeticSolutionsRateindependent2024} are limited to the \emph{rate-independent regime}, where the dissipation of dislocation motion does not depend on the speed of movement (equivalently, the dissipation potential is $1$-homogeneous). However, dislocation motion is usually \emph{rate-dependent}~\cite{BulatovCai06book}, so the work \cite{rindlerEnergeticSolutionsRateindependent2024} must be seen as a slow loading limit of the underlying physical model. 

The present work tackles the rate-dependent case, governed by a $p$-superlinear dissipation potential ($p>1$). In fact, we treat the general theory of moving integral $k$-currents, with dislocation motion corresponding to $k=1$. As such, applications go beyond the original motivation of elasto-plasticity theory and extend also to, for example, membrane transport~\cite{bonicattoTransportCurrentsGeometric2025} and magnetohydrodynamics in the space-time framework, cf.~\cite{bonicatto2025generalfrobeniustheoremtransport}.

A fundamental issue in the rate-dependent regime is that the dissipation potential lacks weak lower semicontinuity in the topology of space-time currents. This geometric degeneracy echoes the phenomenon observed already in \cite{michorVanishingGeodesicDistance2005, michorRiemannianGeometriesSpaces2006} regarding the geometry of spaces of plane curves, and also the examples in \cite{Kampschulte17PhD}. Specifically, a sequence of currents can undergo microscopic oscillations or cancellations to distribute its movement over a larger volume at a lower velocity, thereby artificially decreasing its superlinear cost (as we will illustrate in Examples \ref{ex:zigzag-current} and \ref{ex:thousand-cuts-current}, respectively). When passing to the weak limit of increasingly fast oscillations, the limit current straightens and so ``forgets'' the microstructure. Consequently, the dissipation evaluated on the limit curve jumps up relative to the limit of the dissipation of the approximating motions, violating the standard requirement of weak lower semicontinuity.

Consequently, we need to employ some kind of relaxation. To this end, we introduce the framework of \emph{space-time varicurrents}. Inspired by the theory of Young measures, this notion retains the microscopic orientational structure of the approximating sequence. Moreover, the associated varicurrent topology naturally preserves the lower semicontinuity of the dissipation potential. 

With the necessary lower semicontinuity secured, we formulate a solution concept for a prototypical setup of rate-dependent dislocation evolutions using the so-called \emph{Energy-Dissipation Inequality} (EDI). This is a robust framework extensively developed for metric gradient flows and doubly-nonlinear evolutions \cite{ambrosioGradientFlowsMetric2008}, notably by Mielke, Rossi, and Savar\'e \cite{mielkeVariationalConvergenceGradient2012}. To construct these solutions, we employ the Weighted Energy-Dissipation (WED) principle, introduced for gradient flows by Mielke and Stefanelli \cite{mielkeWeightedEnergydissipationFunctionals2011}. It should be noted that a similar approach had already been used in the context of mean curvature by Ilmanen \cite{ilmanenEllipticRegularizationPartial1994}. Our specific approach is heavily inspired by the extension of the WED principle to metric spaces by Rossi, Savar\'e, Segatti, and Stefanelli \cite{rossiWeightedEnergyDissipationPrinciple2019}. As our space-time varicurrents are inherently global-in-time, they are natural objects for the global-in-time nature of the WED scheme. This alignment allows us to pass to the limit in a sequence of global variational problems, establishing the energy-dissipation inequality and thus obtaining EDI solutions.

Our main contributions are the following:
\begin{itemize}
    \item \textbf{Space-time varicurrents.}
    In Section \ref{sec:varicurrents} we introduce space-time varicurrents, which retain local
    distributions of space-time orientations while having integral
    currents as their barycenters. We develop their slicing theory,
    including slices carrying velocity information, and their
    transformation under changes of time.

    \item \textbf{Compactness and lower semicontinuity.}
    Next, the first part of Section \ref{sec:WED} proves weak${}^*$ lower semicontinuity of the integrated
    dissipation potential and establishes compactness for admissible
    trajectories under uniform bounds.
    These bounds also yield locally uniform convergence of the
    associated current-valued state curves.

    \item \textbf{Existence of rate-dependent evolutions.}
    Then, in the remainder of Section \ref{sec:WED}, for energies satisfying our standing assumptions, we use the 
    weighted energy-dissipation principle to construct global
    evolutions from every admissible initial current. These
    evolutions satisfy the aforementioned Energy-Dissipation Inequality; see Theorem~\ref{thm:energy-dissipation-inequality}.
    
    \item \textbf{Conditional energy balance.}
    Finally, in Section \ref{sec:balance}, we also give conditions under which this inequality becomes
    an energy-dissipation balance
    (Corollary~\ref{cor:energy-dissipation-balance}).
\end{itemize}

\subsection{Examples of relaxation}\label{ssec:examples-currents}

The following two examples illustrate the loss of lower semicontinuity for the superlinear dissipation potential through two different mechanisms: fine oscillations of a moving line (Example \ref{ex:zigzag-current}), and cancellation between lines with opposite orientations (Example \ref{ex:thousand-cuts-current}). We describe the geometry, limiting currents, and dissipational costs here; the complete calculations for both examples are deferred to the Appendix.

\begin{definition}[Zig-zagging function]
We define the zig-zagging function $\zigzag\colon\R\to \R$ by
$$
\zigzag(t):=\begin{cases}
    2t, &t\in[0,\frac12),\\
    2-2t, &t\in[\frac12,1) ,\\
    \zigzag(t-\lfloor t\rfloor), &t\notin [0,1),
\end{cases}
$$
and as a shorthand, a zig-zagging function with period $P>0$ as
$$
\zigzag_{P}(s)=\zigzag\left(\frac{s}{P}\right),
$$
as depicted in Figure \ref{fig:zig-definition}.
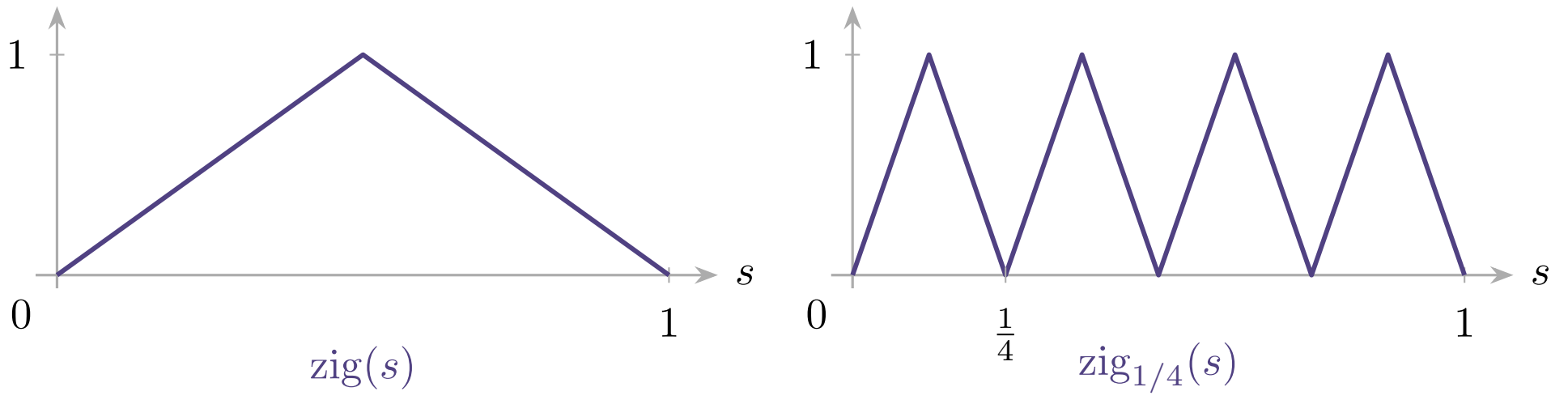
\begin{figure}[htbp]
    \centering
\begin{tikzpicture}[
  x=5cm,y=1.8cm,font=\small,
  >={Stealth[length=1.9mm]},
  graph/.style={draw=slice,line width=1.1pt,line join=round},
  axis/.style={draw=gray!65,line width=.6pt,->}
]
\foreach \N/\shift in {1/0,4/6.5}{
  \begin{scope}[xshift=\shift cm]
    \draw[axis] (-.035,0) -- (1.08,0) node[right,black] {$s$};
    \draw[axis] (0,-.06) -- (0,1.22);
    \draw[gray!65] (-.012,1) -- (.012,1);
    \node[left=3pt] at (0,1) {$1$};
    \node[below left=2pt] at (0,0) {$0$};
    \draw[gray!65] (1,-.035) -- (1,.035);
    \node[below=4pt] at (1,0) {$1$};
    \draw[graph] (0,0)
      \foreach \j in {1,...,\N}{
        -- ({(\j-.5)/\N},1) -- ({\j/\N},0)
      };
  \end{scope}
}
\node[text=slice] at (.5,-.42) {$\zigzag(s)$};
\begin{scope}[xshift=6.5cm]
  \draw[gray!65] (.25,-.035) -- (.25,.035);
  \node[below=4pt] at (.25,0) {$\tfrac14$};
  \node[text=slice] at (.5,-.42) {$\zigzag_{1/4}(s)$};
\end{scope}
\end{tikzpicture}

    \caption{The zigzagging function $\zigzag$ resp.\ $\zigzag_P$.}\label{fig:zig-definition}
\end{figure}
\end{definition}

The following two examples demonstrate the failure of lower semicontinuity of the dissipation potential, when treated in the sense of currents. Later in Section \ref{ssec:examples-revisited}, after developing the theory of varicurrents, we will revisit these examples to see that when understood as varicurrents, the (generalized) dissipation potential behaves correctly. As we will detail in the section on notation, we will fix $1<p<\infty$ throughout.

Given an integral current $T\in\intcur{k}(\overline\Omega)$ and a velocity field $v\colon\operatorname{supp}T\to\R^d$ (defined $\mathcal H^k$-a.e.), we define the \emph{dissipation potential} to be
$$
\mathcal R_p(T,v) := \frac1p\int|v(x)|^p\d\|T\|(x).
$$
For a space-time trajectory $S\in\intcur{1+k}(I\times\overline\Omega)$ admitting a Lebesgue-time disintegration, the geometric derivative $\frac{D}{Dt}S$ supplies the normal velocity of its slices whenever it is defined.\footnote{Precisely, if $\tproj{}_\#\|S\|\ll\mathcal L^1$ and $S$ satisfies \eqref{eqn:NC-current}; see Section~\ref{sec:preliminaries} for the time-disintegration and geometric derivative.} Under these conditions, the integrated dissipation potential over a set of times $A\subset I$ is
\begin{equation}\label{eqn:intro-Rp-current-defn}
\int_A \mathcal R_p\left(S(t), \frac{D}{Dt}S\right) \d t= \frac1p\int_A \int\abs{\frac{D}{Dt} S(t,x)}^p \d \|S(t)\|(x)\d t.
\end{equation}
Here, $S(t)$ denotes the curve at time $t$; see the next section for details.

To include the time-concentration limits appearing in the examples below, we extend this definition to any $S\in\intcur{1+k}(I\times\overline\Omega)$ by setting
\begin{equation}\label{eq:example-current-cost}
\Psi_p(S;A):=\int_{A\times\overline\Omega}f_p(\vec S(t,x))\d\|S\|(t,x),
\qquad A\subset \R \text{ Borel},
\end{equation}
where, for any unit simple $(1+k)$-vector $\eta\in\Wedge^{1+k}\R^{1+d}$,
$$
f_p(\eta)=\begin{cases}
    \frac{|\xproj{\eta}|^p}{p|\eta \mres \mathrm dt|^{p-1}}, &\text{if }|\eta \mres \mathrm dt|\neq 0,
    \\
    \infty, & \text{if }|\eta \mres \mathrm dt|= 0.
\end{cases}
$$

For currents satisfying \eqref{eqn:NC-current} and $\tproj{}_\#\|S\|\ll\mathcal L^1$, this agrees with \eqref{eqn:intro-Rp-current-defn}:
$$
\Psi_p(S;A)=\int_A \mathcal R_p\left(S(t), \frac{D}{Dt}S\right) \d t.
$$

Note, however, that $\Psi_p$ is \emph{not} the dissipation along the evolution (unless $p=1$, which is here excluded), but a time-integrated dissipation potential.
We will see in the following two examples that $\Psi_p$ is not lower semicontinuous (with respect to the weak convergence of currents) since $p>1$, and in fact relaxes to zero under such a convergence. This prevents the use of classical gradient flow methods and necessitates the introduction of the notion of varicurrents.

\begin{example}[{Zig-zag, cf.\ \cite[Thm 3.2]{michorVanishingGeodesicDistance2005}, \cite[Ex 9.1]{Kampschulte17PhD}}]\label{ex:zigzag-current}
Consider the rectangle $\Omega=(0,1)\times(-1,2)\subset\R^2$, and use the coordinate unit vectors $e_0,e_1,e_2$ in the $t,x,y$ directions, respectively. A zigzag line translates upward by one unit during $[0,\Delta t]$ and is stationary afterwards; see Figure~\ref{fig:zigzag-slice}.

Let $0<\Delta t\leq1$, $\alpha\in(0,\pi/2)$ and $N\in\N$. The space-time current $S_{\Delta t,\alpha,N}\in\intcur{2}([0,1]\times\overline\Omega)$ has slices with $N$ zigs and zags of angle $\alpha$. During the moving interval, its slice at time $t$ is the graph of
\begin{equation}\label{eq:zigzag-slice-graph}
y_t(x):=\frac{t}{\Delta t}+h_{\alpha,N}\zigzag_{1/N}(x), \quad t\in[0,\Delta  t] \quad \text{with }h_{\alpha,N}:=\frac{\cot\alpha}{2N}.
\end{equation}
We assume $h_{\alpha,N}\leq1$ so that the trajectory stays in $[0,1]\times\overline\Omega$. For $t\in[\Delta t,1]$, the slice is the final translated zigzag graph $y_{\Delta t}$, and we set $y_t:=y_{\Delta t}$ on this interval.
\begin{figure}[htbp]
\centering
\begin{tikzpicture}[
  x=1cm,y=1cm, font=\small,
  >={Stealth[length=2.2mm]},
  normal/.style={->,velocity,line width=1.1pt},
  tangential/.style={->,tangent,line width=1.1pt},
  construction/.style={gray!65,densely dashed,line width=.55pt}
]
\def\angledeg{60}
\pgfmathsetmacro{\height}{1.5*cot(\angledeg)}
\def\U{1.75} 
\def\T{0.85} 

\coordinate (O) at (0,0);
\coordinate (C) at (4.5,\height);
\coordinate (Utop) at ($(C)+(0,\U)$);
\coordinate (D) at (4.5,0);
\coordinate (L) at (3,0);
\coordinate (Pplus) at (3.75,{\height/2});
\coordinate (Pminus) at (5.25,{\height/2});
\coordinate (Uplus) at ($(Pplus)+(0,\U)$);
\coordinate (Uminus) at ($(Pminus)+(0,\U)$);
\coordinate (Qplus) at
  ($(Pplus)+({-\U*sin(\angledeg)*cos(\angledeg)},{\U*sin(\angledeg)^2})$);
\coordinate (Qminus) at
  ($(Pminus)+({\U*sin(\angledeg)*cos(\angledeg)},{\U*sin(\angledeg)^2})$);

\draw[slice,line width=1.15pt]
  (0,0) -- (1.5,\height) -- (3,0) -- (4.5,\height)
  -- (6,0) -- (7.5,\height) -- (9,0) -- (10.5,\height) -- (12,0);

\draw[->,gray!70] (0,-.65) -- (12.45,-.65) node[right,black] {$x$};
\draw[->,gray!70] (0,-.65) -- (0,2.9) node[above,black] {$y$};
\draw[gray!70] (-.07,0) -- (.07,0) (-.07,\height) -- (.07,\height);
\draw[decorate,decoration={brace,amplitude=3pt},gray!70]
  (-.12,0) -- (-.12,\height)
  node[midway,left=4pt,inner sep=0pt,black] {$h_{\alpha,N}$};
\node[text=slice] at (9,1.75) {$S_{\Delta t,\alpha,N}(t)$};

\draw[construction] (C) -- (D);
\draw[gray!80] (C) ++(210:.47) arc[start angle=210,end angle=270,radius=.47];
\node[font=\scriptsize] at ($(C)+(240:.72)$) {$\alpha$};

\draw[construction] (Qplus) -- (Utop);
\draw[construction] (Qminus) -- (Utop);
\pic[draw=gray!75,angle radius=4pt] {right angle=Pplus--Qplus--Utop};
\pic[draw=gray!75,angle radius=4pt] {right angle=Utop--Qminus--Pminus};

\draw[normal] (Pplus) -- (Qplus) node[left=3pt] {$v_+$};
\draw[normal] (Pminus) -- (Qminus) node[right=3pt] {$v_-$};
\draw[tangential] (Pplus) -- ++({\T*sin(\angledeg)},{\T*cos(\angledeg)})
  node[pos=.65,above left=1pt] {$\tau_+$};
\draw[tangential] (Pminus) -- ++({\T*sin(\angledeg)},{-\T*cos(\angledeg)})
  node[pos=.55,below left=1pt] {$\tau_-$};
\fill[slice] (Pplus) circle (1.3pt);
\fill[slice] (Pminus) circle (1.3pt);

\draw[->,black!80,line width=1pt] (C) -- (Utop);
\draw[decorate,decoration={brace,amplitude=4pt},black!80]
  ($(Utop)+(.18,0)$) -- ($(C)+(.18,0)$)
  node[midway,right=4pt,inner sep=0pt,black] {$\dfrac{1}{\Delta t}$};
\end{tikzpicture}

\caption{The slice $S_{\Delta t,\alpha,N}(t)$ for $\alpha=\pi/3$ and $N=4$ at time $t\in[0,\Delta t]$ with its normal velocity vectors $v_\pm$ and tangent vectors $\tau_\pm$.}\label{fig:zigzag-slice}
\end{figure}
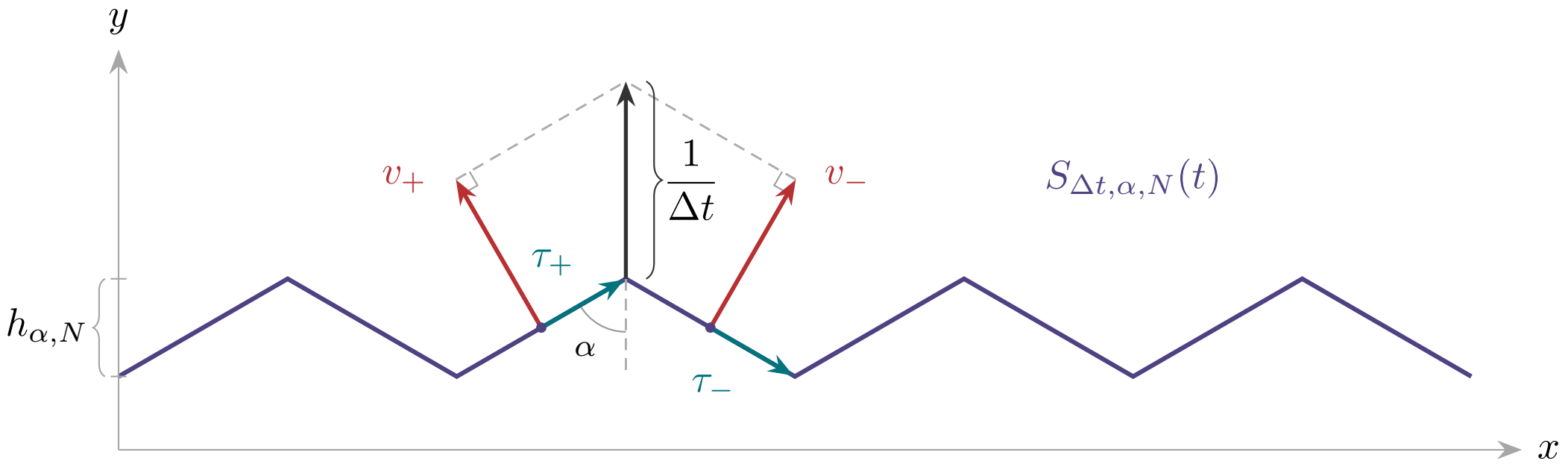

We thus define the space-time current by
$$
S_{\Delta t,\alpha,N}:=S_{\Delta t,\alpha,N}^{\mathrm m}+S_{\Delta t,\alpha,N}^{\mathrm s},
$$
where the respective moving and stationary parts are
\begin{align*}
S_{\Delta t,\alpha,N}^{\mathrm m}
&:=\bigl((t,x)\mapsto\left(t,x,y_t(x)\right)\bigr)_*
\llbracket(0,\Delta t)\times(0,1)\rrbracket,
\\
S_{\Delta t,\alpha,N}^{\mathrm s}
&:=\bigl((t,x)\mapsto\left(t,x,y_{\Delta t}(x)\right)\bigr)_*
\llbracket(\Delta t,1)\times(0,1)\rrbracket.
\end{align*}

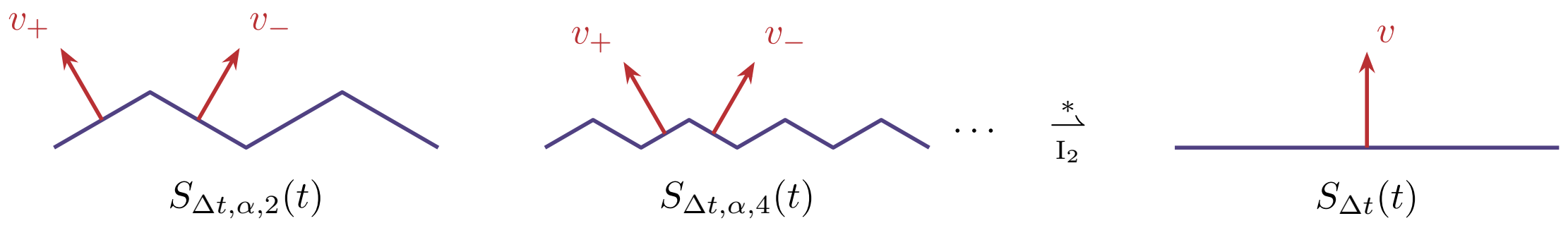
\begin{figure}[htbp]
    \centering
\begin{tikzpicture}[
  x=1cm,y=1cm,font=\small,
  >={Stealth[length=2mm]},
  zigcurve/.style={draw=slice,line width=1.05pt,line join=round},
  zigvelocity/.style={->,draw=velocity,text=velocity,line width=1.05pt}
]
\def\panelwidth{3.6}
\def\angledeg{60}
\def\velocityscale{.9} 

\foreach \N/\offset/\tooth in {2/0/0,4/4.6/1}{
  \begin{scope}[xshift=\offset cm]
    \pgfmathsetmacro{\period}{\panelwidth/\N}
    \pgfmathsetmacro{\toothheight}{\period*cot(\angledeg)/2}
    \draw[zigcurve] (0,0)
      \foreach \j in {1,...,\N}{
        -- ({(\j-.5)*\period},\toothheight)
        -- ({\j*\period},0)
      };
    \coordinate (Pplus) at
      ({(\tooth+.25)*\period},{\toothheight/2});
    \coordinate (Pminus) at
      ({(\tooth+.75)*\period},{\toothheight/2});
    \draw[zigvelocity] (Pplus)
      -- ++({-\velocityscale*sin(\angledeg)*cos(\angledeg)},
            {\velocityscale*sin(\angledeg)^2})
      node[above left=1pt,inner sep=1pt] {$v_+$};
    \draw[zigvelocity] (Pminus)
      -- ++({\velocityscale*sin(\angledeg)*cos(\angledeg)},
            {\velocityscale*sin(\angledeg)^2})
      node[above right=1pt,inner sep=1pt] {$v_-$};
    \node at ({\panelwidth/2},-.48) {$S_{\Delta t,\alpha,\N}(t)$};
  \end{scope}
}

\node at (8.65,.15) {$\cdots$};
\node at (9.5,.15)
  {$\displaystyle\underset{\intcur{2}}{\overset{*}{\rightharpoonup}}$};

\begin{scope}[xshift=10.5cm]
  \draw[zigcurve] (0,0) -- (\panelwidth,0);
  \draw[zigvelocity] ({\panelwidth/2},0)
    -- ++(0,\velocityscale)
    node[above right=1pt,inner sep=1pt] {$v$};
  \node at ({\panelwidth/2},-.48) {$S_{\Delta t}(t)$};
\end{scope}
\end{tikzpicture}

    \caption{Spatial slices of the approximating currents $S_{\Delta t_i,\alpha_i,N_i}$ with $\alpha_i\equiv\pi/3$ and of their weak${}^*$ limit $S_{\Delta t}$.}\label{fig:zigzag-convergence}
\end{figure}

The masses of the slices (i.e.\ the zigzag length) and pointwise normal speed are
$$
\mass{S_{\Delta t,\alpha,N}(t)}=\frac1{\sin\alpha},\quad t\in[0,1],
\qquad |v_\sigma|=\frac{\sin\alpha}{\Delta t},\quad 0<t<\Delta t,
$$
where the sign $\sigma \in\{+1,-1\}$ corresponds to the rising and falling segments respectively.

Computing $\Psi_p$ gives
\begin{equation}\label{eqn:zigzag-finite-cost}
\begin{aligned}
\Psi_p(S_{\Delta t,\alpha,N};[0,1])
&= \frac1p\int_0^{\Delta t}\int|v_\sigma|^p\d\|S_{\Delta t,\alpha,N}(t)\| \d t\\
&=\frac1p\int_0^{\Delta t}\left(\frac{\sin\alpha}{\Delta t}\right)^p\frac1{\sin\alpha}\d t
=\frac1p\left(\frac{\sin\alpha}{\Delta t}\right)^{p-1}.
\end{aligned}
\end{equation}
Take a sequence $\Delta t_i\in(0,1], \alpha_i\in(0,\pi/2), N_i\in\N$,  for $i\in\N$ with
\begin{equation}\label{eqn:limits-zigzag}
\lim_{i\to\infty} N_i=\infty, \quad \lim_{i\to\infty} \alpha_i=\alpha\in[0,\pi/2),\quad \lim_{i\to\infty} \Delta t_i =\Delta t\in[0,1],\quad \lim_{i\to\infty} h_{\alpha_i,N_i}=0.
\end{equation}
Note in particular that we allow the limiting angle $\alpha =0$ and limiting time $\Delta t=0$. As proved in Appendix \ref{sec:zigzag-appendix}, the integral currents converge as illustrated by the spatial slices in Figure~\ref{fig:zigzag-convergence}:
$$
S_{\Delta t_i,\alpha_i, N_i}\weakstarto  S_{\Delta t}=S_{\Delta t}^{\mathrm m}+S_{\Delta t}^{\mathrm s}\quad \text{in }\intcur{2}([0,1]\times\overline\Omega), 
$$
where the moving part is a straight horizontal line traversing from $y=0$ to $y=1$ in time $\Delta t$, and the stationary part stays at $y=1$, namely\footnote{The minus in the $\Delta t=0$ case corresponds to the current having orientation $\mathrm dy 
\wedge \mathrm dx= -\mathrm dx \wedge \mathrm dy$.}
\begin{align*}
S_{\Delta t}^{\mathrm m}&=\begin{cases}
\left((t,x)\mapsto(t,x,t/\Delta t)\right)_*
\llbracket(0,\Delta t)\times(0,1)\rrbracket
&\text{for }\Delta t>0,
\\
-\llbracket\{0\}\times(0,1)\times(0,1)\rrbracket
&\text{for }\Delta t=0,
\end{cases}
\\
S_{\Delta t}^{\mathrm s}&=\llbracket(\Delta t,1)\times(0,1)\times\{1\}\rrbracket.
\end{align*}
If $\Delta t>0$, the cost satisfies
$$
\Psi_p (S_{\Delta t_i,\alpha_i,N_i};[0,1])
= \tfrac1p\left(\frac{\sin\alpha_i}{\Delta t_i}\right)^{p-1} \overset{i\to\infty}{\to} \tfrac1p\left(\frac{\sin\alpha}{\Delta t}\right)^{p-1} <\tfrac1p \left(\frac{1}{\Delta t}\right)^{p-1} =\Psi_p(S_{\Delta t};[0,1]).
$$
In particular, the inequality is strict for every limiting angle $0\leq\alpha<\pi/2$.

If $\Delta t=0$, the limiting motion takes place instantaneously. Here we assume in addition that the time and angle vanish at a comparable rate:
\begin{equation}\label{eqn:beta-zigzag}
\lim_{i\to\infty}\frac{\Delta t_i}{\sin \alpha_i}= \beta\in(0,\infty).
\end{equation}
Then, we have
$$
\Psi_p(S_{\Delta t_i,\alpha_i,N_i};[0,1])
= \tfrac1p\left(\frac{\sin\alpha_i}{\Delta t_i}\right)^{p-1} \overset{i\to\infty}{\to} \tfrac1p \beta^{1-p} <\infty=\Psi_p(S_{0};[0,1]).
$$
\end{example}

\begin{example}[Thousand Cuts]\label{ex:thousand-cuts-current}
We use the coordinate unit vectors $e_0,e_1,e_2$ in the $t,x,y$ directions, respectively. 

In the same domain $\Omega=(0,1)\times(-1,2)$, a horizontal line starts at $y=0$ and finishes translated upwards by one unit, similarly as above, but this movement proceeds with $N$ intermediate pairs with opposite orientations. Let $0<\Delta t\leq1$ and $N\in\N$. During $(0,\Delta t)$, $N$ intermediate pairs are created, so there are $N+1$ upward-moving lines with spatial orientation $e_1$ and $N$ downward-moving lines with the opposite orientation $-e_1$. At time $\Delta t$, all pairs cancel and only the top line at $y=1$ remains.

Writing
$$
y_j^\pm(t):=\frac{2j\pm t/\Delta t}{2N+1},
$$
with $j=0,\ldots,N$ for the upward lines $y^+_j$, and $j=1,\ldots,N$ for the downward lines $y^-_j$, the full space-time current $S_{\Delta t,N}\in\intcur{2}([0,1]\times\overline\Omega)$ has the slices (see Figure \ref{fig:thousand-cuts})
$$
S_{\Delta t,N}(t)=\begin{cases}
\llbracket (0,1) \times \{0\} \rrbracket, &t=0,
\\
\displaystyle\sum_{j=0}^N\llbracket(0,1)\times\{y_j^+(t)\}\rrbracket
-\displaystyle\sum_{j=1}^N\llbracket(0,1)\times\{y_j^-(t)\}\rrbracket,
&0<t<\Delta t,
\\
\llbracket (0,1) \times \{1\} \rrbracket, &\Delta t\leq t\leq1.
\end{cases}
$$
We use the superscripts $\mathrm m$ and $\mathrm s$ for the moving and stationary parts, respectively, and write
\begin{align*}
S_{\Delta t,N}
&=S_{\Delta t,N}^{\mathrm m}+S_{\Delta t,N}^{\mathrm s},\\
S_{\Delta t,N}^{\mathrm m}
&:=S_{\Delta t,N}\mres((0,\Delta t)\times\overline\Omega),\\
S_{\Delta t,N}^{\mathrm s}
&:=\llbracket(\Delta t,1)\times(0,1)\times\{1\}\rrbracket.
\end{align*}
Here $S_{\Delta t,N}^{\mathrm s}$ has orientation $e_0\wedge e_1$ and zero velocity. The moving slices satisfy
$$
\mass{S_{\Delta t,N}(t)}=2N+1,\qquad
v^\pm=\pm\frac{e_2}{(2N+1)\Delta t},\qquad
|v^\pm|=\frac1{(2N+1)\Delta t},
\quad 0<t<\Delta t.
$$
Consequently, the extended dissipation potential \eqref{eq:example-current-cost} is
\begin{equation}\label{eqn:cuts-finite-cost}
\Psi_p(S_{\Delta t,N};[0,1])
=\Psi_p(S_{\Delta t,N}^{\mathrm m};[0,1]) =\frac1p\int_0^{\Delta t}\frac{2N+1}{((2N+1)\Delta t)^p}\d t
=\frac1p\bigl((2N+1)\Delta t\bigr)^{1-p}.
\end{equation}

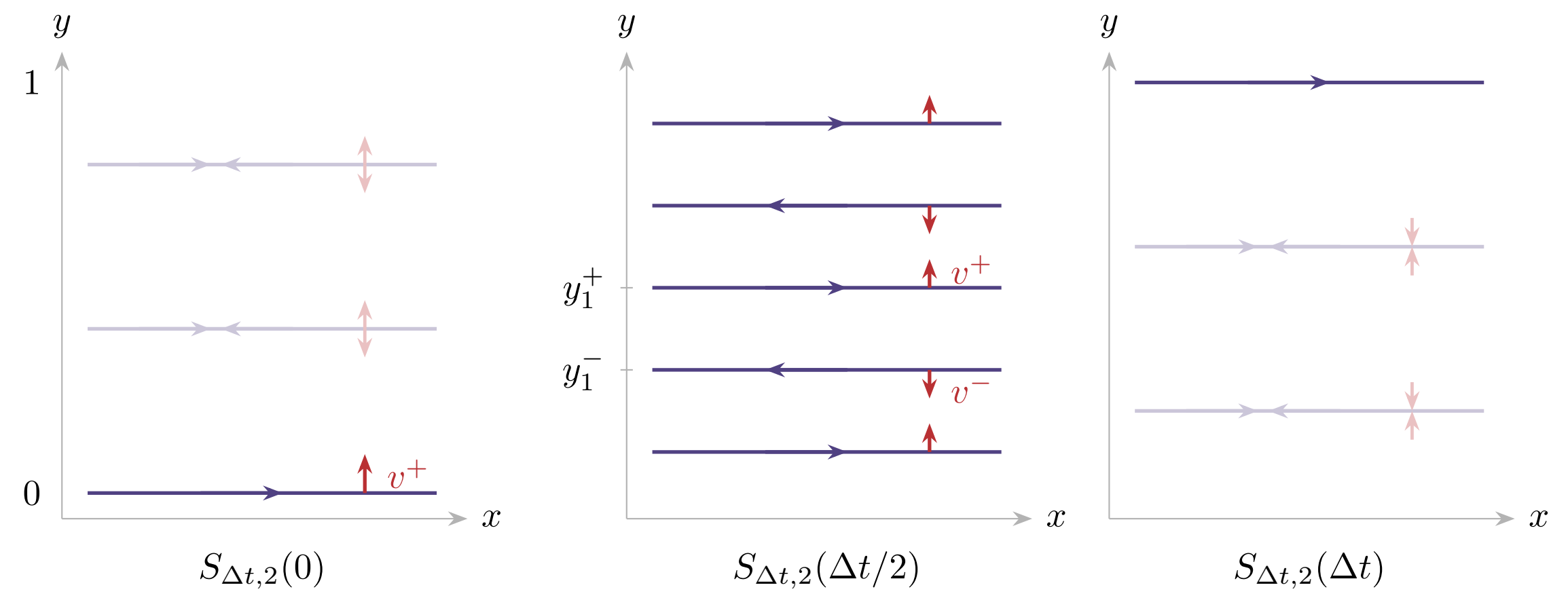
\begin{figure}
    \centering
\begin{tikzpicture}[
  x=1cm,y=1cm,font=\small,>={Stealth[length=1.9mm]},
  cut/.style={draw=slice,line width=1pt},
  ghost/.style={draw=slice!30,line width=1pt},
  normal/.style={->,draw=velocity,text=velocity,line width=1pt},
  ghostnormal/.style={->,draw=velocity!30,line width=.9pt}
]
\def\H{4}
\def\W{3.4}
\def\V{.28}

\begin{scope}
  \draw[cut] (0,0) -- (\W,0);
  \draw[cut,->] (1.1,0) -- (1.9,0);
  \draw[normal] (2.7,0) -- ++(0,{\V+0.1}) node[midway,right=2pt] {$v^+$};
  \foreach \j in {1,2}{
    \pgfmathsetmacro{\Y}{\H*2*\j/5}
    \draw[ghost] (0,\Y) -- (\W,\Y);
    \draw[ghost,->] (.5,\Y) -- (1.2,\Y);
    \draw[ghost,->] (2,\Y) -- (1.3,\Y);
    \draw[ghostnormal] (2.7,\Y) -- ++(0,\V);
    \draw[ghostnormal] (2.7,\Y) -- ++(0,-\V);
  }
  \node at (1.7,-.75) {$S_{\Delta t,2}(0)$};
\end{scope}

\begin{scope}[xshift=5.5cm]
  \foreach \j in {0,1,2}{
    \pgfmathsetmacro{\Y}{\H*(2*\j+.5)/5}
    \draw[cut] (0,\Y) -- (\W,\Y);
    \draw[cut,->] (1.1,\Y) -- (1.9,\Y);
    \draw[normal] (2.7,\Y) -- ++(0,\V);
  }
  \foreach \j in {1,2}{
    \pgfmathsetmacro{\Y}{\H*(2*\j-.5)/5}
    \draw[cut] (0,\Y) -- (\W,\Y);
    \draw[cut,->] (1.9,\Y) -- (1.1,\Y);
    \draw[normal] (2.7,\Y) -- ++(0,-\V);
  }
  \node[text=velocity,right=2pt] at (2.7,2.18) {$v^+$};
  \node[text=velocity,right=2pt] at (2.7,{\H*(2-.5)/5-\V+0.1}) {$v^-$};
  \draw[gray!60] (-.31,{\H*(2+.5)/5}) -- (-.19,{\H*(2+.5)/5});
  \node[left=2pt] at (-.25,{\H*(2+.5)/5}) {$y_1^+$};
  \draw[gray!60] (-.31,{\H*(2-.5)/5}) -- (-.19,{\H*(2-.5)/5});
  \node[left=2pt] at (-.25,{\H*(2-.5)/5}) {$y_1^-$};
  \node at (1.7,-.75) {$S_{\Delta t,2}(\Delta t/2)$};
\end{scope}

\begin{scope}[xshift=10.2cm]
  \foreach \j in {1,2}{
    \pgfmathsetmacro{\Y}{\H*(2*\j-1)/5}
    \draw[ghost] (0,\Y) -- (\W,\Y);
    \draw[ghost,->] (.5,\Y) -- (1.2,\Y);
    \draw[ghost,->] (2,\Y) -- (1.3,\Y);
    \draw[ghostnormal] (2.7,{\Y-\V}) -- (2.7,\Y);
    \draw[ghostnormal] (2.7,{\Y+\V}) -- (2.7,\Y);
  }
  \draw[cut] (0,\H) -- (\W,\H);
  \draw[cut,->] (1.1,\H) -- (1.9,\H);
  \node at (1.7,-.75) {$S_{\Delta t,2}(\Delta t)$};
\end{scope}

\foreach \offset in {0,5.5,10.2}{
  \begin{scope}[xshift=\offset cm]
    \draw[->,gray!60] (-.25,-.25) -- (3.7,-.25) node[right,black] {$x$};
    \draw[->,gray!60] (-.25,-.25) -- (-.25,4.3) node[above,black] {$y$};
  \end{scope}
}
\node[left=2pt] at (-.25,0) {$0$};
\node[left=2pt] at (-.25,\H) {$1$};
\end{tikzpicture}

    \caption{Slices of the Thousand Cuts for $N=2$ at times $t=0,\Delta t/2,\Delta t$, and their corresponding velocity vectors $v^\pm$.}
    \label{fig:thousand-cuts}
\end{figure}
Now take a sequence $N_i\in\N$, $\Delta t_i\in(0,1]$ such that
\begin{equation}\label{eqn:limits-cuts}
 \lim_{i\to\infty}N_i= \infty,
 \quad
 \lim_{i\to\infty}\Delta t_i=0,
 \quad
 \lim_{i\to\infty}(2N_i+1)\Delta t_i=\beta\in[0,\infty).
\end{equation}
We denote
\begin{equation}\label{eqn:ell-cuts}
\ell_{\Delta t,N}:=\sqrt{1+(2N+1)^2\Delta t^2},
\qquad
\ell_\beta:=\sqrt{1+\beta^2}.
\end{equation}

Thus $\ell_{\Delta t_i,N_i}\to\ell_\beta$.
As space-time currents, the moving parts converge to
$$
S_{\Delta t_i, N_i}^{\mathrm m}\weakstarto S_0^{\mathrm m}:=\llbracket\{0\}\times(0,1)\times(0,1)\rrbracket
\quad\text{with orientation }e_2\wedge e_1.
$$
The stationary tails converge to $S_0^{\mathrm s}=\llbracket(0,1)\times(0,1)\times\{1\}\rrbracket$ with orientation $e_0\wedge e_1$, so the full-current limit is
$$
S_{\Delta t_i,N_i}\weakstarto S_0=S_0^{\mathrm m}+S_0^{\mathrm s}
\quad\text{in }\intcur{2}([0,1]\times\overline\Omega).
$$
In the $(t,y)$-plane, the moving part forms a zigzag, see Figure~\ref{fig:thousand-cuts-space-time}. The calculations in the Appendix reveal this explicitly with parametrization of $t$ by $y$ via the zigzag function.

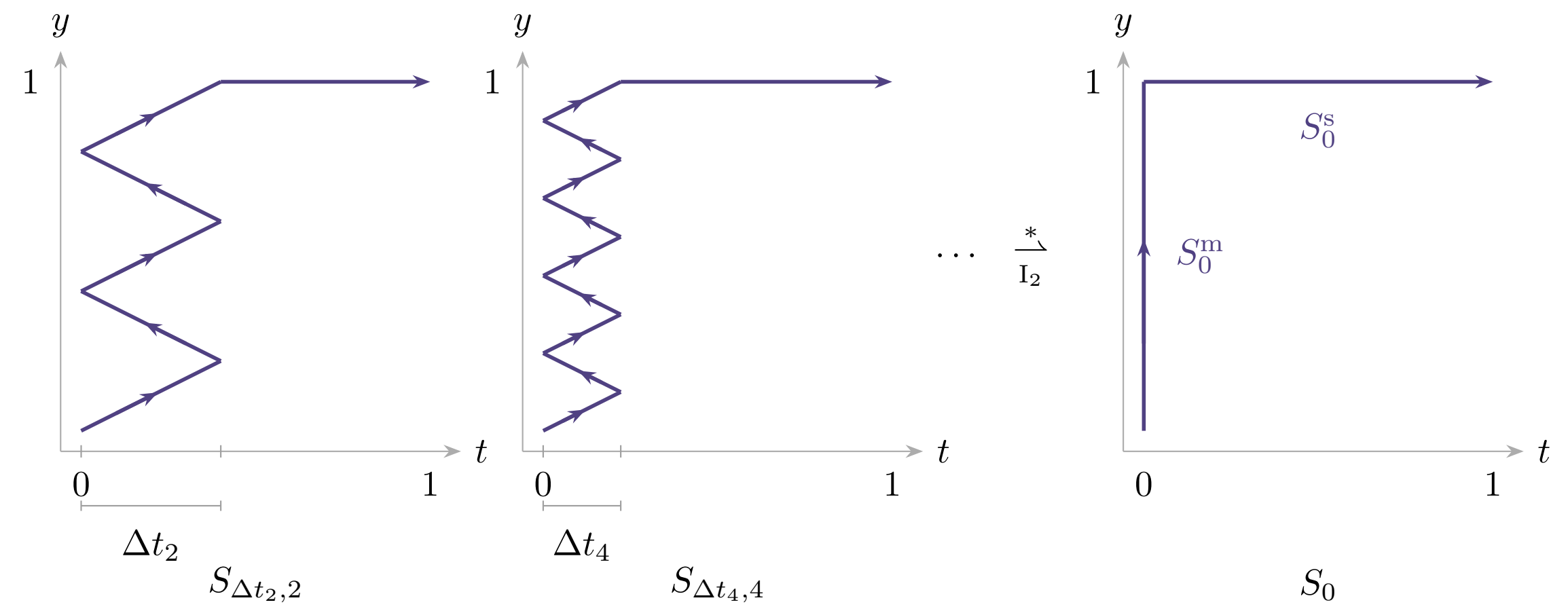
\begin{figure}
    \centering
\begin{tikzpicture}[
  x=1cm,y=1cm,font=\small,>={Stealth[length=1.7mm]},
  moving/.style={draw=slice,line width=1.05pt},
  tail/.style={draw=slice,line width=1.05pt,->}
]
\def\W{3.4}
\def\H{3.4}
\def\betaval{2}
\foreach \N/\offset in {2/0,4/4.5}{
  \begin{scope}[xshift=\offset cm]
    \pgfmathsetmacro{\D}{\W*\betaval/(2*\N+1)}
    \foreach \j in {0,...,\N}{
      \coordinate (A) at (0,{\H*2*\j/(2*\N+1)});
      \coordinate (B) at (\D,{\H*(2*\j+1)/(2*\N+1)});
      \draw[moving] (A) -- (B);
      \draw[moving,->] ($(A)!.25!(B)$) -- ($(A)!.55!(B)$);
    }
    \foreach \j in {1,...,\N}{
      \coordinate (A) at (\D,{\H*(2*\j-1)/(2*\N+1)});
      \coordinate (B) at (0,{\H*2*\j/(2*\N+1)});
      \draw[moving] (A) -- (B);
      \draw[moving,->] ($(A)!.25!(B)$) -- ($(A)!.55!(B)$);
    }
    \draw[tail] (\D,\H) -- (\W,\H);
    \draw[->,gray!65] (-.2,-.2) -- (3.7,-.2) node[right,black] {$t$};
    \draw[->,gray!65] (-.2,-.2) -- (-.2,3.7) node[above,black] {$y$};
    \node[left=2pt] at (-.2,\H) {$1$};
    \node[below=2pt] at (\W,-.2) {$1$};
    \draw[gray!75] (0,-.14) -- (0,-.26);
    \draw[gray!75] (\D,-.14) -- (\D,-.26);
    \node[below=2pt] at (0,-.2) {$0$};
    \draw[gray!75] (0,-.73) -- (\D,-.73);
    \draw[gray!75] (0,-.68) -- (0,-.78);
    \draw[gray!75] (\D,-.68) -- (\D,-.78);
    \node[below=3pt] at ({\D/2},-.73) {$\Delta t_{\N}$};
    \node at (1.7,-1.50) {$S_{\Delta t_{\N},\N}$};
  \end{scope}
}
\node at (8.55,1.7) {$\cdots$};
\node at (9.25,1.7)
  {$\displaystyle\underset{\intcur{2}}{\overset{*}{\rightharpoonup}}$};
\begin{scope}[xshift=10.35cm]
  \draw[moving] (0,0) -- (0,\H);
  \draw[moving,->] (0,{.25*\H}) -- (0,{.55*\H});
  \draw[tail] (0,\H) -- (\W,\H);
  \draw[->,gray!65] (-.2,-.2) -- (3.7,-.2) node[right,black] {$t$};
  \draw[->,gray!65] (-.2,-.2) -- (-.2,3.7) node[above,black] {$y$};
  \node[left=2pt] at (-.2,\H) {$1$};
  \node[below=2pt] at (\W,-.2) {$1$};
  \node[below=2pt] at (0,-.2) {$0$};
  \node[text=slice,right=5pt] at (0,1.7) {$S_0^{\mathrm m}$};
  \node[text=slice,below=5pt] at (1.7,\H) {$S_0^{\mathrm s}$};
  \node at (1.7,-1.50) {$S_0$};
\end{scope}
\end{tikzpicture}

    \caption{Thousand Cuts currents $S_{\Delta t_i,N_i}$ and their weak${}^*$ limit $S_0$, viewed in the $(t,y)$-plane. The $x$-direction is suppressed.}
    \label{fig:thousand-cuts-space-time}
\end{figure}

These convergences are proved in Appendix \ref{sec:thousand-cuts-appendix}. We have $\Psi_p(S_0;[0,1])=\infty$. Since $(2N_i+1)\Delta t_i\to\beta$, for $\beta>0$ we obtain
$$
\lim_{i\to\infty}\Psi_p(S_{\Delta t_i,N_i};[0,1])
=\frac1p\beta^{1-p}<\infty=\Psi_p(S_0;[0,1]).
$$
\end{example}

In the finite-$\Psi_p$ regimes of the Zig-zag and Thousand Cuts examples, the functional $\Psi_p$ jumps strictly upward in the limit when evaluated on space-time currents. In the Zig-zag example this holds for $\Delta t>0$ and for $\Delta t=0$ with $\beta>0$; in the Thousand Cuts example it holds for $\beta>0$. The weak current limit loses the microscopic orientations and their associated normal velocities. This loss of lower semicontinuity prevents the use of variational methods. To retain this information, we introduce varicurrents in the following subsections. We will return to these examples in Section \ref{ssec:examples-revisited} to demonstrate how varicurrents resolve the issue here.

\subsection*{Acknowledgments}
The authors would like to thank Thomas Hudson and Harry Turnbull for discussions related to this work. This project has received funding from the UKRI Frontier Research Guarantee (ERC guarantee) grant EP/Z000297/1 (ERC CONCENTRATE). M.K.\ was supported by the Ministry of Education, Youth and Sport of the Czech Republic through the project Ferroic Multifunctionalities (FerrMion) [reg.\ no.\ CZ.02.01.01/00/22\_008/0004591] as well as by the Charles University through the grant PRIMUS/24/SCI/020 and the Research Centre program No.\ UNCE/24/SCI/005.  

\subsection*{AI disclosure}
AI was only used for proofreading and text editing purposes to streamline the presentation. This includes simplifying technical steps in the proofs, checking consistency of used notation, and drawing the figures. The mathematical content has been developed by the authors without AI assistance.

\section{Space-time varicurrents} \label{sec:varicurrents}

\subsection{Notation and preliminaries}\label{sec:preliminaries}

Before introducing varicurrents, we recall some notions from measure theory, as well as (space-time) integral currents. For more details, we refer the reader to \cite{rindlerEnergeticSolutionsRateindependent2024}, \cite{bonicattoTransportCurrentsGeometric2025} and references therein.

Throughout this paper, $\Omega\subset\R^d$ is a bounded Lipschitz domain. The time interval $I\subset \R$ is closed, either bounded or $I=[0,\infty)$. The exponent $1<p<\infty$ is fixed throughout, and $p'=p/(p-1)$ is the dual exponent.

\subsubsection*{Measures.}
We write $\mathcal L^1$ for one-dimensional Lebesgue
measure and $\mathcal M_{\loc}^+(X)$ for positive
Radon measures on $X$ that are finite on compact sets.
Weak${}^*$ convergence in $\mathcal M_{\loc}^+(X)$
is understood against test functions in $C_c(X)$.
We denote the pushforward measure by $\varphi_\#\mu$. Restriction of a measure $\mu$ to a set $A$ is denoted $\mu\mres A$.
The notation
$\mu=\kappa(\mathrm dy)\otimes\mu^y(\mathrm dz)$
denotes disintegration with respect to $\kappa$;
for Lebesgue measure in time, we write
$\mu=\mathrm dt\otimes\mu^t$.

\subsubsection*{Integral currents}
We denote by $\intcur{k}(\overline\Omega)$ integral
$k$-currents in $\R^d$ supported in $\overline\Omega$. General $k$-currents in $\overline\Omega$ are denoted by $\mathcal D_k(\overline\Omega)$. We use analogous notation in other dimensions and domains.
The notation $\intcur{1+k}^{\loc}(I\times\overline\Omega)$
requires locally finite mass and boundary mass; this is only relevant in the case $I=[0,\infty)$, as $\overline\Omega$ is compact.

The symbol $\Wedge^k\R^d$ denotes the $k$-th exterior power of $\R^d$; its elements are referred to as $k$-vectors. For an integral current $T$, we denote by $\vec T$ its orienting multivector and $\|T\|$ its mass measure, so that
$$T= \vec T \, \|T\|,\qquad \text{i.e.}\quad\langle T,\omega\rangle = \int \left\langle \vec T (x), \omega(x)\right\rangle\d\|T\|(x).$$
The mass of $T$, denoted by $\mass{T}$, is the total mass of $\|T\|$. For a positive measure $\mu$, we likewise use $\mass{\mu}$ for its total mass. The symbol $\llbracket M\rrbracket$ denotes integration of a form over the oriented rectifiable set $M$ with unit multiplicity. The symbol $\mathcal D^k$ denotes smooth compactly supported differential $k$-forms. Weak${}^*$ convergence of $k$-currents means convergence tested against $\mathcal D^k$. We denote the boundary of $T$ by $\partial T$ and the pushforward of a current by $\varphi_*T$. We denote by $\mathbf F$ the ambient integral flat norm
in $\R^d$.

We distinguish the current pushforward $(\cdot)_*$ from the measure pushforward $(\cdot)_\#$ introduced above. This will be especially important for varicurrents later, for which both $\varphi_\#\mu$ and $\timepush{\varphi}{\mu}$ are well-defined but distinct (cf. Definition \ref{def:time-pushforward}).

On classes with common compact support and uniformly bounded
mass and boundary mass, flat convergence is equivalent to
weak${}^*$ convergence, and every sequence has an integral
limit along a subsequence.
Mass is weak${}^*$ lower semicontinuous and boundary is weak${}^*$
continuous.
We also use the local versions of these statements in space-time;
see \cite[Propositions 2.4--2.5]
{rindlerSpacetimeIntegralCurrents2023}.

\subsubsection*{Space-time currents}
We use the terminology for space-time integral currents
as in \cite{rindlerSpacetimeIntegralCurrents2023,
bonicattoTransportCurrentsGeometric2025}.
We recall the conventions and results used below.

The time and spatial projections are
$$
\tproj{(t,x)}=t,
\qquad
\xproj{(t,x)}=x,\qquad (t,x)\in\R^{1+d}.
$$
We extend the definition of $\tproj{}$ and $\xproj{}$ to $m$-vectors $\eta\in\Wedge^m\R^{1+d}$ by setting
$$
\xproj{\eta}:=(\Wedge^m\mathbf p)\eta,
\qquad
\tproj{\eta}:=(\Wedge^m\mathbf t)\eta,
$$
note that $\tproj{}$ will only be used for $m=1$. We use the contraction convention
$$
\langle\eta\mres\alpha,\beta\rangle
=\langle\eta,\alpha\wedge\beta\rangle.
$$
Spatial multivectors are identified with their lift to space-time by $\Wedge^k\R^d\simeq \Wedge^k(\{0\}\times\R^d)$. We will also consider $\mathbf t, \mathbf p$ to be defined on product spaces which have $\R^{1+d}$ as their first component, so, for instance, $\mathbf p(t,x,z) := x$ for $(t,x,z) \in \R^{1+d}\times Z$ with an arbitrary set $Z$.

For $S\in\intcur{1+k}^{\loc}(I\times\overline\Omega)$, the slice with respect to $\tproj{}$ at a.e.\  $t\in I$ is denoted by $S|_t\in\intcur{k}(\{t\}\times\overline\Omega)$. We denote its spatial projection by $S(t):=\mathbf p_*(S|_t)\in\intcur{k}(\overline\Omega)$.

If
$$
\partial S\mres(\interior{I}\times\R^d)=0,
$$
then $\partial S(t)=0$ for a.e.\ $t$, and the cylinder formula gives
\begin{equation}\label{eqn:current-cylinder-formula}
\partial(S\mres(s,t))=S|_t-S|_s,
\qquad
\partial\,\mathbf p_*(S\mres(s,t))=S(t)-S(s) \quad \text{for a.e.\ }s<t \text{ in } I.
\end{equation}

\begin{definition}[Variation of a space-time current {\cite[Section~3.1, (3.1)]{rindlerSpacetimeIntegralCurrents2023}}]
Let $S \in \intcur{1+k}^{\loc}(I \times\overline\Omega)$. The \emph{variation} of $S$ over a Borel set of times $A\subset I$ is
$$
\variation{S;A}=\int_{A\times \mathbb{R}^d}|\xproj{\vec{S}}| \d\norm{S}.
$$
\end{definition}

\begin{lemma}[{\cite[Section~3.2, (3.5)]{rindlerSpacetimeIntegralCurrents2023}}]\label{lem:flatdist-leq-var}
If $S \in \intcur{1+k}^{\loc}(I \times\overline\Omega)$ satisfies $\partial S\mres(\interior{I}\times\R^d)=0$, then
\begin{equation}\label{eq:flat-distance-by-current-variation}
\mathbf F(S(t)-S(s))
\leq\variation{S;(s,t)}
\quad\text{for a.e.\ }s<t.
\end{equation}
\end{lemma}

\begin{definition}[Negligible-Criticality condition {\cite[Definition~7.1]{bonicattoTransportCurrentsGeometric2025}}]
We say that $S\in\intcur{1+k}^{\loc}(I\times\R^d)$ satisfies the \emph{negligible-criticality condition} if
\begin{equation*}\tag{NC}\label{eqn:NC-current}
    \norm{S}\mres \left\{(t, x) \in I \times \mathbb{R}^d: |\vec S(t,x) \mres\mathrm dt|=0\right\} =0.
\end{equation*}
\end{definition}

\begin{definition}[Geometric derivative {\cite[Section~7.3, (7.2)]{bonicattoTransportCurrentsGeometric2025}}]\label{def:geometric-derivative}
Let $S\in\intcur{1+k}^{\loc}(I\times\R^d)$ satisfy
\eqref{eqn:NC-current}. For $\|S\|$-a.e.\ $(t,x)\in I\times\R^d$, write
\[
\vec S(t,x)=\xi(t,x)\wedge\vec S|_t(t,x),
\qquad
\xi(t,x)\perp\vec S|_t(t,x),
\qquad
\tproj{\xi(t,x)}>0.
\]
Then $\xi(t,x)$ is uniquely defined and the \emph{geometric derivative} of $S$ is
\[
\frac{D}{Dt}S(t,x)
:=
\frac{\xproj{\xi(t,x)}}{|\tproj{\xi(t,x)}|}
\in\R^d,\qquad \text{for }\|S\|\text{-a.e.\ }(t,x)\in I\times\R^d.
\]
\end{definition}

\begin{proposition}[Time disintegration of \eqref{eqn:NC-current} currents, cf.~{\cite[Section~4]{bonicattoTransportCurrentsGeometric2025}}]\label{prop:time-disintegr-NC-current}
Let $S\in\intcur{1+k}^{\loc}(I\times\overline\Omega)$ satisfy \eqref{eqn:NC-current}. Then it holds that $\tproj{}_\#\|S\|\ll \mathcal L^1$ and
$$
\langle S,\omega\rangle =\int_I\int_{\overline\Omega}\left \langle \left(e_0+\frac{D}{Dt}S(t,\cdot)\right)\wedge \vec S|_t,\omega(t,\cdot)\right \rangle \d\|S(t)\|\d t, \quad\omega\in\mathcal D^{1+k}(I\times\overline\Omega),
$$
where $e_0$ denotes the unit vector in the time direction.
\end{proposition}

\begin{proposition}[Coarea formula for time slices {\cite[Section~2.4 and Lemma~4.1]{bonicattoTransportCurrentsGeometric2025}}]\label{prop:coarea-time-slices}
Let $S\in\intcur{1+k}^{\loc}(I\times\overline\Omega)$.
For every nonnegative Borel function $g$,
$$
\int_I\int_{\overline\Omega}
g(t,x)\,\d\|S(t)\|(x)\,\d t=\int_{I\times\overline\Omega}
g(t,x)|\vec S(t,x)\mres\mathrm dt|\,\d\|S\|(t,x)
.
$$
Moreover,  for every $\omega\in \mathcal D^k(I\times\overline\Omega)$
$$
\int_I\langle S(t),\omega(t,\cdot)\rangle\,\d t
=
\langle S,\mathrm dt\wedge\omega\rangle.
$$
\end{proposition}

\begin{proposition}[Compactness of space-time integral currents, {\cite[Theorem~3.7]{rindlerSpacetimeIntegralCurrents2023}}]
\label{prop:current-compactness}
Let $(S_j)_{j\in\N}\subset
\intcur{1+k}^{\loc}(I\times\overline\Omega)$ satisfy
\[
\operatorname{supp}\partial S_j
\subset \partial I\times\overline\Omega,
\qquad
\sup_{j\in\N}\mass{\partial S_j}<\infty,
\]
and
\[
\sup_{j\in\N}\mass{S_j\mres J}<\infty
\qquad
\text{for every compact interval }J\subset I.
\]
Then there exist
$S\in\intcur{1+k}^{\loc}(I\times\overline\Omega)$
with
$\operatorname{supp}\partial S
\subset\partial I\times\overline\Omega$
and a subsequence, not relabeled, such that
\begin{align*}
S_j&\weakstarto S
&&\text{in }\intcur{1+k}^{\loc}(I\times\overline\Omega),\\
S_j(t)&\weakstarto S(t)
&&\text{in }\intcur{k}(\overline\Omega)
\quad\text{for a.e.\ }t\in I.
\end{align*}
Moreover,
\[
\mass{S(t)}
\leq\liminf_{j\to\infty}\mass{S_j(t)}
\qquad\text{for a.e.\ }t\in I.
\]
\end{proposition}

\subsection{Definition and convergence of varicurrents}

The examples in Section~\ref{ssec:examples-currents} show why a space-time current alone may not retain the information needed to evaluate the dissipation potential after the passage to a weak limit. A space-time varicurrent records a distribution of orienting $(1+k)$-vectors at each space-time point, whose barycenter gives the associated current. We use simple unit vectors so that each recorded orientation represents an oriented tangent $(1+k)$-plane, its weight is carried by the measure. Requiring the barycenter to be an integral current preserves the current structure while allowing the measure to record oscillations and cancellations invisible to it.

\begin{definition}[Orienting $(1+k)$-vectors]
We define the set of \emph{orienting $(1+k)$-vectors} as
$$\orient:=\left\{\eta\in\Wedge^{1+k}\R^{1+d}: \eta\text{ simple, unit}\right\}.$$
\end{definition}

We use the domain $\Omega$ and time interval $I$ fixed in Section~\ref{sec:preliminaries}.

\begin{definition}[Space-time varicurrent]\label{def:space-time-varicurrent}
For every positive locally finite measure $\mu \in \mathcal M_{\loc}^+(I \times \overline\Omega \times \orient)$ we define the corresponding space-time $(1+k)$-current $\stcurr{\mu}\in\mathcal D_{1+k}^{\loc}(I\times\overline\Omega)$ by setting 
$$\langle \stcurr{\mu},\omega\rangle:=\int \langle \eta,\omega(t,x)\rangle\d\mu(t,x,\eta),\quad \omega\in\mathcal D^{1+k}(I\times\overline\Omega).$$
We call $\mu$ a \emph{space-time integral $(1+k)$-varicurrent} if it satisfies
$$
\stcurr{\mu}\in\intcur{1+k}^{\loc}(I\times\overline\Omega).
$$
We denote the space of all space-time integral varicurrents by
$$
\varicur{1+k}^{\loc}(I\times\overline\Omega):=\left\{\mu\in \mathcal M_{\loc}^+(I \times \overline\Omega \times \orient) :\stcurr{\mu}\in\intcur{1+k}^{\loc}(I\times\overline\Omega)\right\}.
$$

\end{definition}
For a Borel set of times $A\subset I$, we adopt the convention
$$
\stcurr{\mu}\mres A:=\stcurr{\mu}\mres(A\times\overline\Omega),
\qquad
\mu\mres A:=\mu\mres(A\times\overline\Omega\times\orient).
$$

\begin{lemma}[Properties of $\stcurr{\cdot}$]\label{lem:barycenter-current-properties}
The operation $\stcurr{\cdot}\colon \mathcal M_{\loc}^+(I \times \overline\Omega \times \orient) \to\mathcal D_{1+k}^{\loc}(I\times\overline\Omega)$ satisfies the following properties.
\begin{enumerate}
    \item (Linearity) For $\mu,\nu\in\mathcal M_{\loc}^+(I \times\overline\Omega\times\orient)$ and $a,b\geq 0$ it holds that $$\stcurr{a\mu+b\nu}=a\stcurr{\mu}+b\stcurr{\nu}.$$
    \item  (Mass-measure domination) For $\mu\in \mathcal M_{\loc}^+(I \times\overline\Omega\times\orient)$, the mass measure of $\stcurr{\mu}$ satisfies $$\|\stcurr{\mu}\|\leq (\tproj{},\xproj{})_\# \mu\quad \text{as measures in } \mathcal M^+_{\loc}(I\times\overline\Omega),$$ where the right-hand side is a shorthand for $((t,x,\eta) \mapsto (t,x))_\# \mu$. In particular, the current $\stcurr{\mu}$ has locally finite mass and for every compact $J\subset I$ it holds that $$ \mass{\stcurr\mu \mres J}\leq \mu(J\times\overline\Omega\times\orient).$$
    \item (Weak${}^*$ continuity) For a sequence $\mu_i,\mu\in\mathcal M_{\loc}^+(I \times\overline\Omega\times\orient)$, $i\in\N$ it holds that $$\mu_i\weakstarto\mu \text{ in }\mathcal M_{\loc}^+(I \times\overline\Omega\times \orient) \quad\implies \quad\stcurr{\mu_i}\weakstarto\stcurr{\mu}\text{ in }\mathcal D_{1+k}^{\loc}(I\times\overline\Omega).$$
\end{enumerate}
\end{lemma}
\begin{proof}
Linearity is clear from the definition. To show the mass-measure domination, let $\mu\in \mathcal M_{\loc}^+(I \times\overline\Omega\times\orient)$ and $\omega\in\mathcal D^{1+k}(I\times\overline\Omega)$. Then we estimate, using that $|\eta|=1$ for all $\eta\in\orient$,
$$
\langle \stcurr{\mu},\omega\rangle = \int \langle \eta,\omega(t,x)\rangle\d\mu = \int \omega(t,x) \d\mu(t,x,\eta)=\int \omega \d(\tproj{},\xproj{})_\#\mu.
$$
As $\omega$ was arbitrary, this shows the inequality $\|\stcurr{\mu}\|\leq (\tproj{},\xproj{})_\#\mu$ as measures.

For the weak${}^*$ continuity, take such a converging sequence $$\mu_i\weakstarto\mu \text{ in }\mathcal M_{\loc}^+(I \times\overline\Omega\times \orient)$$ and let $\omega\in\mathcal D^{1+k}(I\times\overline\Omega)$. Compute
$$\langle \stcurr{\mu_i},\omega\rangle = \int \langle\eta,\omega(t,x)\rangle\d\mu_i(t,x,\eta)\to\int \langle\eta,\omega(t,x)\rangle\d\mu(t,x,\eta)=\langle \stcurr{\mu},\omega\rangle,$$ as the integrand is continuous and compactly supported. This verifies the convergence $\stcurr{\mu_i}\weakstarto\stcurr{\mu}$ in $\mathcal D_{1+k}(I\times\overline\Omega)$.
\end{proof}

\begin{remark}
    If $\mu_i\in\varicur{1+k}(I\times\overline\Omega)$ weakly${}^*$ converge to $\mu\in\mathcal M^+(I\times\overline\Omega\times\orient)$, it does not automatically follow that $\mu\in\varicur{1+k}(I\times\overline\Omega)$. Indeed, the previous lemma shows that $\stcurr{\mu}$ is a current of finite mass, but it is not necessarily integral. Integrality requires additional boundary control, see Lemma~\ref{lem:compactness-varicurrents}.
\end{remark}

\begin{definition}[Elementary varicurrent]
For an $S\in\intcur{1+k}^{\loc}(I\times\overline\Omega)$, we define the corresponding \emph{elementary varicurrent} $\eym{S}\in\varicur{1+k}^{\loc}(I\times\overline\Omega)$ as 
$$
\eym{S}(\mathrm dt,\mathrm dx,\mathrm d\eta):=\|S\|(\mathrm dt,\mathrm dx)\otimes\delta_{\vec S(t,x)}(\mathrm d\eta),
$$
that is,
$$
\int g(t,x,\eta) \d \eym{S}(t,x,\eta)=  \int g\left(t,x,\vec S(t,x)\right) \d\|S\|(t,x), \quad g\in C_c(I\times\overline\Omega\times\orient).
$$
Conversely, $\mu\in\varicur{1+k}(I\times\overline\Omega)$ is called \emph{elementary} if $\mu = \eym{S}$ for some $S\in\intcur{1+k}^{\loc}(I\times\overline\Omega)$.
\end{definition}

As a consistency check, it is readily apparent that $\stcurr{\eym{S}}=S$ for any $S\in\intcur{1+k}^{\loc}(I\times\overline\Omega)$. Indeed, fixing $\omega\in \mathcal D^{1+k}(I\times\overline\Omega)$ it follows from the definitions that 
\begin{align*}
    \langle \stcurr{\eym{S}},\omega\rangle = \int\langle \eta,\omega(t,x)\rangle\d\eym{S}(t,x,\eta)=\int\langle \vec S(t,x),\omega(t,x)\rangle\d\|S\|(t,x)=\langle S,\omega\rangle.
\end{align*}

\subsection{Disintegration/reduction}

A varicurrent can carry information that is invisible to its barycenter: different orientations may cancel without the measure itself vanishing. To describe how the measure is related to the current, the following decomposition splits the varicurrent into a part represented relative to the barycenter mass and a residual part with zero barycenter. The residual can be removed without changing the associated current.

\begin{lemma}[Disintegration representation]\label{lem:disintegration-representation}
Let $\mu\in\varicur{1+k}^{\loc}(I\times\overline\Omega)$. There exists a decomposition 
\begin{equation}\label{eqn:disintegration-decomposition}
\mu= \|\stcurr{\mu}\| \otimes \mu^{t,x} + \nu\quad \text{ with }\stcurr{\nu}=0,
\end{equation}
where $\mu^{t,x}\in \mathcal M^+(\orient)$ for $\|\stcurr{\mu}\|$-a.e. $(t,x)$ and $\nu\in\varicur{1+k}^{\loc}(I\times\overline\Omega)$.
For any such $(\mu^{t,x},\nu)$ it holds
\begin{equation}\label{eqn:orienter-barycenter-spacetime}
\vstcurr{\mu}(t,x)=\int_{\orient}\eta\d\mu^{t,x}(\eta),\quad \text{for }\|\stcurr{\mu}\|\text{-a.e. }(t,x).
\end{equation}
Moreover, there exists a decomposition $(\mu^{t,x}_0,\nu_0)$, unique up to a $\|\stcurr{\mu}\|$-null set of $(t,x)$, such that  $\nu_0$ is minimal among all such decompositions. This unique decomposition is characterized by the additional property
$$(\tproj{},\xproj{})_\#\nu_0\perp \|\stcurr{\mu}\|.$$
\end{lemma}
\begin{proof}
Let us use the disintegration theorem (see e.g.\ \cite[Theorem 2.28]{ambrosioFunctionsBoundedVariation2000}) on $\mu$ to obtain
$$
\mu(\mathrm dt,\mathrm dx,\mathrm d\eta)=\lambda(\mathrm dt,\mathrm dx)\otimes \rho^{t,x}(\mathrm d\eta)
$$
for some  $\lambda \in \mathcal M_{\loc}^+(I\times\overline\Omega)$, and $\rho^{t,x}\in\mathcal M^+(\orient)$ with $\mass{\rho^{t,x}}=1$, for $\lambda$-a.e.\ $(t,x)\in I\times\overline\Omega$. Note that $\lambda=(\tproj{},\xproj{})_\#\mu$, so by Lemma \ref{lem:barycenter-current-properties} it holds that $\|\stcurr{\mu}\|\leq\lambda$.
Take the Lebesgue decomposition of $\lambda$ with respect to $\|\stcurr{\mu}\|\in\mathcal M^+_{\loc}(I\times\overline\Omega)$, that is, $\lambda^a,\lambda^s\in\mathcal M_{\loc}^+(I\times\overline\Omega)$ satisfying
$$
\lambda = \lambda^a+\lambda^s, \quad \lambda^a \ll \|\stcurr{\mu}\|, \quad \lambda^s \perp \|\stcurr{\mu}\|.
$$
We now define 
$$
\mu_0^{t,x} := \frac{\mathrm d\lambda^a}{\mathrm d\|\stcurr{\mu}\|}(t,x)\rho^{t,x}\quad \|\stcurr{\mu}\|\text{-a.e.\ }(t,x), \quad \text{and} \quad  \nu_0:=\lambda^s\otimes \rho^{t,x}.
$$
With this, the decomposition
$$
\mu= \|\stcurr{\mu}\| \otimes \mu_0^{t,x}+\nu_0
$$
holds, and it remains to prove $\stcurr{\nu_0}=0$. For this, take $\omega\in\mathcal D^{1+k}(I\times\overline\Omega)$. We have
\begin{align*}
    \langle S[\mu],\omega\rangle = &
    \int\int\left\langle\eta,\omega(t,x)\right\rangle\d\mu_0^{t,x}(\eta)\d\|\stcurr{\mu}\|(t,x)
   +
   \int\int\left\langle\eta,\omega(t,x)\right\rangle\d\rho^{t,x}(\eta)\d\lambda^s(t,x).
\end{align*}
On the other hand,
$$
\langle \stcurr{\mu},\omega\rangle =\int\left\langle\vstcurr{\mu}(t,x),\omega(t,x)\right\rangle \d\|\stcurr{\mu}\|(t,x).
$$
Comparing the last two equalities gives us
\begin{align*}
&\int\left\langle\vstcurr{\mu}(t,x)-\int\eta\d\mu_0^{t,x}(\eta),\omega(t,x)\right\rangle \d\|\stcurr{\mu}\|(t,x)\\
&\quad=\int\int\left\langle\eta,\omega(t,x)\right\rangle\d\rho^{t,x}(\eta)\d\lambda^s(t,x)=\langle\stcurr{\nu_0},\omega\rangle.
\end{align*}
As the measures $\|\stcurr{\mu}\|$ and $\lambda^s$ are mutually singular, this implies 
$\stcurr{\nu_0}=0$. Thus the constructed decomposition $\mu_0^{t,x}$, $\nu_0$ satisfies \eqref{eqn:disintegration-decomposition}, proving the existence statement.

Now let $\mu^{t,x}$, $\nu$ be any decomposition satisfying \eqref{eqn:disintegration-decomposition}.
We find that
\begin{align*}
0&=\langle \stcurr{\nu},\omega\rangle\\
&=\left\langle\stcurr{\mu-\|\stcurr{\mu}\|\otimes \mu^{t,x}},\omega\right\rangle\\
&=\int\left\langle\vstcurr{\mu}(t,x)-\int\eta\d\mu^{t,x}(\eta),\omega(t,x)\right\rangle \d\|\stcurr{\mu}\|(t,x)
\end{align*}
which shows that \eqref{eqn:orienter-barycenter-spacetime} holds.

It remains to prove the minimality and uniqueness assertion. Let \(N\subset I\times\overline\Omega\) be a Borel set such that
\begin{equation}\label{eqn:singular-support-lambda-s}
\|\stcurr{\mu}\|(N)=0,
\qquad
\lambda^s\bigl((I\times\overline\Omega)\setminus N\bigr)=0.
\end{equation}
Then we have 
$$ \|\stcurr{\mu}\|\otimes\mu^{t,x}\mres N\times\orient =0, \quad \|\stcurr{\mu}\|\otimes\mu_0^{t,x}\mres N\times\orient =0,\quad \nu_0\mres(N\times \mathcal O)=\nu_0$$so that
$$
\nu\mres(N\times\orient)=\mu\mres(N\times\orient)
=
\nu_0\mres(N\times\orient)
=
\nu_0,
$$
and consequently $\nu_0\leq\nu$, which proves that $\nu_0$ is minimal.

Finally, let us show the characterization of minimality. Write
$$
\nu -\nu_0  = \|\stcurr{\mu}\|\otimes (\mu_0^{t,x}-\mu^{t,x}),
$$
from which it follows that $$(\tproj{},\xproj{})_\#(\nu-\nu_0)\ll\|\stcurr{\mu}\|.$$
From the construction of $\nu_0$, we have $$(\tproj{},\xproj{})_\#\nu_0=\lambda^s \perp \|\stcurr{\mu}\|.$$
The last two statements imply that
$$(\tproj{},\xproj{})_\#(\nu-\nu_0)=0 \iff (\tproj{},\xproj{})_\#\nu \perp\|\stcurr{\mu}\|.$$
Finally, we already know that the measure $\nu-\nu_0$ is positive, hence it follows that
$$
(\tproj{},\xproj{})_\#(\nu-\nu_0)=0 \iff \nu - \nu_0=0.
$$
Altogether, this shows that $\nu=\nu_0$ if and only if $(\tproj{},\xproj{})_\#\nu \perp\|\stcurr{\mu}\|$. This finishes the proof.
\end{proof}
\begin{remark}[Reduction interpretation]\label{rem:reduction-interpretation}

In the context of Lemma \ref{lem:disintegration-representation}, one may interpret
$$
\tilde\mu = \|\stcurr{\mu}\|\otimes \mu^{t,x} \in \varicur{1+k}^{\loc}(I\times\overline\Omega)
$$
as a reduction of $\mu\in\varicur{1+k}^{\loc}(I\times\overline\Omega)$. The mass singular with respect to the barycenter mass measure $\|\stcurr{\mu}\|$ is contained in $\nu$, and is deleted when passing to $\tilde\mu$. A zero-barycenter contribution that is absolutely continuous with respect to this mass measure may or may not be deleted. The lemma gives a maximal reduction $\tilde\mu_0$, which retains the largest measure, and a corresponding minimal residual measure $\nu_0$, which deletes the least. This minimal residual removes only the mass singular with respect to the barycenter mass measure.

One may also wish to remove as much of the zero-barycenter contribution as possible, so that the reduced measure retains less information irrelevant to the current. Such a minimal reduction need not have a canonical choice: different parts of the orientation measure can cancel the same barycenter contribution. For instance, consider the case
$$
\mu=\mu_1+\mu_2+\mu_3, \quad 
\stcurr{\mu_2+\mu_3}=0, \quad\stcurr{\mu_1+\mu_3}=0, \quad \mu_1\neq\mu_2.
$$
(E.g.\ take $S\in\intcur{1+k}(I\times\overline\Omega)$, $\mu_1:=\eym{S}$, $\mu_2$ non-elementary with $\stcurr{\mu_2}=S$, and $\mu_3:=\eym{-S}$.) Deleting $\mu_2+\mu_3$ leaves $\mu_1$, whereas deleting $\mu_1+\mu_3$ leaves $\mu_2$. Thus both $\mu_1$ and $\mu_2$ are valid reductions with the same barycenter, but the construction does not single out either one as preferred. Whether a useful canonical notion of minimal reduction can be specified is a separate question; we use only the decomposition established above.
\end{remark}

\begin{remark}\label{rem:elementary-varicur-disintegration}
If $\mu\in\varicur{1+k}^{\loc}(I\times\overline\Omega)$ is an elementary varicurrent, then it has exactly one decomposition, namely
$$
\mu^{t,x}=\delta_{\vstcurr{\mu}(t,x)}, \qquad \nu=0.
$$
\end{remark}
\begin{lemma}\label{lem:mass-elem-varicur}
Let $\mu\in\varicur{1+k}(I\times\overline\Omega)$, and let $\mu^{t,x},\nu$ be defined from any decomposition in Lemma~\ref{lem:disintegration-representation}. Then it holds
$$\mass{\mu^{t,x}} \geq 1 \quad \text{for }\|\stcurr{\mu}\|\text{-a.e.\ }(t,x).$$ 
Moreover, $\mu$ is an elementary varicurrent if and only if equality holds and $\nu=0$. Equivalently, $\mu$ is elementary if and only if
$$
(\tproj{},\xproj{})_\#\mu = \|\stcurr{\mu}\|.
$$
\end{lemma}
\begin{proof}
Recall that orienting $(1+k)$-vectors have unit length. Using the orienter identity \eqref{eqn:orienter-barycenter-spacetime}, we obtain for $\|\stcurr{\mu}\|$-a.e.\ $(t,x)$ that 
$$1 = \left|\vstcurr{\mu}(t,x)\right| = \left| \int_{\orient} \eta \d\mu^{t,x}(\eta) \right| \leq \int_{\orient} |\eta| \d\mu^{t,x}(\eta) = \mass{\mu^{t,x}},$$
which proves the lower bound.
Furthermore, because the unit ball in $\Wedge^{1+k}\R^{1+d}$ is strictly convex, the triangle inequality used above holds as an equality if and only if the measure $\mu^{t,x}$ is entirely supported on a single point in $\orient$. Thus, the equality holds if and only if $\mu^{t,x}=\delta_{\vstcurr{\mu}(t,x)}$. Along with the condition $\nu=0$ this is equivalent to $\mu$ being elementary (see Remark \ref{rem:elementary-varicur-disintegration}).
\end{proof}

\subsection{Slicing of varicurrents}

We next introduce Lebesgue-time slices of space-time varicurrents, in analogy with the theory of space-time currents~\cite{rindlerSpacetimeIntegralCurrents2023}. The construction used here assumes both \eqref{eqn:NC-varicurrent} and absolute continuity of the time marginal. These are the conditions that we will use in our variational constructions (Section \ref{sec:WED}), although they are not strictly necessary for slicing, cf. Remark \ref{rem:more-general-slicing}.

The slices must retain velocity information. We therefore use \emph{velocity-attached slices}, which record the joint distribution of spatial orientations and normal velocities. In the case of currents, it is sufficient to define the velocity as a separate object (see Definition \ref{def:geometric-derivative}). Here a joint measure is needed because, loosely speaking, we need to know which velocity belongs to which tangent plane.

\begin{definition}[Spatial orienting vectors]
We define the set of all \emph{purely spatial orienting $(1+k)$-vectors} to be 
$$\orientSpatial:=\orient\cap \Wedge^{1+k}(\{0\}\times \R^d).$$
We also introduce the set of unit simple spatial $k$-vectors
\begin{equation}\label{eq:spatial-k-orientations}
\orientSpatial^k
:=
\left\{
\tau\in\Wedge^k\R^d:
\tau\text{ is simple and }|\tau|=1
\right\}.
\end{equation}
\end{definition}
Geometrically, a purely spatial orientation represents an (oriented) tangent plane contained in a single time slice. Motion along this plane traverses a spatial distance over zero elapsed time, so we associate it with infinite speed. This interpretation matters for the superlinear dissipation potential: such orientations cannot carry positive mass in a trajectory of finite cost.

We first decompose each non-purely spatial orientation into spatial and temporal parts. This separates the orientation of the spatial slice from the direction in which it moves and prepares the definition of normal velocity.

\begin{lemma}[Space-time decomposition]\label{lem:space-time-decomposition}
Any non-purely spatial orienting $(1+k)$-vector $\eta \in \orient \setminus \orientSpatial$ can be decomposed as 
\begin{align*}
\eta=\xi(\eta)\wedge\tau(\eta),\quad \text{with } \ \xi(\eta)\in \R^{1+d}, \tau(\eta)\in\Wedge^{k}(\{0\}\times\R^d) \text{ such that}
\\
|\xi(\eta)|=|\tau(\eta)|=1,\ \xi(\eta)\perp\tau(\eta),\ \tproj{\xi(\eta)}>0.
\end{align*}
Moreover, this decomposition is unique. 
\end{lemma}

\begin{proof}
Let $\eta \in \orient \setminus \orientSpatial$. Then we put $\tau(\eta) := \frac{\eta \mres \d t}{|\eta \mres \mathrm{d}t|}$, and  $\xi(\eta) := \eta \mres \tau(\eta)^\flat$, where $\flat\colon \Wedge^k\R^{1+d}\to(\Wedge^k\R^{1+d})^*$ is the canonical isomorphism. This can be readily seen to satisfy the required decomposition. The uniqueness is provided by fixing the sign $\tproj{\xi(\eta)}>0$.
\end{proof}

The preceding decomposition allows us to extract the normal velocity from a space-time orientation. Its vector $\xi(\eta)$ has a positive time component, so the ratio of its spatial and temporal components describes the spatial displacement per unit time. This is the velocity that will enter the dissipation potential.

\begin{definition}[Velocity vector]\label{def:velocity-vector}
For $\eta\in\orient\setminus\orientSpatial$, we define its corresponding \emph{velocity vector} $v(\eta)\in\R^d$
by setting
$$
v(\eta):=
    \dfrac{\xproj{\xi(\eta)}}{\tproj{\xi(\eta)}}\in\R^{d}.
$$
Note that here we interpret $\tproj{\xi(\eta)}$ as a \emph{scalar}, more precisely as $\xi(\eta) \mres \mathrm dt$.
Additionally,  for $\eta\in\orientSpatial$ we use the convention $|v(\eta)|:=\infty$, whilst $v(\eta)$ itself remains undefined.
\end{definition}

\begin{remark}
For $S\in\intcur{1+k}^{\loc}(I\times\overline\Omega)$, the velocity corresponding to the orientation of $\vec S$ is exactly the geometric derivative $\frac{D}{Dt}S$ at non-critical points:
$$v(\vec S(t,x))=\frac{D}{Dt}S(t,x)\quad \text{for $\|S\|$-a.e. }(t,x)\in \operatorname{supp}\|S\|\ \text{with }\vec S(t,x)\in\orient\setminus\orientSpatial.$$
\end{remark}

\begin{lemma}\label{lem:velocity-cont-lsc}
The maps $\eta\mapsto v(\eta)\colon \orient \setminus \orientSpatial \to \R^d$ and $\eta \mapsto |v(\eta)|\colon\orient\to [0,\infty]$ are continuous.
\end{lemma}
\begin{proof}
In the case $\eta\notin\orientSpatial$, it holds that $$\tau(\eta) = \frac{\eta\mres \mathrm dt}{|\eta\mres \mathrm dt|}.$$
Given this, we may compute $$\xi(\eta)=\eta\mres \tau^\flat = \eta\mres \left(\frac{\eta\mres \mathrm dt}{|\eta\mres \mathrm dt|}\right)^{\flat}$$
hence we see that $\xi(\eta)$ depends continuously on $\eta$, and thus $\eta\mapsto v(\eta)$ is continuous on $\orient\setminus\orientSpatial$.

To address the continuity near $\orientSpatial$, observe that $\eta$ approaching a point in $\orientSpatial$ implies $\tproj{\xi(\eta)}\to 0$. Thus also $|\xproj{\xi(\eta)}|\to 1$, due to  $|\tproj{\xi(\eta)}|^2+|\xproj{\xi(\eta)}|^2=1$ as $\xi$ is unit, and consequently we obtain $|v(\eta)|=|\xproj{\xi(\eta)}|/|\tproj{\xi(\eta)}|\to \infty$ as $\eta$ approaches $\orientSpatial$. This concludes the proof.
\end{proof}

\begin{remark}
The convention $|v(\eta)|=\infty$ for $\eta\in\orientSpatial$ expresses motion over zero elapsed time, or an infinitesimal jump. It specifies only an infinite scalar speed, while the direction of the velocity vector $v(\eta)$ itself remains undefined. The purely spatial tangent plane itself does not determine a direction of traversal for that jump. 
\end{remark}

The following identities relate the velocity to the spatial and temporal orientation parts.

\begin{lemma}\label{lem:v-tau-p-x-identities}
Let $\eta\in\orient\setminus\orientSpatial$. Then its space-time decomposition $\xi(\eta)$, $\tau(\eta)$ from Lemma \ref{lem:space-time-decomposition}, and the velocity vector from Definition \ref{def:velocity-vector} satisfy the following properties:
\begin{align*}
v(\eta)&\perp\operatorname{span}\tau(\eta),
&\xi(\eta)&=\frac{e_0+v(\eta)}{\sqrt{1+|v(\eta)|^2}},\\
|\eta \mres \mathrm dt|&=\frac{1}{\sqrt{1+|v(\eta)|^2}},
&|\xproj{\eta}|&=|v(\eta)|\,|\eta \mres \mathrm dt|.
\end{align*}
\end{lemma}
\begin{proof}
    These are immediate from the definitions.
\end{proof}

For a space-time current, the negligible-criticality condition \eqref{eqn:NC-current} requires zero mass on the critical set, where the orientation is purely spatial and a finite geometric derivative is unavailable. The varicurrent condition expresses the same restriction at the level of the orientation measure: it assigns zero mass to purely spatial orientations. Thus the velocity vector is defined at almost every orientation.

\begin{definition}[Negligible-criticality condition]\label{def:NC-varicurrent}
We say that a space-time varicurrent $\mu\in\varicur{1+k}(I\times\overline\Omega)$ satisfies the \emph{negligible-criticality condition} if
\begin{equation*}\tag{V-NC}\label{eqn:NC-varicurrent}
\mu\left(I\times\overline\Omega\times\orientSpatial \right)=0.
\end{equation*}
Further, by the convention introduced in Definition \ref{def:velocity-vector}, \eqref{eqn:NC-varicurrent} is equivalent to
$$|v(\eta)|<\infty, \quad\text{for $\mu$-a.e.\ }(t,x,\eta)\in I\times\overline\Omega\times\orient.$$
\end{definition}

\begin{remark}
For elementary varicurrents, this agrees with the current condition, namely
$$\mu \text{ satisfies \eqref{eqn:NC-varicurrent}}\quad\iff \quad \stcurr{\mu}\text{ satisfies \eqref{eqn:NC-current}}.$$
\end{remark}

We can now record the spatial orientation and its attached normal velocity together. For the velocity-attached slices, we use the phase space
\begin{equation}\label{eq:slice-phase-space}
\mathcal Z_k
:=
\left\{
(x,\tau,v)
\in
\overline\Omega\times\orientSpatial^k\times\R^d:
 v\perp\tau
\right\}.
\end{equation}
The space $\mathcal Z_k$ is the phase space of spatial $k$-orientations
with an attached normal velocity.

\begin{definition}[Velocity-attached time slices]\label{def:velocity-attached-slices}
Let
$
\mu\in\varicur{1+k}^{\loc}(I\times\overline\Omega)
$
satisfy \eqref{eqn:NC-varicurrent} and $\tproj{}_\#\mu\ll \mathcal L^1$.
We define the coarea-weighted measure
\begin{equation}\label{eq:Sigma-mu-defn}
\Sigma_\mu
:=
(\tproj{},\xproj{},\tau,v)_\#
\bigl(|\eta \mres \mathrm dt|\mu\bigr)
\in
\mathcal M^+_{\loc}(I\times\mathcal Z_k).
\end{equation}
More explicitly, for every nonnegative Borel function
$g:I\times \mathcal Z_k \to[0,+\infty]$,
\begin{equation}\label{eq:weighted-phase-space-measure-test}
\int g(t,x,\tau,v)\,\d\Sigma_\mu(t,x,\tau,v)
=
\int
 g\bigl(t,x,\tau(\eta),v(\eta)\bigr)
 |\eta \mres \mathrm dt|\,\d\mu(t,x,\eta).
\end{equation}
Then we define the measures
\begin{equation*}
\slice{\mu}{t}\in\mathcal M^+(\mathcal Z_k), \quad \text{for a.e.\ }t\in I
\end{equation*}
obtained by disintegration of $\Sigma_\mu$ in time with respect to the Lebesgue measure, i.e.\ it is the Lebesgue-a.e.\ uniquely defined family of measures that satisfies
\begin{equation}\label{eq:definition-velocity-attached-slices}
\Sigma_\mu
=
\mathrm dt\otimes\slice{\mu}{t}.
\end{equation}
We call $\slice{\mu}{t}$ the \emph{velocity-attached time slice} of $\mu$
at time $t$.
\end{definition}

\begin{remark}\label{rem:more-general-slicing}
The hypotheses in the definition specify the class of slices used in this paper, rather than the most general possible slicing construction. If \eqref{eqn:NC-varicurrent} fails, purely spatial orientations have $|\eta \mres \mathrm dt|=0$ and do not contribute to the weighted measure $\Sigma_\mu$. If the time marginal is not absolutely continuous, a Lebesgue-time disintegration can describe only its absolutely continuous part and leaves out singular-time mass (cf.~\cite[Section~4]{bonicattoTransportCurrentsGeometric2025}). This description in the general case is possible, but remains outside the scope of this work.
\end{remark}
\begin{remark}\label{rem:slice-measures-converge}
The map $\mu\mapsto\Sigma_\mu$ is weak$^*$ continuous whenever these measures are defined. Indeed, for $g\in C_c(I\times\mathcal Z_k)$ we define its weighted pullback by extending to 0 on critical orientations:
$$
(t,x,\eta)\mapsto\begin{cases}
|\eta \mres \mathrm dt|g(t,x,\tau(\eta),v(\eta)),& \eta\in\orient\setminus\orientSpatial,
\\
0,&\eta\in\orientSpatial.
\end{cases}
$$
This is continuous as $|\eta \mres \mathrm dt|\to0$ implies $|v(\eta)|\to\infty$, eventually ending up outside the compact velocity support of $g$. Testing $\mu_i\weakstarto\mu$ against this proves $\Sigma_{\mu_i}\weakstarto\Sigma_\mu$. This continuity will be used to pass the slice-mass bound to the limit in Section~\ref{ssec:mass-bound-time}.
\end{remark}

Forgetting the velocity in a velocity-attached slice gives its orientation marginal, an analogue of current slices without remembering the geometric derivative. The velocity-attached slices remain necessary for evaluating the dissipation potential.

\begin{definition}[Slices without velocity]
Let $\mu\in\varicur{1+k}^{\loc}(I\times\overline\Omega)$ satisfy \eqref{eqn:NC-varicurrent} and $\tproj{}_\#\mu\ll \mathcal L^1$. The velocity-free spatial slice is obtained by integrating out the velocity
variable. That is, for a.e. $t\in I$, define $\slicenov{\mu}{t}\in \mathcal M^+(\overline\Omega\times\orientSpatial^k)$ by setting
$$
\slicenov{\mu}{t}:= \bigl((x,\tau,v)\mapsto(x,\tau)\bigr)_\#\slice{\mu}{t},
$$
more explicitly,
$$\int_{\overline\Omega\times\orientSpatial^k}g(x,\tau)\d\slicenov{\mu}{t}(x,\tau)=\int_{\mathcal Z_k}g(x,\tau)\d \slice{\mu}{t}(x,\tau,v),\quad g\in  C(\overline\Omega\times\orientSpatial^k).$$
\end{definition}
The two slice measures have the same total mass:
\begin{equation}\label{eq:slice-masses-with-without-velocity}
\mass{\slicenov{\mu}{t}}
=
\mass{\slice{\mu}{t}}
\qquad\text{for a.e. }t\in I.
\end{equation}
We next relate slicing a varicurrent to slicing its barycenter current. The expected compatibility is not simply a consequence of averaging: extracting the spatial orientation uses the nonlinear map $\eta\mapsto\tau(\eta)$, so it does not in general commute with taking a barycenter. For $\nu\in\mathcal M^+(\orient)$, in general
$$
\tau\left(\int\eta\,\d\nu(\eta)\right)
\neq
\int\tau(\eta)\,\d\nu(\eta).
$$
The coarea weight $|\eta \mres \mathrm dt|$ in the definition of slicing restores the required compatibility of the two slicing notions. This is expressed in the following lemma.
\begin{lemma}[Orientation of slices]
Let $\mu\in\varicur{1+k}^{\loc}(I\times\overline\Omega)$ with \eqref{eqn:NC-varicurrent} and $\tproj{}_\#\mu\ll\mathcal L^1$. Then for a.e.\ $t\in I$ it holds

\begin{equation}\label{eq:barycentre-slice-current-identity-weak}
\left\langle
\stcurr{\mu}(t),
\omega
\right\rangle
=
\int_{\overline\Omega\times\orientSpatial^k}
\langle\tau,\omega(x)\rangle
\,\d\slicenov{\mu}{t}(x,\tau)\quad\text{for all }\omega\in\mathcal D^k(\overline\Omega).
\end{equation}
In particular, it holds for a.e.\ $t\in I$ that
$$
\mass{\stcurr{\mu}(t)}\leq\mass{\slicenov{\mu}{t}}.
$$
\end{lemma}
\begin{proof}
Let $\zeta\in \mathcal D^k(I\times\overline\Omega)$. By the slicing identity in Proposition \ref{prop:coarea-time-slices} for $\stcurr{\mu}\in\intcur{1+k}^{\loc}(I\times\R^d)$ we have
$$
\int_I\langle \stcurr{\mu}(t), \zeta(t,\cdot)\rangle \d t=  \langle\stcurr{\mu},\mathrm dt \wedge\zeta\rangle.
$$
By the definition of $\stcurr{\mu}$ we compute
\begin{align*}
    \langle\stcurr{\mu},\mathrm dt \wedge\zeta\rangle = \int\langle \eta , \mathrm dt \wedge \zeta(t,x)\rangle \d\mu(t,x,\eta).
\end{align*}
As $\mu$ satisfies \eqref{eqn:NC-varicurrent}, we may decompose $\eta$ as in Lemma \ref{lem:space-time-decomposition} to obtain
$$
\int\langle \eta , \mathrm dt \wedge \zeta(t,x)\rangle \d\mu(t,x,\eta)=\int \langle\tau(\eta),\zeta(t,x)\rangle \tproj{\xi(\eta)} \d \mu(t,x,\eta).
$$
By the definition of $\Sigma_\mu$ and $\slicenov{\mu}{t}$ we obtain
$$
\int \langle\tau(\eta),\zeta(t,x)\rangle \tproj{\xi(\eta)} \d \mu(t,x,\eta)=\int \langle\tau,\zeta(t,x)\rangle \d \Sigma_\mu(t,x,\tau,v)=  \int_I\int \langle\tau,\zeta(t,x)\rangle \d \slice{\mu}{t}(x,\tau,v)\d t.
$$
Collecting all these, we have shown that 
$$
\int_I\langle \stcurr{\mu}(t), \zeta(t,\cdot)\rangle \d t=  \int_I\int \langle\tau,\zeta(t,x)\rangle \d \slicenov{\mu}{t}(x,\tau)\d t.
$$
Picking all $\zeta$ in the form $\zeta(t,x)=\psi(t)\omega(x)$ for $\psi\in C^\infty_c(I)$ and $\omega\in\mathcal D^k(\overline\Omega)$ proves the statement.
\end{proof}

\subsection{Compactness and Kinematics}

The direct method in the later WED construction (Section \ref{sec:WED}) requires limits to remain within the class of varicurrents. Weak compactness of the underlying measures alone is not enough, since the barycenter of the limit must still be an integral current. The following lemma combines mass bounds for the measures with boundary control for their barycenters. It provides the weak${}^*$ compactness needed for this purpose.

\begin{lemma}[Compactness in $\varicur{1+k}$]\label{lem:compactness-varicurrents}
Let $(\mu_i)_{i\in\N}\subset \varicur{1+k}^{\loc}(I\times\overline\Omega)$ be a sequence of varicurrents with
\begin{align*}
&\sup_{i\in\N}\mass{\mu_i\mres J}<\infty \quad \text{for every compact }J\subset I,
\\
&\operatorname{supp}\partial \stcurr{\mu_i}  \subset \partial I\times \overline\Omega, \text{ and }
\\
&\sup_{i\in\N}\mass{\partial \stcurr{\mu_i} } <\infty.
\end{align*}
Then there exists a $\mu\in\varicur{1+k}^{\loc}(I\times\overline\Omega)$ and a (non-relabeled) subsequence, such that
\begin{align*}
    \mu_i &\weakstarto\mu,  &\text{in }\mathcal M^+_{\loc}(I\times\overline\Omega\times\orient),
    \\
    \stcurr{\mu_i} &\weakstarto\stcurr{\mu}&\text{in }  \intcur{1+k}^{\loc}(I\times\overline\Omega),
    \\
    \stcurr{\mu_i}(t) &\weakstarto\stcurr{\mu}(t)&\text{in }  \intcur{k}(\overline\Omega)\text{ for a.e. }t\in I,
\end{align*}
and, moreover, for intervals $J\subset I$ with endpoints outside of a countable set
\[
\mass{\mu\mres J}\leq\liminf_{i\to\infty}\mass{\mu_i\mres J}
    \quad\text{and}\quad \mass{\stcurr{\mu}(t)}\leq\liminf_{i\to\infty} \mass{\stcurr{\mu_i}(t)}\quad \text{for a.e. }t\in I.
\]
\end{lemma}
\begin{proof}
By the assumed local mass bound on $\mu_i$ we obtain (for a subsequence) $\mu\in\mathcal M_{\loc}^+(I\times\overline\Omega\times\orient)$ with
$$\mu_i\weakstarto\mu\quad \text{in }\mathcal M_{\loc}^+(I\times\overline\Omega\times\orient)$$
and by Lemma \ref{lem:barycenter-current-properties} also
$$
\stcurr{\mu_i}\weakstarto\stcurr{\mu}\quad \text{in }\mathcal D_{1+k}^{\loc}(I\times\overline\Omega).
$$
By Lemma \ref{lem:barycenter-current-properties}
the local mass bound on $\mu_i$ implies the  local mass bound on $\stcurr{\mu_i}$, that is,
$$
\sup_{i\in\N}\mass{\stcurr{\mu_i}\mres J}<\infty \quad \text{for every compact }J\subset I.
$$
Using these bounds, by the Federer-Fleming compactness theorem, we see that for some subsequence it holds
$$\stcurr{\mu_{i_j}}\weakstarto S\quad \text{in }\intcur{1+k}^{\loc}(I\times\overline\Omega)$$
for  some $S\in\intcur{1+k}^{\loc}(I\times\overline\Omega)$. This implies $\stcurr{\mu}=S$ and consequently $\stcurr{\mu}\in\intcur{1+k}^{\loc}(I\times\overline\Omega)$, by which we verified that $\mu\in\varicur{1+k}^{\loc}(I\times\overline\Omega)$. The weak${}^*$ convergence for almost every slice is obtained by applying \cite[Theorem 3.7]{rindlerSpacetimeIntegralCurrents2023}, after possibly passing to another subsequence.

As before, for intervals $J\subset I$ with $\mu(\partial J\times\overline\Omega\times\orient)=0$ it holds
$$\mu_i\mres J \weakstarto\mu\mres J,$$
and the mass inequalities follow from the weak${}^*$ lower semicontinuity of mass on measures, resp. on integral currents.
\end{proof}

We now define the variation of a varicurrent. As for currents, this quantity measures the spatial area swept out by a trajectory.
\begin{definition}[Variation of a varicurrent]\label{def:variation-of-varicurrent}
For a space-time varicurrent $\mu\in\varicur{1+k}(I\times\overline\Omega)$ we define its \emph{variation} as
\begin{equation}\label{eqn:variation-of-varicurrent}
\variation{\mu;A} := \int_{A\times\overline\Omega\times\orient}|\xproj{\eta}|\d\mu(t,x,\eta), \quad \text{ for any Borel }A\subset I.
\end{equation}
\end{definition}

\begin{lemma}[Inequality of variations]\label{lem:variation-ineq}
For a space-time varicurrent $\mu\in\varicur{1+k}(I\times\overline\Omega)$, the following inequality holds
$$
\variation{\stcurr{\mu};A}\leq \variation{\mu;A}.
$$
\end{lemma}
\begin{proof}
Take a decomposition of $\mu$ from Lemma~\ref{lem:disintegration-representation}, so that
$$
\mu=\|\stcurr{\mu}\|\otimes\mu^{t,x} + \nu,\quad \stcurr{\nu}=0.
$$
With this, from the definitions of variation and linearity of $\xproj{}$ and using the orienter identity \eqref{eqn:orienter-barycenter-spacetime}, we obtain
\begin{align*}
\variation{\stcurr{\mu};A}
&=\int_{A\times\overline\Omega}\left|\xproj{\vstcurr{\mu}(t,x)}\right|\d\|\stcurr{\mu}\|(t,x)
\\
&=\int_{A\times\overline\Omega}\left|\xproj{\int_{\orient}\eta\d\mu^{t,x}(\eta)}\right|\d\|\stcurr{\mu}\|(t,x)
\\
&\leq\int_{A\times\overline\Omega}\int_{\orient}|\xproj{\eta}|\d\mu^{t,x}(\eta)\d\|\stcurr{\mu}\|(t,x)
\\
&\leq\int_{A\times\overline\Omega\times \orient}|\xproj{\eta}|\d\mu(t,x,\eta)
= \variation{\mu;A}. \qedhere
\end{align*}
\end{proof}

\subsubsection{Mass bound in time}\label{ssec:mass-bound-time}
In this auxiliary subsection we give a precise meaning to a uniform mass bound in time for a time-dependent measure. Applied to varicurrent slices, the bound combines with finite dissipation to control temporal variation and hence continuity of the barycenter current; see Lemma~\ref{lem:p-continuity-in-time}. The Thousand Cuts example (Example~\ref{ex:thousand-cuts-current}) shows why the bound matters: bounded costs can coexist with unbounded slice masses, resulting in an instantaneous limiting motion.

\begin{definition}[$L^\infty$ norm for measures]
For any $\Sigma\in\mathcal M_{\loc}^+(I\times\mathcal Z_k)$ with $\tproj{}_\#\Sigma \ll\mathcal L^1(I)$ we denote $$\norm{\Sigma}_{L^\infty(I;\mathcal M)}:=\operatorname*{ess\,sup}_{t\in I} \mass{\Sigma^t},$$ where we use the disintegration $$\Sigma=\mathrm dt\otimes \Sigma^t.$$
\end{definition}

\begin{definition}[Measure with bounded mass in time]
We define $$L_{w^*}^\infty(I;\mathcal M(\mathcal Z_k)):=\left\{ \Sigma\in\mathcal M_{\loc}^+(I\times\mathcal Z_k):\tproj{}_\#\Sigma \ll\mathcal L^1(I), \norm{\Sigma}_{L^\infty(I;\mathcal M)}<\infty \right\}.$$
\end{definition}

\begin{remark}[Weak${}^*$  measurability]
For $\Sigma\in\mathcal M_{\loc}^+(I\times\mathcal Z_k)$ with $\tproj{}_\#\Sigma \ll\mathcal L^1(I)$ it follows by the disintegration theorem that  $t\mapsto \Sigma^t$ is weak${}^*$ measurable, i.e. 
$$
\text{for all }\varphi\in C_b(\mathcal Z_k)\text{ the scalar function }t\mapsto \langle\Sigma^t,\varphi\rangle \text{ is measurable}.
$$
However, as already the simple Example \ref{ex:point-trajectory} shows, $t\mapsto \Sigma^t$ is generally not strongly measurable. For this reason, we use the $w^*$ subscript.
\end{remark}

The next characterization in case of finite $L^\infty(I;\mathcal M)$ will let us handle weak convergence without first comparing the disintegrated measures at individual times.

\begin{lemma}[Alternative norm]\label{lem:alt-norm}
For $\Sigma\in\mathcal M_{\loc}^+(I\times\mathcal Z_k)$ denote $$\norm{\Sigma}_*:=\sup\left\{\int\varphi\d\Sigma : \varphi\in C_c(I\times\mathcal Z_k), \norm{\varphi}_{L^1(I;L^\infty)}\leq 1\right\},$$ with the norm $$\norm{\varphi}_{L^1(I;L^\infty)}:=\int_I\norm{\varphi(t,\cdot)}_{L^\infty}\d t,\quad \varphi\in C_c(I\times\mathcal Z_k).$$
If a measure $\Sigma\in\mathcal M_{\loc}^+(I\times\mathcal Z_k)$ satisfies $\norm{\Sigma}_*<\infty$, then $\Sigma\in L_{w^*}^\infty(I;\mathcal M(\mathcal Z_k))$ with $$\norm{\Sigma}_{L^\infty(I;\mathcal M)}=\norm{\Sigma}_*.$$
\end{lemma}

\begin{proof}
Let $\Sigma\in\mathcal M_{\loc}^+(I\times\mathcal Z_k)$ with $\norm{\Sigma}_*<\infty$.

First, show that 
$$\tproj{}_\#\Sigma\ll\mathcal L^1(I)\quad \text{and} \quad \frac{\mathrm d(\tproj{}_\#\Sigma)}{\mathrm d\mathcal L^1}\in L^\infty(I).$$
For this, take $\phi\in C_c(I)$ and velocity cutoffs $\chi_m\in C_c(\mathcal Z_k)$, $0\leq\chi_m\leq1$ increasing to $1$. Then $|\phi(t)|\chi_m(x,\tau,v)\in C_c(I\times\mathcal Z_k)$ and its $L^1(I;L^\infty)$ norm is at most $\|\phi\|_{L^1(I)}$. Hence monotone convergence gives
$$
\int_I|\phi(t)|\d(\tproj{}_\#\Sigma)(t)
=\lim_{m\to\infty}\int |\phi(t)|\chi_m(x,\tau,v)\d\Sigma
\leq\|\Sigma\|_*\|\phi\|_{L^1(I)}.
$$
Thus $\tproj{}_\#\Sigma\leq\|\Sigma\|_*\mathcal L^1$, proving both assertions.

Second, we now show that it holds 
$$\norm{\Sigma}_{L^\infty(I;\mathcal M)}=\norm{\Sigma}_*.$$
Let us take $\varphi\in C_c(I\times\mathcal Z_k), \norm{\varphi}_{L^1(I;L^\infty)}\leq 1$ and estimate as follows:
$$
\int\varphi\d\Sigma=\int_I\int \varphi(t,\cdot)\d\Sigma^t\d t
\leq\int_I \mass{\Sigma^t}\cdot\norm{\varphi(t,\cdot)}_{L^\infty}\d t\leq \norm{\Sigma}_{L^\infty(I;\mathcal M)}\cdot 1.
$$
Taking the supremum over all such $\varphi$ yields the inequality
$$\norm{\Sigma}_{L^\infty(I;\mathcal M)}\geq\norm{\Sigma}_*.$$

For the opposite inequality, suppose that $\norm{\Sigma}_{L^\infty(I;\mathcal M)}<\infty$ (otherwise there is nothing to prove), and fix $\varepsilon>0$. We have that $\left(t\mapsto \mass{\Sigma^t}\right)\in L^\infty(I)$, which implies that we can find $\phi\in C_c(I)$ with $\int_I|\phi(t)|\d t\leq 1$ and $$\norm{\Sigma}_{L^\infty(I;\mathcal M)} - \varepsilon\leq \int_I \phi(t)\mass{\Sigma^t}\d t.$$
We continue this inequality by
$$
\int_I \phi(t)\mass{\Sigma^t}\d t=
\int_I\int \phi(t)\d\Sigma^t\d t = \int\phi\circ\tproj{}\d\Sigma
\leq\norm{\Sigma}_*,
$$
where the first equality holds since the measures $\Sigma^t$ are positive, and the final inequality follows by applying the same cutoffs $\chi_m$ and passing to the limit. Since $\varepsilon>0$ was arbitrary,
we conclude that 
$$\norm{\Sigma}_{L^\infty(I;\mathcal M)}\leq\norm{\Sigma}_*.$$

Altogether, this concludes  the proof.
\end{proof}

The preceding characterization makes the mass bound stable under weak${}^*$ convergence, because it can be tested at the level of the space-time measure. The next two propositions record this lower semicontinuity and the corresponding compactness property. These facts will allow the same mass constraint to be imposed before and after passage to a limit.

\begin{proposition}\label{prop:wlsc-linfty-norm}
 If $\Sigma^i,\Sigma\in \mathcal M_{\loc}^+(I\times\mathcal Z_k)$ with $\tproj{}_\#\Sigma^i\ll\mathcal L^1(I)$ satisfy $$\Sigma^i\weakstarto\Sigma\quad \text{in }\mathcal M_{\loc}^+(I\times \mathcal Z_k), \text{ and }\liminf_{i\to\infty} \norm{\Sigma^i}_{L^\infty(I;\mathcal M)}<\infty,$$
then it holds that $\tproj{}_\#\Sigma\ll\mathcal L^1(I)$ with $$\norm{\Sigma}_{L^\infty(I;\mathcal M)}\leq\liminf_{i\to\infty} \norm{\Sigma^i}_{L^\infty(I;\mathcal M)}.$$
\end{proposition}

\begin{proof}
We utilize the $\|\cdot\|_*$ norm defined in Lemma~\ref{lem:alt-norm}. Writing out the definition of this norm and using the weak${}^*$ convergence in $\Sigma^i\weakstarto\Sigma\quad\text{in } \mathcal M_{\loc}^+(I\times\mathcal Z_k)$, we obtain
$$
\norm{\Sigma}_{*}=\sup_{\norm{\varphi}_{L^1(L^\infty)}\leq1} \int\varphi\d \Sigma =\sup_{\norm{\varphi}_{L^1(L^\infty)}\leq1}\lim_{i\to\infty}\int\varphi\d\Sigma^i.
$$
To estimate this, we switch $\sup$ and $\liminf$ with an inequality as follows: Given $(a_i^\alpha)_{i\in\N}^{\alpha\in A}$, for any index set $A$, it holds 
$$\sup_{\alpha\in A}\liminf_{i\to\infty}a_i^\alpha=\sup_{\alpha\in A}\lim_{i\to\infty}\inf_{j\geq i} a_j^\alpha\leq \sup_{\alpha\in A}\lim_{i\to\infty}\inf_{j\geq i}\sup_{\beta\in A}a_j^\beta=\liminf_{i\to\infty}\sup_{\beta\in A}a_i^\beta.$$
Using this, we estimate
$$\sup_{\norm{\varphi}_{L^1(L^\infty)}\leq1}\lim_{i\to\infty}\int\varphi\d\Sigma^i
\leq \liminf_{i\to\infty}\sup_{\norm{\varphi}_{L^1(L^\infty)}\leq1} \int\varphi\d\Sigma^i =\liminf_{i\to\infty} \norm{\Sigma^i}_{L^\infty(I;\mathcal M)}<\infty.$$

Now, by the previous lemma, we have $$\norm{\Sigma}_{L^\infty(I;\mathcal M)}=\|\Sigma\|_*,$$
which, together with the above, completes the proof.
\end{proof}

\begin{proposition}[Compactness in $L^\infty(\mathcal M)$]\label{prop:compactness-linfty-measures}
Let $\Sigma^i\in L_{w^*}^\infty(I;\mathcal M(\mathcal Z_k))$  be a sequence with $$\sup_{i\in\N}\norm{\Sigma^i}_{L^\infty(I;\mathcal M)}<\infty.$$
Then there exists a (non-relabeled) subsequence and a limit $\Sigma\in L_{w^*}^\infty(I;\mathcal M(\mathcal Z_k))$ such that $$\Sigma^i\weakstarto\Sigma \quad \text{in }L_{w^*}^\infty(I;\mathcal M(\mathcal Z_k)).$$
\end{proposition}
\begin{proof}
We can estimate $$\|\Sigma^i\|_{\mathcal M}\leq T \|\Sigma^i\|_{L^\infty(\mathcal M)}$$ so the sequence $(\Sigma^i)_{i\in\N}$ is bounded in $\mathcal M(I\times\mathcal Z_k)$. Thus, there is a subsequence such that $$\Sigma^i\weakstarto\Sigma\quad\text{in }\mathcal M(I\times\mathcal Z_k)$$
which along with the closedness above gives the desired result.
\end{proof}

\subsubsection{Slice mass bound}

We now return from the auxiliary measure setting to varicurrents and apply the preceding definition to $\Sigma_\mu$. Since its Lebesgue-time disintegration consists of the velocity-attached slices, the resulting norm controls their masses uniformly in time.

\begin{definition}[$L^\infty$ slice bound]
For any $\mu\in\varicur{1+k}^{\loc}(I\times\overline\Omega)$ with \eqref{eqn:NC-varicurrent} and $\tproj{}_\#\mu \ll\mathcal L^1(I)$ we denote $$\norm{\mu}_{\LinftySlices(I)}:=\|\Sigma_\mu\|_{L^\infty(\mathcal M)},$$
where $\Sigma_\mu\in\mathcal M^+_{\loc}(I\times\mathcal Z_k)$ is defined by \eqref{eq:Sigma-mu-defn}, or more explicitly,
$$\norm{\mu}_{\LinftySlices(I)}=\operatorname*{ess\,sup}_{t\in I} \mass{\slice{\mu}{t}}.$$
We denote the space of varicurrents with uniform slice mass bound as
$$
\LinftySlices(I):=\{\mu\in\varicur{1+k}^{\loc}(I\times\overline\Omega): \mu\text{ satisfies }\eqref{eqn:NC-varicurrent}, \tproj{}_\#\mu\ll\mathcal L^1, \|\mu\|_{\LinftySlices}<\infty\}.
$$
\end{definition}

The next estimate bounds space-time mass by slice mass and variation, it is the varicurrent analogue of \cite[Lemma 3.5]{rindlerSpacetimeIntegralCurrents2023}.

\begin{lemma}\label{lem:spacetime-mass-bdd-slice-variation}
For $
\mu\in\varicur{1+k}^{\loc}(I\times\overline\Omega)
$
satisfying \eqref{eqn:NC-varicurrent} and $\tproj{}_\#\mu\ll \mathcal L^1$ it holds for all subintervals $J\subset I$
$$
\mass{\mu\mres J}\leq  |J| \|\mu\|_{\LinftySlices(J)}+\variation{\mu;J}.
$$
\end{lemma}
\begin{proof}
Utilizing the elementary inequality $1\leq|\eta\mres \mathrm dt|+|\xproj{\eta}|$ for $\eta\in\orient$, we obtain
\begin{align*}
\mass{\mu\mres J}&\leq  \int_{J\times\overline\Omega\times\orient} |\eta \mres \mathrm dt|+|\xproj{\eta}|\d\mu = \mass{\Sigma_\mu \mres J\times \mathcal Z_k} +\variation{\mu;J}\leq |J|\|\mu\|_{\LinftySlices}+\variation{\mu;J}.\qedhere
\end{align*}
\end{proof}

\subsection{The Dissipation Potential and Lower Semicontinuity}

The velocity-attached slices allow us to evaluate the $p$-power dissipation potential against the distribution of velocities. At each time, this defines a functional of the slice. Integrating it over time evaluates it over the trajectory. This is weak${}^*$ lower semicontinuous on varicurrents, as established in Lemma~\ref{lem:wlsc-dissipation}, unlike the case of currents.

\begin{definition}[Dissipation potential]\label{def:dissipation-potential}
Fixing $p>1$, we first define the instantaneous dissipation potential evaluated on a spatial slice measure $\nu \in \mathcal M^+(\mathcal Z_k)$ as $$\mathcal{R}_p(\nu) :=\frac 1p \norm{v}_{L^p(\nu)}^p=  \frac{1}{p}\int_{\mathcal Z_k}|v|^p  \d\nu(x,\tau,v).$$
Here $\|v\|_{L^p(\nu)}^p:=\int_{\mathcal Z_k}|v|^p\d\nu$. For $\mu\in\varicur{1+k}^{\loc}(I\times\overline\Omega)$ satisfying \eqref{eqn:NC-varicurrent} and $\tproj{}_\#\mu\ll\mathcal L^1(I)$, and a measurable function $f$ on $I\times\mathcal Z_k$, we use the notation
$$
\|f\|_{L^p(J;L^p(\mu(\cdot)))}
:=\left(\int_J\int_{\mathcal Z_k}|f(t,x,\tau,v)|^p
\d\slice{\mu}{t}(x,\tau,v)\d t\right)^{1/p}.
$$
In particular, we have the identities
$$
\|v\|_{L^p(\mu(t))}^p=p\mathcal R_p(\slice{\mu}{t}),
\qquad
\|v\|_{L^p(J;L^p(\mu(\cdot)))}^p
=p\int_J\mathcal R_p(\slice{\mu}{t})\d t.
$$
\end{definition}

\begin{remark}[Representation of integrated dissipation potential]
For $\mu$ satisfying \eqref{eqn:NC-varicurrent} and $\tproj{}_\#\mu\ll\mathcal L^1(I)$, the definition of the slices gives
\begin{equation}\label{eqn:dissipation-potential-space-time}
\int_J\mathcal R_p(\slice{\mu}{t})\d t
=\frac1p\|v\|_{L^p(J;L^p(\mu(\cdot)))}^p =\int_{J\times\overline\Omega\times\orient}\frac1p|v(\eta)|^p|\eta \mres \mathrm dt|\d\mu(t,x,\eta).
\end{equation}
For elementary varicurrents associated with $S\in\intcur{1+k}(I\times\R^d)$ satisfying \eqref{eqn:NC-current}, this agrees with the current formula:
$$
\mathcal R_p(\slice{\eym{S}}{t})=\mathcal R_p\left(S(t),\frac{\mathrm D}{\mathrm Dt}S(t,\cdot)\right)
$$
where the expression on the right is understood as in \eqref{eqn:intro-Rp-current-defn}.
\end{remark}

\begin{example}[Trajectory of a point]\label{ex:point-trajectory}
A moving point $x\colon[0,\infty)\to\overline{\Omega}$ with
$$
x\in W^{1,p}_{\loc}([0,\infty);\overline\Omega)
$$
is represented by the current
$$
S=(t\mapsto(t,x(t)))_*\llbracket(0,\infty)\rrbracket \in \intcur{1}^{\loc}([0,\infty)\times\overline\Omega),
$$
with slices
$$
S(t)= \llbracket x(t) \rrbracket = \delta_{x(t)} \in\intcur{0}(\overline\Omega),\quad t\in[0,\infty).
$$
The associated elementary varicurrent is
$$
\mu=\eym{S}=\underbrace{\mathcal L^1(\mathrm dt)\otimes\sqrt{1+|\dot x(t)|^2}\cdot \delta_{x(t)}(\mathrm dx)}_{\|S\|(\mathrm dt,\mathrm dx)}\otimes \delta_{\underbrace{\frac{e_0+\dot x(t)}{\sqrt{1+|\dot x(t)|^2}}}_{=\vec S(t,x)}}(\mathrm d\eta)\in\varicur{1}^{\loc}([0,\infty)\times\overline\Omega).
$$
This elementary varicurrent satisfies \eqref{eqn:NC-varicurrent} and $\tproj{}_\#\|S\|\ll \mathcal L^1$. The dissipation potential of the slice is then
$$
\mathcal R_p(\slice{\mu}{t})=\tfrac1p|\dot x(t)|^p\quad \text{for a.e. }t\in(0,\infty),
$$
and in the integrated form
$$
\int_A\mathcal R_p(\slice{\mu}{t})\d t=\int_A\tfrac1p|\dot x(t)|^p\d t.
$$
\end{example}

The following lower-semicontinuity result supplies the property missing in the introductory current examples. It allows us to pass to the limit in the integrated dissipation term when constructing WED minimizers and their limiting evolutions in Section~\ref{sec:WED}.

\begin{lemma}[Lower semicontinuity of the dissipation potential]\label{lem:wlsc-dissipation}
For any $(\mu_n)\subset\varicur{1+k}(I\times\overline\Omega)$, $\mu\in\varicur{1+k}(I\times\overline\Omega)$ satisfying \eqref{eqn:NC-varicurrent} and $\tproj{}_\#\mu_n,\tproj{}_\#\mu\ll \mathcal L^1$  with
$$\mu_n\weakstarto\mu ,$$
it holds for all subintervals $J\subset I$ with $\partial J\cap N=\emptyset$, where $N$ is the countable set $N=\left\{a\in I:\mu(\{a\}\times\overline\Omega\times\orient)>0\right\}$, that
$$ \int_J \mathcal R_p(\slice{\mu}{t})\d t \leq \liminf_{n\to\infty}\int_J \mathcal R_p(\slice{\mu_n}{t})\d t.$$
\end{lemma}
\begin{proof}
Pick an interval $J\subset I$ such that $\partial J\cap N =\emptyset$. As $\mu$ is locally finite in time, $N$ is countable. By \cite[Proposition 1.62]{ambrosioFunctionsBoundedVariation2000}, we also have that
$$\mu_n\mres J\weakstarto\mu\mres J.$$
By Lemmas~\ref{lem:velocity-cont-lsc} and \ref{lem:v-tau-p-x-identities}, the integrand
$$
\eta \mapsto|v(\eta)|^p |\eta \mres \mathrm dt|
$$
extended by $+\infty$ on $\orientSpatial$, is lower semicontinuous on $\orient$ because $p>1$. Thus, by the Portmanteau theorem (see e.g.\ \cite{bogachevMeasureTheory2007})
\begin{align*}
 \int_J \mathcal R_p(\slice{\mu}{t})\d t &=\int_{J\times\overline\Omega\times\orient}\frac1p|v(\eta)|^p |\eta \mres \mathrm dt|\d\mu(t,x,\eta)\\
&\leq
\liminf_{n\to\infty} \int_{J\times\overline\Omega\times\orient}\frac1p|v(\eta)|^p |\eta \mres \mathrm dt|\d\mu_n(t,x,\eta)\\
&=
\liminf_{n\to\infty}\int_J \mathcal R_p(\slice{\mu_n}{t})\d t ,
\end{align*}
finishing the proof.
\end{proof}

Together with the slice-mass bound, finite dissipation controls temporal variation. The Thousand Cuts example (Example~\ref{ex:thousand-cuts-current}) shows why the mass bound cannot simply be dropped.

\begin{lemma}[Continuity in time]\label{lem:p-continuity-in-time}
Let $\mu\in\varicur{1+k}([a,b]\times\overline\Omega)$ satisfy \eqref{eqn:NC-varicurrent} and
$$
\int_a^b\mathcal R_p(\slice{\mu}{t})\d t<\infty, \quad \tproj{}_\#\mu\ll\mathcal L^1([a,b]), \quad \|\mu\|_{\LinftySlices((a,b))} <\infty.
$$
Then it holds 
$$
\variation{\mu;[a,\cdot]}\in C^{0,1-1/p}([a,b])
$$
with
\begin{align*}
[\variation{\mu;[a,\cdot]}]_{C^{0,1-1/p}}\leq
\|v\|_{L^p((a,b);L^p(\mu(\cdot)))}
\|\mu\|_{\LinftySlices((a,b))}^{1-1/p}.
\end{align*}
\end{lemma}
\begin{proof}
We compute for $a\leq t_1<t_2\leq b$ by using \eqref{eqn:NC-varicurrent}, the slice integral \eqref{eq:weighted-phase-space-measure-test}, \eqref{eq:definition-velocity-attached-slices}, and the H\"older inequality:
\begin{align*}
\variation{\mu;[t_1,t_2]}
&=
\int_{[t_1,t_2]\times\overline\Omega\times\orient}|\xproj\eta|\d\mu(t,x,\eta)
=
\int_{[t_1,t_2]\times\overline\Omega\times\orient}|v(\eta)||\eta \mres \mathrm dt|\d\mu(t,x,\eta)
\\&=
\int_{t_1}^{t_2}\int_{\mathcal Z_k}|v| \d\slice{\mu}{t}(x,\tau,v)\d t
\\&\leq
\left(\int_{t_1}^{t_2}\int_{\mathcal Z_k}|v|^p\d\slice{\mu}{t}(x,\tau,v)\d t\right)^{1/p}
\left(\int_{t_1}^{t_2}\int_{\mathcal Z_k}1\d\slice{\mu}{t}(x,\tau,v)\d t\right)^{1-1/p}
\\&\leq
\|v\|_{L^p((a,b);L^p(\mu(\cdot)))}
\|\mu\|_{\LinftySlices((a,b))}^{1-1/p}|t_2-t_1|^{1-1/p}.
\end{align*}
This completes the proof.
\end{proof}

\begin{lemma}[Weak${}^*$ continuity of variation]\label{lem:lower-semicontinuity-variation}
Let $\mu_n\in\varicur{1+k}^{\loc}(I\times\overline\Omega)$ and $\mu\in\varicur{1+k}^{\loc}(I\times\overline\Omega)$ with
$$
\mu_n\weakstarto\mu.
$$
Then for all compact intervals $J\subset  I$ whose endpoints avoid a countable set it holds
$$
\variation{\mu;J}=\lim_{n\to\infty}\variation{\mu_n;J}.
$$
\end{lemma}
\begin{proof}
By local finiteness of $\mu$, we have $\mu(\partial J\times\overline\Omega\times\orient)=0$ provided the endpoints of $J$ avoid a countable set of times where $\mu$ has positive mass. For such $J$,
$$
\mu_n\mres J\weakstarto\mu\mres J.
$$
Consequently, as the integrand $\eta\mapsto |\xproj{\eta}|$ is continuous and bounded, and the interval $J$ is compact, we have
$$
\variation{\mu;J}=\int_{J\times\overline\Omega\times\orient}|\xproj{\eta}|\d\mu(t,x,\eta)=\lim_{n\to\infty}\int_{J\times\overline\Omega\times\orient}|\xproj{\eta}|\d\mu_n(t,x,\eta)= \lim_{n\to\infty}\variation{\mu_n;J},
$$
concluding the proof.
\end{proof}

We can now combine the mass and variation estimates with lower semicontinuity of the dissipation. Bounds on slice mass and dissipation give compactness of trajectories, while the following result also ensures that the limiting measure satisfies the hypotheses needed for slicing. This stability is essential when these trajectories are used as competitors in the variational construction.
\begin{proposition}[Compactness of the dissipation potential]\label{prop:compactness-dissipation-potential}
Let $(\mu_i)_{i\in\N}\subset\varicur{1+k}(I\times\overline\Omega)$ satisfy \eqref{eqn:NC-varicurrent},
$\tproj{}_\#\mu_i\ll\mathcal L^1$, and
\begin{align*}
    \sup_{i\in\N}\|\mu_i\|_{\LinftySlices(I)}&<\infty,
    \\
    \sup_{i\in\N}\int_J \mathcal R_p(\slice{\mu_i}{t})\d t&<\infty\quad \text{for all compact intervals }J\subset I.
    \\
   \partial \stcurr{\mu_i} \mres \interior{I} \times \R^d &= 0.
\end{align*}
Then there exists a limit $\mu\in\varicur{1+k}(I\times\overline\Omega)$ such that, for a non-relabeled subsequence it holds
$$
\mu_i\weakstarto\mu,
$$
where the limit $\mu$ satisfies \eqref{eqn:NC-varicurrent},
$\tproj{}_\#\mu\ll\mathcal L^1$, and
$\partial\stcurr{\mu}\mres(\interior{I}\times\R^d)=0$. Moreover, for every compact interval $J\subset I$,
$$
\|\mu\|_{\LinftySlices(J)} \leq\liminf_{i\to\infty}\|\mu_i\|_{\LinftySlices(J)} 
$$
and 
$$
\int_J\mathcal R_p(\slice{\mu}{t})\d t\leq \liminf_{i\to\infty}\int_J\mathcal R_p(\slice{\mu_i}{t})\d t .
$$
\end{proposition}
\begin{proof}
By Lemmas \ref{lem:spacetime-mass-bdd-slice-variation}
and \ref{lem:p-continuity-in-time},
\[
\mass{\mu_i\mres J}
\leq
|J|\|\mu_i\|_{\LinftySlices(I)}
+
\|\mu_i\|_{\LinftySlices(I)}^{1-1/p}|J|^{1-1/p}
\left(p\int_J\mathcal R_p(\slice{\mu_i}{t})\,\d t\right)^{1/p}.
\]
Thus Lemma \ref{lem:compactness-varicurrents} gives,
along a subsequence,
$\mu_i\weakstarto\mu$ for some
$\mu\in\varicur{1+k}^{\loc}(I\times\overline\Omega)$.

By Lemma \ref{lem:velocity-cont-lsc}, $|v|$ is continuous
with values in $[0,\infty]$. Since $p>1$ and
\[
|v(\eta)|^p|\eta \mres \mathrm dt|
=
\frac{|v(\eta)|^p}{\sqrt{1+|v(\eta)|^2}}
\qquad\text{on }\orient\setminus\orientSpatial,
\]
the dissipation density is lower semicontinuous, with value
$+\infty$ on $\orientSpatial$.
The Portmanteau theorem, applied with compactly supported
time cutoffs, and the uniform local dissipation bounds
imply that this density is locally $\mu$-integrable.
In particular, $\mu$ satisfies \eqref{eqn:NC-varicurrent}.

Furthermore,
\[
\tproj{}_\#\bigl(|\eta \mres \mathrm dt|\mu_i\bigr)
\leq
\left(\sup_j\|\mu_j\|_{\LinftySlices(I)}\right)\mathcal L^1.
\]
This inequality passes to the limit. Since
$|\eta \mres \mathrm dt|>0$ $\mu$-a.e., it follows that
$\tproj{}_\#\mu\ll\mathcal L^1$, and consequently, the slices of $\mu$ are defined.

For every compact interval $J\subset I$, time absolute
continuity gives $\mu_i\mres J\weakstarto\mu\mres J$.
Passing to the limit in
\[
\tproj{}_\#\bigl(|\eta \mres \mathrm dt|(\mu_i\mres J)\bigr)
\leq
\|\mu_i\|_{\LinftySlices(J)}\mathcal L^1\mres J
\]
along a subsequence realizing the liminf yields
\[
\|\mu\|_{\LinftySlices(J)}
\leq\liminf_{i\to\infty}\|\mu_i\|_{\LinftySlices(J)}.
\]
Finally, Lemma \ref{lem:wlsc-dissipation} gives the $\mathcal R_p$-lower semicontinuity on every compact interval.
\end{proof}

\subsection{Time pushforward}

Changing the time parametrization of a varicurrent must transform both its space-time position and the orientation, as its steepness in the time direction changes. An ordinary measure-pushforward does not account for the change in the tangent planes. We therefore extend the time-pushforward construction for currents, cf.~\cite[Lemma 3.4]{rindlerEnergeticSolutionsRateindependent2024}, to the case of varicurrents.

The tangent $(1+k)$-vector needs to stay unit so its correct transformation is to tilt it in the time direction.  A time dilation transforms the direction correctly, but it changes the size of an orienting $(1+k)$-vector. We thus achieve the desired tilt by correcting for the $(1+k)$-area factor after time-dilation. This is explained in what follows.

\begin{definition}[Tilt of orienting $(1+k)$-vectors]
For $\lambda >0$ we define the \emph{time-dilation}  $D_\lambda\colon\R^{1+d}\to\R^{1+d}$ by 
$$D_\lambda(t,x)=(\lambda t,x), \quad t\in\R,x\in\R^d,$$
and denote the induced dilation on $(1+k)$-vectors by $\Wedge^{1+k} D_\lambda\colon\Wedge^{1+k}\R^{1+d}\to\Wedge^{1+k}\R^{1+d}$.

We define the corresponding \emph{tilt transformation} of orienting $(1+k)$-vectors $\lambda_\sim\colon\orient\to\orient$ by
$$
\lambda_\sim \eta := \frac{1}{J_\lambda(\eta)}\Wedge^{1+k}D_\lambda\eta, \quad \text{where } J_\lambda(\eta)=\sqrt{\lambda^2|\eta \mres \mathrm dt|^2+|\xproj{\eta}|^2},\quad \eta\in\orient.
$$
\end{definition}
The next lemma identifies $J_\lambda(\eta)$ as the $(1+k)$-area factor for a time dilation by $\lambda$. This is the weight that accompanies the normalized orientation in the pushforward construction.
\begin{lemma}\label{lem:Jlambda-is-jacobian}
For every $\lambda>0$ and $\eta\in\orient$, the factor $J_\lambda(\eta)$ is the tangential Jacobian of $D_\lambda$ on the oriented $(1+k)$-plane represented by $\eta$. Moreover, it holds
\[
J_{\lambda^{-1}}(\lambda_\sim\eta)=\frac{1}{J_\lambda(\eta)}
\qquad\text{and}\qquad
(\lambda^{-1})_\sim=(\lambda_\sim)^{-1}.
\]
\end{lemma}
\begin{proof}
As $|\eta|=1$, we immediately see that 
$$
J_\lambda(\eta)=\frac{|(\Wedge^{1+k}D_\lambda)\eta|}{|\eta|}
$$
which is the definition of the tangential Jacobian of $D_\lambda$ on the plane spanned by $\eta$.

To get the formulas for the inverse, we compute
$$
J_{\lambda^{-1}}(\lambda_\sim\eta)=\sqrt{\lambda^{-2}|\lambda_\sim\eta\mres \mathrm dt|^2+|\xproj{\lambda_\sim\eta}|^2}=\sqrt{\lambda^{-2}\left|\frac{\lambda \eta \mres \mathrm dt}{J_\lambda(\eta)}\right|^2+\left|\frac{\xproj{\eta}}{J_\lambda(\eta)}\right|^2}=\frac{1}{J_\lambda(\eta)},
$$
where we used $\sqrt{|\eta \mres \mathrm dt|^2+|\xproj{\eta}|^2}=|\eta|=1$, and consequently
\begin{align*}
(\lambda^{-1})_{\sim}(\lambda_\sim \eta)
&=  \frac{1}{J_{\lambda^{-1}}(\lambda_\sim\eta)} \frac{1}{J_\lambda(\eta)}D_{\lambda^{-1}}D_\lambda\eta =\eta. \qedhere
\end{align*}
\end{proof}

Having identified the area factor, we record that the tilt preserves the spatial slice orientation and scales its velocity.

\begin{lemma}[Slice orientation and velocity under tilting]
\label{lem:tilt-slice-orientation-velocity}
Let $\lambda>0$ and let
$\eta\in\orient\setminus\orientSpatial$. Then
$\lambda_\sim\eta\in\orient\setminus\orientSpatial$ and, using the decomposition of Lemma \ref{lem:space-time-decomposition} we have
\begin{equation}\label{eq:tilt-tau-velocity}
\tau(\lambda_\sim\eta)=\tau(\eta),
\qquad
v(\lambda_\sim\eta)=\frac{v(\eta)}{\lambda}.
\end{equation}
\end{lemma}
\begin{proof}
Since $\tau(\eta)$ has no temporal component, we have $D_\lambda\tau(\eta)=\tau(\eta)$. Consequently,
$$
\lambda_\sim\eta= \frac{1}{J_\lambda(\eta)} D_\lambda\xi(\eta)\wedge D_\lambda\tau(\eta) = \frac{1}{J_\lambda(\eta)}D_\lambda\xi(\eta)\wedge \tau(\eta)
$$
from which it follows that
$$
\xi(\lambda_\sim\eta)=\frac{1}{J_\lambda(\eta)}D_\lambda\xi(\eta)\quad\text{and}\quad \tau(\lambda_\sim\eta)=\tau(\eta).
$$
From this we obtain
$$
v(\lambda_\sim\eta)=\frac{\xproj{\xi(\lambda_\sim\eta)}}{\tproj{\xi(\lambda_\sim\eta)}}=\frac{\xproj{\xi(\eta)}}{\lambda\tproj{\xi(\eta)}}=\frac{v(\eta)}{\lambda},
$$
which finishes the proof.
\end{proof}

Using this, we define time pushforward of a varicurrent. The local dilation factor is given by the derivative of the time map.

\begin{definition}[Time pushforward]\label{def:time-pushforward}
Let $\mu\in\varicur{1+k}^{\loc}(I \times\overline\Omega)$ and $\varphi\in C^1(I)$, $\dot \varphi>0$. Then we define $\timepush{\varphi}{\mu}$ by setting
\begin{equation}\label{eqn:time-pushforward}
\int g(s,x,\xi)\d\timepush{\varphi}{\mu}(s,x,\xi):=\int g\left(\varphi(t),x,\dot \varphi(t)_\sim\eta\right)J_{\dot\varphi(t)}(\eta)\d\mu(t,x,\eta)
\end{equation}
for all $g\in C_c(\varphi(I)\times\overline\Omega\times\orient)$.
\end{definition}
The following lemma verifies that the factor $J_{\dot\varphi(t)}(\eta)$ in \eqref{eqn:time-pushforward} is exactly the area factor associated with the pushforward by $\overline\varphi$.
\begin{lemma}
\label{lem:time-reparametrization-area-factor}
Let $\varphi\in C^1([a,b])$ satisfy $\dot\varphi>0$, and define
\[
\overline\varphi(t,x):=(\varphi(t),x).
\]
Then $$
\nabla\overline\varphi(t,x)=D_{\dot\varphi(t)},
$$
and the tangential  Jacobian of $\overline\varphi$ on the oriented plane $\eta$ satisfies
\[
|\nabla^\eta\overline\varphi(t,x)|
=J_{\dot\varphi(t)}(\eta).
\]
Moreover, the induced orientation is
\[
\frac{\nabla\overline\varphi(t,x)\eta}
     {|\nabla\overline\varphi(t,x)\eta|}
=\dot\varphi(t)_\sim\eta.
\]
\end{lemma}
\begin{proof}
The form of $\nabla \overline\varphi$ is immediate and the rest follows from Lemma \ref{lem:Jlambda-is-jacobian}.
\end{proof}
From this we obtain the formula for the time-pushforward of a space-time current, cf. \cite[Lemma 3.4]{rindlerSpacetimeIntegralCurrents2023}.
\begin{corollary}[Time-pushforward of space-time current]\label{cor:pushforward-current}
Let $\varphi\in C^1(I)$ with $\dot\varphi>0$. Then for $S\in\intcur{1+k}^{\loc}(I\times\overline\Omega)$ it holds
\begin{align*}
\langle \timepush{\varphi}{S},\omega\rangle&=\int \langle \dot\varphi(t)_\sim \vec S(t,x),\omega(\varphi(t),x)\rangle J_{\dot\varphi(t)}(\vec S(t,x)) \d\| S\|(t,x)\quad \text{for all }\omega\in \mathcal D^{1+k}(\varphi(I)\times\overline\Omega).
\end{align*}
\end{corollary}
\begin{proof}
As $\timepush{\varphi}{S}$ denotes the pushforward of $S$ by $\overline\varphi$, we have by the definition of current pushforward
$$
\langle \timepush{\varphi}{S},\omega\rangle = \int \langle [\nabla\overline\varphi(t,x)]\vec S(t,x),\omega(\varphi(t),x)\rangle\d\|S\|(t,x).
$$
Utilizing the identities of Lemma \ref{lem:time-reparametrization-area-factor} completes the proof.
\end{proof}

The next proposition describes how time pushforward transforms the slices, their orientations, and velocities.

\begin{proposition}[Time pushforward and slicing]\label{prop:time-pushforward}
Let $\mu\in\varicur{1+k}^{\loc}(I\times\overline\Omega)$ satisfy $\tproj{}_\#\mu\ll\mathcal L^1(I)$ and \eqref{eqn:NC-varicurrent}, and let $\varphi\in C^1(I)$, $\dot\varphi>0$. Then $\timepush{\varphi}{\mu}\in\varicur{1+k}^{\loc}(\varphi(I)\times\overline\Omega)$ satisfies $\tproj{}_\#(\timepush{\varphi}{\mu})\ll\mathcal L^1(\varphi(I))$ and \eqref{eqn:NC-varicurrent}. The velocity-attached time slices of $\timepush{\varphi}{\mu}$ are given by
$$
\slice{\timepush{\varphi}{\mu}}{\varphi(t)}=\left((x,\tau,v)\mapsto\left(x,\tau,\frac{v}{\dot\varphi(t)}\right)\right)_\#\slice{\mu}{t} \quad \text{for a.e. }t\in I,
$$
equivalently, it holds
$$
\int g(x,\tau,v)\d\slice{\timepush{\varphi}{\mu}}{\varphi(t)}(x,\tau,v)
=\int g\left(x,\tau,\frac{v}{\dot\varphi(t)}\right)\d\slice{\mu}{t}(x,\tau,v)
$$
for a.e.\ $t\in I$ and every $g\in C_b(\mathcal Z_k)$.
Consequently, it holds 
$$
\mass{\slice{\timepush{\varphi}{\mu}}{\varphi(t)}}=\mass{\slice{\mu}{t}}, \quad \text{a.e. }t\in I
$$
and
$$
\slicenov{\timepush{\varphi}{\mu}}{\varphi(t)}=\slicenov{\mu}{t},\quad \text{a.e. }t\in I.
$$
Moreover, $\mu\mapsto \varphi_*\mu\colon \varicur{1+k}^{\loc}(I\times\overline\Omega)\to\varicur{1+k}^{\loc}(\varphi(I)\times\overline\Omega)$ is a bijection with
$$
\timepush{(\varphi^{-1})}{(\timepush{\varphi}{\mu})}=\mu.
$$
\end{proposition}
\begin{proof}
The tilt preserves noncritical orientations by Lemma \ref{lem:tilt-slice-orientation-velocity}, so the pushforward preserves \eqref{eqn:NC-varicurrent}. As $\varphi$ is $C^1$, the property $\tproj{}
_\#(\timepush{\varphi}{\mu})\ll\mathcal L^1$ follows. The barycenter is the current pushforward (as verified below), and hence is locally integral. We may now use the slice disintegration. Utilizing its definition and Lemma \ref{lem:tilt-slice-orientation-velocity}, we get for $h\in C_c(\varphi(I)\times\mathcal Z_k)$
\begin{align*}
&\int_{\varphi(I)}\int_{\mathcal Z_k} h(s,x,\tau,v)\d\slice{\timepush{\varphi}{\mu}}{s}(x,\tau,v)\d s\\
&\quad=\int h(s,x,\tau(\zeta),v(\zeta))|\zeta \mres\mathrm dt| \d\timepush{\varphi}{\mu}(s,x,\zeta)
\\&\quad=
\int h (\varphi(t),x,\tau(\dot\varphi(t)_\sim\eta),v(\dot\varphi(t)_\sim\eta))|\tproj{(\dot\varphi(t)_\sim\eta)}|J_{\dot\varphi(t)}(\eta)\d\mu(t,x,\eta)
\\&\quad=
\int h \left(\varphi(t),x,\tau(\eta),\frac{v(\eta)}{\dot\varphi(t)}\right)\dot\varphi(t)|\eta \mres \mathrm dt|\d\mu(t,x,\eta)
\\
&\quad=\int_I \dot\varphi(t)\int_{\mathcal Z_k} h\left(\varphi(t),x,\tau,\frac{v}{\dot\varphi(t)}\right)\d\slice{\mu}{t}(x,\tau,v)\d t.
\end{align*}
Changing variables $s=\varphi(t)$, the last expression transforms as
\begin{align*}
&\int_I \dot\varphi(t)\int_{\mathcal Z_k} h\left(\varphi(t),x,\tau,\frac{v}{\dot\varphi(t)}\right)\d\slice{\mu}{t}(x,\tau,v)\d t\\
&\quad=\int_{\varphi(I)} \int_{\mathcal Z_k} h\left(s,x,\tau,\frac{v}{\dot\varphi(\varphi^{-1}(s))}\right)\d\slice{\mu}{\varphi^{-1}(s)}(x,\tau,v)\d s
\end{align*}
The statement now follows by taking all $h$ of the form $h(s,x,\tau,v)=\psi(s)g(x,\tau,v )$, $\psi\in C_c(\varphi(I))$, $g\in C_c(\mathcal Z_k)$, or equivalently by the uniqueness of Lebesgue disintegration in time. The claims about slice mass and velocity-free slices are a direct consequence.

To verify the formula for the inverse, take any $g\in C_c(I\times\overline\Omega\times\orient)$ and compute
\begin{align*}
\langle \timepush{(\varphi^{-1})}{(\timepush{\varphi}{\mu})},g\rangle 
&=
\int g(t,x,\eta)\d\timepush{(\varphi^{-1})}{(\timepush{\varphi}{\mu})}(t,x,\eta)
\\
&=
\int g(\varphi^{-1}(s),x,\dot\varphi^{-1}(s)_\sim\zeta)\d\timepush{\varphi}{\mu}(s,x,\zeta)
\\
&=\int g(\varphi^{-1}(\varphi(t)),x,\dot\varphi(t)_\sim\dot\varphi(t)^{-1}\eta)\d\mu(t,x,\eta)
\\
&=\int g(t,x,\eta)\d\mu(t,x,\eta)=\langle \mu,g\rangle. \qedhere
\end{align*}
\end{proof}

We now verify that time-pushforward commutes with taking the barycenter current.

\begin{proposition}[Pushforward and barycenter current]
Let $\mu\in\varicur{1+k}^{\loc}(I\times\overline\Omega)$ and $\varphi\in C^1(I)$, $\dot \varphi>0$. Then it holds
$$
\stcurr{\timepush{\varphi}{\mu}}=\timepush{\varphi}{\stcurr{\mu}}
$$
and consequently,
$$
\slice{\stcurr{\timepush{\varphi}{\mu}}}{\varphi(t)}=\slice{\stcurr{\mu}}{t}, \quad \text{for a.e. }t\in I.
$$
\end{proposition}
\begin{proof}
Take $\omega\in \mathcal D^{1+k}(\varphi(I)\times\overline\Omega)$ and compute
\begin{align*}
\langle \stcurr{\timepush{\varphi}{\mu}},\omega\rangle 
&=
\int\langle \zeta,\omega(s,x)\rangle \d\timepush{\varphi}{\mu}(s,x,\zeta)
\\
&=\int \langle \dot \varphi(t)_\sim\eta,\omega(\varphi(t),x)\rangle J_{\dot\varphi(t)}(\eta)\d\mu(t,x,\eta)
\\
&=\int\langle D_{\dot\varphi(t)}\eta,\omega(\varphi(t),x)\rangle \d\mu(t,x,\eta)
\\
&=\int\langle \eta,D_{\dot\varphi(t)}\omega(\varphi(t),x)\rangle \d\mu(t,x,\eta)
\\
&=\int\langle \vstcurr{\mu}(t,x),D_{\dot\varphi(t)}\omega(\varphi(t),x)\rangle \d\|\stcurr{\mu}\|(t,x)
\\
&=\int\langle D_{\dot\varphi(t)}\vstcurr{\mu}(t,x),\omega(\varphi(t),x)\rangle \d\|\stcurr{\mu}\|(t,x)
\\
&=\int\langle \overrightarrow{\timepush{\varphi}{\stcurr{\mu}}}(s,x),\omega(s,x)\rangle \d\|\timepush{\varphi}{\stcurr{\mu}}\|(s,x)= \langle \timepush{\varphi}{\stcurr{\mu}},\omega\rangle,
\end{align*}
where in the second-to-last equality we used Corollary \ref{cor:pushforward-current}.
\end{proof}

We will also need to traverse a trajectory in the opposite time direction while preserving the orientations of its spatial slices.

\begin{lemma}[Time reversal]\label{lem:time-reversal} 
Let $\mu\in\varicur{1+k}([0,a]\times\overline\Omega)$ satisfy
\eqref{eqn:NC-varicurrent} and
$\tproj{}_\#\mu\ll\mathcal L^1$, and define $D_{-1}\colon\R^{1+d}\to\R^{1+d}$ by $D_{-1}(t,x)=(-t,x)$.
Define the \emph{time reversal} of $\mu$ by
$$
\mu^{\mathrm{rev}}
:=
\left((t,x,\eta)\mapsto
\left(a-t,x,-(\Wedge^{1+k}D_{-1})\eta\right)\right)_\#\mu.
$$
Then $\mu^{\mathrm{rev}}\in
\varicur{1+k}([0,a]\times\overline\Omega)$ satisfies
\eqref{eqn:NC-varicurrent} and
$\tproj{}_\#\mu^{\mathrm{rev}}\ll\mathcal L^1$, and its velocity-attached slices are given by
\begin{equation}\label{eq:time-reversal-slices}
\slice{\mu^{\mathrm{rev}}}{a-t}
=
\left((x,\tau,v)\mapsto(x,\tau,-v)\right)_\#
\slice{\mu}{t}
\quad\text{for a.e. }t\in [0,a].
\end{equation}
In particular, the velocity-free slices and current slices are
traversed in the opposite order, with their orientations preserved, that is,
$$
\slicenov{\mu^{\mathrm{rev}}}{a-t}
=\slicenov{\mu}{t},
\qquad
\stcurr{\mu^{\mathrm{rev}}}(a-t)
=\stcurr{\mu}(t)
\quad\text{for a.e. }t\in [0,a].
$$
\end{lemma}

\begin{proof}
The map $D_{-1}$ is orthogonal, so
$-(\Wedge^{1+k}D_{-1})\eta\in\orient$ whenever
$\eta\in\orient$, and thus $\mu^{\mathrm{rev}}\in \mathcal M^+([0,a]\times\overline\Omega\times\orient)$.
The definition of the barycenter gives
$$
\stcurr{\mu^{\mathrm{rev}}}
=-((t,x)\mapsto (a-t,x))_*\stcurr{\mu}.
$$
In particular, this is an integral current, so
$\mu^{\mathrm{rev}}\in\varicur{1+k}([0,a]\times\overline\Omega)$.

For $\eta\in\orient\setminus\orientSpatial$, the decomposition in
Lemma \ref{lem:space-time-decomposition} gives
$$
-(\Wedge^{1+k}D_{-1})\eta
=(-D_{-1}\xi(\eta))\wedge\tau(\eta).
$$
It holds  $\tproj{(-D_{-1}\xi(\eta))}>0$, thus $-(\Wedge^{1+k}D_{-1})\eta \in\orient\setminus\orientSpatial$, and by uniqueness of the decomposition we therefore have
$$
\tau\left(-(\Wedge^{1+k}D_{-1})\eta\right)=\tau(\eta),
\qquad
v\left(-(\Wedge^{1+k}D_{-1})\eta\right)=-v(\eta).
$$

The reflection clearly preserves \eqref{eqn:NC-varicurrent} and $\tproj{}_\#\mu^{\mathrm{rev}}\ll \mathcal L^1$. For any $g\in C_c([0,a]\times\mathcal Z_k)$, the definition of the
velocity-attached slices consequently gives
\begin{align*}
&\int_0^a\int_{\mathcal Z_k}
g(s,x,\tau,v)\d\slice{\mu^{\mathrm{rev}}}{s}(x,\tau,v)\d s
=
\int_{[0,a]\times\overline\Omega\times\orient}
g(s,x,\tau(\zeta),v(\zeta))
|\zeta \mres\mathrm dt|\d\mu^{\mathrm{rev}}(s,x,\zeta)
\\
&\quad=
\int_{[0,a]\times\overline\Omega\times\orient}
g(a-t,x,\tau(\eta),-v(\eta))
|\eta \mres \mathrm dt|\d\mu(t,x,\eta)\\
&\quad=
\int_0^a\int_{\mathcal Z_k}
g(a-t,x,\tau,-v)\d\slice{\mu}{t}(x,\tau,v)\d t
\\
&\quad=
\int_0^a\int_{\mathcal Z_k}
g(s,x,\tau,-v)\d\slice{\mu}{a-s}(x,\tau,v)\d s.
\end{align*}
By uniqueness of Lebesgue disintegration, this shows \eqref{eq:time-reversal-slices}.
The assertions about velocity-free slices and current slices follow immediately.
\end{proof}

\begin{remark}[Reversal versus time pushforward]
The extra minus sign in the definition of $\mu^{\mathrm{rev}}$ preserves the spatial slice orientations while reversing the  time order in which those slices are traversed. This differs from geometric pushforward by a decreasing time map. If we extended the time-pushforward construction to decreasing time maps, pushforward under the reflection $t\mapsto a-t$ would instead transform the orientation by
$$
\eta\mapsto(\Wedge^{1+k}D_{-1})\eta,
$$
without the extra minus sign. Its current slices would then have the opposite orientation, cf.~\cite[Proposition~3.6]{bonicattoSpacetimeConnectivityTheorem2026}. We thus distinguish this operation, called \emph{reversal}, from \emph{time pushforward}.
\end{remark}

Returning to increasing time changes, we now record their effect on the dissipation potential. The scaling law below will be used to compare differently parametrized trajectories in the variational arguments.

\begin{lemma}[Dissipation of a pushforward]\label{lem:dissipation-pushforward}
Let $\mu\in\varicur{1+k}^{\loc}(I\times\overline\Omega)$ satisfy \eqref{eqn:NC-varicurrent} and $\tproj{}_\#\mu\ll\mathcal L^1(I)$, and let $\varphi\in C^1(I)$ with $\dot\varphi>0$. Then
$$
\mathcal R_p(\slice{\timepush{\varphi}{\mu}}{\varphi(t)})
=\dot\varphi(t)^{-p}\mathcal R_p(\slice{\mu}{t})\quad\text{for a.e. }t\in I,
$$
and, for every Borel set $A\subset I$,
$$
\int_{\varphi(A)}\mathcal R_p(\slice{\timepush{\varphi}{\mu}}{s})\d s
=\int_A\dot\varphi(t)^{1-p}\mathcal R_p(\slice{\mu}{t})\d t.
$$
\end{lemma}
\begin{proof}
Proposition \ref{prop:time-pushforward} transforms each slice by $v\mapsto v/\dot\varphi(t)$, with no factor multiplying the slice measure. The first identity follows from the $p$-homogeneity of $\mathcal R_p$. Integration and the change of variables $s=\varphi(t)$ give the second. Positive $C^1$ time changes thus preserve local integrated-dissipation bounds. Together with Proposition \ref{prop:time-pushforward}, they preserve NC, time absolute continuity, and the slice-mass bound. Current boundaries and continuous states transform by the corresponding current pushforward; in particular, time changes fixing zero preserve the initial state.
\end{proof}

\subsection{Revisiting the examples}\label{ssec:examples-revisited}

We revisit the examples of Section~\ref{ssec:examples-currents}, interpreting them as varicurrents. To include their singular-time limits, we extend the notation \eqref{eq:example-current-cost} to varicurrents by
\begin{equation}\label{eq:example-varicurrent-cost}
\Psi_p(\mu;A):=\int_{A\times\overline\Omega\times\orient}f_p(\eta)\d\mu(t,x,\eta).
\end{equation}
For $S\in\intcur{1+k}(I\times\overline\Omega)$ and its elementary varicurrent $\eym{S}$, we have $$\Psi_p(\eym{S};A)=\Psi_p(S;A).$$ For $\mu\in\varicur{1+k}(I\times\overline\Omega)$
 satisfying \eqref{eqn:NC-varicurrent} and $\tproj{}_\#\mu\ll\mathcal L^1$, \eqref{eqn:dissipation-potential-space-time} gives $$\Psi_p(\mu;A)=\int_A\mathcal R_p(\slice{\mu}{t})\d t.$$ The extended cost \eqref{eq:example-varicurrent-cost} also applies to singular varicurrents, for which it need not be representable by a Lebesgue-time slice integral.

\begin{examplecontinued}[Zig-zag, continued]{ex:zigzag-current}

Recall the sequence $S_{\Delta t_i,\alpha_i,N_i}$ from Example~\ref{ex:zigzag-current}. We now consider its elementary varicurrents and retain the decomposition into moving and stationary parts. Set
$$
\mu_{\Delta t,\alpha,N}:=\eym{S_{\Delta t,\alpha,N}}.
$$
We decompose the varicurrent into moving and stationary parts as
$$
\begin{aligned}
\mu_{\Delta t,\alpha,N}^{\mathrm m}&:=\eym{S_{\Delta t,\alpha,N}^{\mathrm m}},
&\mu_{\Delta t,\alpha,N}^{\mathrm s}&:=\eym{S_{\Delta t,\alpha,N}^{\mathrm s}},\\
\mu_{\Delta t,\alpha,N}&=\mu_{\Delta t,\alpha,N}^{\mathrm m}+\mu_{\Delta t,\alpha,N}^{\mathrm s}.
\end{aligned}
$$
We first consider the limits of the moving parts under \eqref{eqn:limits-zigzag}.

For $\Delta t>0$ and $\alpha>0$, we have
$$
\mu_{\Delta t_i,\alpha_i,N_i}^{\mathrm m}\weakstarto\mu_{\Delta t,\alpha}^{\mathrm m},
$$
where
\begin{equation}\label{eqn:zigzag-moving-limit}
\begin{aligned}
\mu_{\Delta t,\alpha}^{\mathrm m}(\mathrm dt,\mathrm dx,\mathrm dy,\mathrm d\eta)
&=\mathcal H^2((0,\Delta t)\times(0,1)\times\{t/\Delta t\})(\mathrm dt, \mathrm dx,\mathrm dy)\\
&\qquad\otimes\frac{\sqrt{1+\beta^2}}{2\Delta t}
\left(\delta_{\eta_{\Delta t,\alpha}^+}+\delta_{\eta_{\Delta t,\alpha}^-}\right)(\mathrm d\eta),\\
\eta_{\Delta t,\alpha}^{\pm}
&:=\frac{e_0\mp\frac{\sin\alpha\cos\alpha}{\Delta t}e_1
+\frac{\sin^2\alpha}{\Delta t}e_2}{\sqrt{1+(\sin\alpha/\Delta t)^2}}
\wedge(\sin\alpha\,e_1\pm\cos\alpha\,e_2).
\end{aligned}
\end{equation}
These are the moving orientations on $0<t<\Delta t$. The stationary tail instead has orientations $e_0\wedge\tau_\pm$ and zero velocity. For the moving limit,
$$
\Psi_p(\mu_{\Delta t,\alpha}^{\mathrm m};[0,1])
=\frac1p\left(\frac{\sin\alpha}{\Delta t}\right)^{p-1}
=\lim_{i\to\infty}\Psi_p(\mu_{\Delta t_i,\alpha_i,N_i}^{\mathrm m};[0,1]).
$$
The assumption $\alpha>0$ concerns the limiting angle and is needed for a finite varicurrent limit. Indeed, if $\Delta t_i\to\Delta t>0$ but $\alpha_i\to0$, then $\beta_i=\Delta t_i/\sin\alpha_i\to\infty$ and
$$
\mu_{\Delta t_i,\alpha_i,N_i}^{\mathrm m}([0,1]\times\overline\Omega\times\orient)
=\mass{S_{\Delta t_i,\alpha_i,N_i}^{\mathrm m}}
=\sqrt{1+\beta_i^2}\longrightarrow\infty.
$$
Since the ambient space is compact, these moving varicurrents cannot converge to a finite measure. Their underlying currents still converge to $S_{\Delta t}^{\mathrm m}$, as proved in Appendix \ref{sec:zigzag-appendix}.

For $\Delta t=0$, assume $\beta_i\to\beta\in(0,\infty)$. Then
$$
\mu_{\Delta t_i,\alpha_i,N_i}^{\mathrm m}\weakstarto\mu_\beta^{\mathrm m},
$$
where
\begin{equation}\label{eqn:zigzag-concentrated-limit}
\begin{aligned}
\mu_\beta^{\mathrm m}(\mathrm dt,\mathrm dx,\mathrm dy,\mathrm d\eta)
&=\delta_0(\mathrm dt)\otimes\mathcal H^2((0,1)\times(0,1))(\mathrm dx,\mathrm dy)\\
&\qquad\otimes\frac{\sqrt{1+\beta^2}}2
\left(\delta_{\eta_\beta^+}+\delta_{\eta_\beta^-}\right)(\mathrm d\eta),\\
\eta_\beta^\pm&:=\frac{\pm \beta\,e_0\wedge e_2+e_2\wedge e_1}{\sqrt{1+\beta^2}}.
\end{aligned}
\end{equation}
The total moving mass is $\sqrt{1+\beta^2}$, and the barycenter current is $\stcurr{\mu_\beta^{\mathrm m}}=S_0^{\mathrm m}$ with orientation $e_2\wedge e_1$. Under this scaling, $\alpha_i\to0$, so the limit depends only on $\beta>0$. Moreover,
$$
|\eta_\beta^\pm \mres\mathrm dt|=\frac{\beta}{\sqrt{1+\beta^2}},\qquad
|v(\eta_\beta^\pm)|=\frac1\beta,
$$
and
$$
\Psi_p(\mu_\beta^{\mathrm m};[0,1])
=\frac1p\beta^{1-p}
=\lim_{i\to\infty}\frac1p\left(\frac{\sin\alpha_i}{\Delta t_i}\right)^{p-1}.
$$
Thus $\mu_\beta^{\mathrm m}$ satisfies \eqref{eqn:NC-varicurrent}, but its mass is concentrated at time zero; its cost is not a Lebesgue-time slice integral. Its purely spatial barycenter $S_0^{\mathrm m}$ fails \eqref{eqn:NC-current} and has infinite extended cost.

The equality of the limiting and approximating costs concerns the moving varicurrents. The full trajectories have the same costs because their tails are stationary. However, for $\beta>0$ the tail masses $(1-\Delta t_i)/\sin\alpha_i$ diverge, so the full elementary varicurrents have no finite measure limit on the stationary time region.

The slice masses $1/\sin\alpha_i$ also diverge, placing this sequence outside the fixed-$\gamma$ compactness result of Proposition \ref{prop:compactness-dissipation-potential}.
\end{examplecontinued}

\begin{examplecontinued}[Thousand Cuts, continued]{ex:thousand-cuts-current}

Recall the moving parts $S_{\Delta t_i,N_i}^{\mathrm m}$ and stationary tails $S_{\Delta t_i,N_i}^{\mathrm s}$ from Example~\ref{ex:thousand-cuts-current}, under \eqref{eqn:limits-cuts}. We again use $\mathrm m$ and $\mathrm s$ for the moving and stationary parts. Set $\mu_{\Delta t_i,N_i}^{\mathrm m}:=\eym{S_{\Delta t_i,N_i}^{\mathrm m}}$. We first identify the moving varicurrent limit and then include the stationary tail. We have
$$
\eym{S_{\Delta t_i,N_i}^{\mathrm m}}\weakstarto\mu_\beta^{\mathrm m}
\quad\text{in }\varicur{2}([0,1]\times\overline\Omega),
$$
where
\begin{align*}
\mu_\beta^{\mathrm m}(\mathrm dt,\mathrm dx,\mathrm dy,\mathrm d\eta)
&=\delta_0(\mathrm dt)\otimes\mathcal H^2((0,1)\times(0,1))(\mathrm dx,\mathrm dy)\\
&\qquad\otimes\frac{\ell_\beta}2
\left(\delta_{\eta_\beta^+}+\delta_{\eta_\beta^-}\right)(\mathrm d\eta),\\
\text{with}\quad\eta_\beta^\pm
&:=\frac{(\pm\beta e_0+e_2)\wedge e_1}{\ell_\beta}.
\end{align*}
Thus the total moving mass is $\mass{\mu_\beta^{\mathrm m}}=\ell_\beta$, and the barycenter current is $\stcurr{\mu_\beta^{\mathrm m}}=S_0^{\mathrm m}$, with orientation $e_2\wedge e_1$. The two orientations occur in equal proportions, due to uniformly alternating orientations.

Combining this moving limit with the stationary tail gives the full-varicurrent limit
$$
\eym{S_{\Delta t_i,N_i}}\weakstarto\mu_\beta^{\mathrm m}+\eym{S_0^{\mathrm s}}.
$$
By \eqref{eqn:cuts-finite-cost}, the costs satisfy
\begin{align*}
\lim_{i\to\infty}\Psi_p(\eym{S_{\Delta t_i,N_i}^{\mathrm m}};[0,1])
&=\lim_{i\to\infty}\frac1p\bigl((2N_i+1)\Delta t_i\bigr)^{1-p}\\
&=\Psi_p(\mu_\beta^{\mathrm m};[0,1])
=\begin{cases}
\frac1p\beta^{1-p},&\beta>0,\\
\infty,&\beta=0.
\end{cases}
\end{align*}
For $\beta>0$, we have $|\eta_\beta^\pm \mres \mathrm dt|=\beta/\ell_\beta$ and $|v(\eta_\beta^\pm)|=1/\beta$. Thus $\mu_\beta^{\mathrm m}$ satisfies \eqref{eqn:NC-varicurrent}, despite its mass concentration at $t=0$. Its purely spatial barycenter $S_0^{\mathrm m}$ fails \eqref{eqn:NC-current} and has infinite $\Psi_p$.

At $\beta=0$, both orientations coincide with $e_2\wedge e_1$ and their combined weight is one; the limit fails \eqref{eqn:NC-varicurrent}, its cost is infinite, and the approximating costs tend to infinity.

Since $N_i\to\infty$, the slice masses satisfy
$$
\mass{S_{\Delta t_i,N_i}(t)}=2N_i+1\longrightarrow\infty
\quad\text{for }0<t<\Delta t_i.
$$
In particular, the uniform-in-time slice-mass bound diverges. This sequence violates the slice-mass hypothesis of Proposition~\ref{prop:compactness-dissipation-potential}, even when $\beta>0$ and the costs remain bounded.
\end{examplecontinued}

\section{The Weighted Energy-Dissipation principle (WED)}\label{sec:WED}

We now use the varicurrent framework to construct evolutions satisfying an energy-dissipation inequality. For each $\varepsilon>0$, the WED principle minimizes an exponentially weighted combination of energy and dissipation potential over complete admissible trajectories with a prescribed initial current. This choice fits the space-time formulation: the competitor is a varicurrent describing the whole evolution, and its slices provide the states and velocities used to evaluate the functional. The slice-mass and dissipation bounds established above give compactness for these minimizers. An auxiliary time discretization will enter later in the estimate for the dual dissipation rate. We will study their energy identities and then pass to the limit as $\varepsilon\to0$ to obtain the desired evolution.

We first specify the energy assumptions and the admissible trajectories, and formulate the solution concept. We then construct WED minimizers and use their value function to derive the estimates needed for the limiting argument.

\subsection{Assumptions on the Energy and Dissipation}

Let the state space for the spatial configurations be the space of integral $k$-currents $\intcur{k}(\overline\Omega)$. We consider an energy functional $\mathcal W \colon \intcur{k}(\overline\Omega) \to (-\infty, \infty]$ satisfying the following assumptions:
\begin{enumerate}[label=(E\arabic*), ref=(E\arabic*)]
    \item\label{item:E-bdd-below} \textbf{Proper and bounded below:} The domain of the energy is non-empty, and it satisfies that $\mathcal W(T) \geq 0$ for all $T \in \intcur{k}(\overline\Omega)$.
    \item\label{item:E-wlsc} \textbf{Weak${}^*$ lower semicontinuity:} For any sequence $(T_n)_{n \in \N} \subset \intcur{k}(\overline\Omega)$ such that $T_n \weakstarto T$, it holds that
    $$\mathcal W(T) \leq \liminf_{n\to\infty} \mathcal W(T_n).$$
    \item\label{item:E-coercivity} \textbf{Coercivity:}
For every $E>0$ there exists $M_E<\infty$ such that
$$
\mass{T}\leq M_E
\quad\text{whenever }
T\in\intcur{k}(\overline\Omega),\quad
\partial T=0,\quad
\mathcal W(T)\leq E.
$$
\item \label{item:E-bdd-by-mass} \textbf{Mass-energy bound:} For any $M>0$ there exists $E_M<\infty$ such that $$\mathcal W(T)\leq E_M\quad \text{whenever }T\in\intcur{k}(\overline{\Omega}), \partial T=0, \mass{T}\leq M.$$
\end{enumerate}

For the avoidance of confusion, the energy $\mathcal W$ considered here is akin to the \emph{elastic} energy in applications to elasto-plasticity, not the \emph{total} energy, which is often denoted by $\mathcal E$. So, we here do not consider the situation involving also an external loading. The reason for this is that the present work focuses on the case of a ``pure'' gradient flow-like system and more concrete applications are left to future works.

Recall from Definition~\ref{def:dissipation-potential} that the dissipation potential is evaluated on the velocity-attached slice. With the fixed exponent $1<p<\infty$, it is given by
$$\mathcal R_p(\slice{\mu}{t}) = \int_{\mathcal Z_k} \frac{1}{p}|v|^p \d\slice{\mu}{t}(x, \tau,v).$$
For $\mu$ satisfying \eqref{eqn:NC-varicurrent} and $\tproj{}_\#\mu \ll \mathcal L^1$, the dissipation potential over a time interval $J$ is
$$\int_J\mathcal R_p(\slice{\mu}{t})\d t.$$ Proposition \ref{prop:compactness-dissipation-potential} establishes preservation of these properties and lower semicontinuity of the integrated dissipation potential under the bounds used below.

\subsection{Admissible Trajectories}\label{sec:admissible-trajectories}
We impose a fixed uniform bound $\gamma>0$ on slice mass, using the framework of Subsection~\ref{ssec:mass-bound-time}. Together with dissipation bounds, this constraint supplies the compactness and time regularity needed for the variational construction. We restrict to trajectories whose barycenter states are boundaryless. Their current slices lie in the space of admissible states
$$\statespace:=\left\{T\in\intcur{k}(\overline\Omega):\partial T=0, \mass{T}\leq \gamma\right\}.$$
This space contains the current states of admissible trajectories, including their initial states. It is compact in the weak${}^*$ topology of currents, and metrizable by the flat metric $\mathbf F$. This follows from the Federer-Fleming compactness theorem recalled in Section \ref{sec:preliminaries} due to the uniform bounds on mass and boundary.

\begin{remark}
    All our results also hold when the setting is restricted to a suitable closed subclass of $X_\gamma$. In particular this includes classes of the form
    $$X_{\gamma,E_0}:= \left\{ T \in X_\gamma: \mathcal W(T) \le E_0\right\}$$
    for given $E_0\in \R$, which are closed by lower semicontinuity and can be useful for some applications. Notably, in this class the mass-energy bound (E4) holds trivially.
\end{remark}

The admissible trajectories form the class $\young$ below. The next lemma provides the continuous current representative needed to prescribe initial data.
$$
\begin{aligned}
\young:= \Big\{&\mu\in\varicur{1+k}^{\loc}([0,\infty)\times\overline\Omega) : \mu\text{ satisfies }\eqref{eqn:NC-varicurrent},
\\
&\tproj{}_\#\mu\ll \mathcal L^1([0,\infty)),
\|\mu\|_{\LinftySlices((0,\infty))}\leq \gamma,
\
\partial \stcurr{\mu}\mres (0,\infty) \times \R^d = 0,
\\
&\int_J \mathcal W(\stcurr{\mu}(t))\d t<\infty,\ \int_J \mathcal R_p(\slice{\mu}{t})\d t<\infty \quad\forall J\subset[0,\infty)\text{ compact}\Big\}.
\end{aligned}
$$

\begin{lemma}[Slicing for all times]\label{lem:all-time-slices-Y}
Let $\mu\in\young$. Then there exists a unique map
\[
T_\mu\in C^{0,1-1/p}_{\loc}
\bigl([0,\infty);(\intcur{k}(\overline\Omega),\mathbf F)\bigr),
\]
which satisfies
\begin{equation}\label{eq:all-time-slice-ae-agreement}
T_{\mu}(t)=\xproj{}_{*}\stcurr{\mu}|_t
\qquad\text{for a.e. }t>0.
\end{equation}
Moreover, it holds that $$T_\mu(t)\in\statespace\quad\text{for every }t\geq 0,$$ and the H\"older continuity estimate
\begin{equation}\label{eq:state-holder-on-J}
\mathbf F\bigl(T_{\mu}(t)-T_{\mu}(s)\bigr)
\leq
\|v\|_{L^p((s,t);L^p(\mu(\cdot)))}
\gamma^{1-1/p}|t-s|^{1-1/p} \quad \text{for every }0
\leq s <t<\infty.
\end{equation}
Furthermore, the cylinder formula
\begin{equation}\label{eq:all-time-cylinder-formula}
\partial\,\xproj{}_{*}
\left(\stcurr{\mu}\mres(s,t)\right)
=T_\mu(t)-T_\mu(s)\quad \text{for every }0
\leq s <t<\infty
\end{equation}
holds, and in particular, the initial value $\stcurr{\mu}(0)$ is represented by the initial boundary
\begin{equation}\label{eq:initial-trace-boundary}
(\partial \stcurr{\mu})\mres\bigl(\{0\}\times\R^d\bigr)
=-\delta_0\times \stcurr\mu(0).
\end{equation}
\end{lemma}
\begin{proof}
By standard slicing results (see Section \ref{sec:preliminaries}),
\[
T_\mu(t):=\xproj{}_{*}\stcurr{\mu}|_t\in\intcur{k}(\overline\Omega)\]
exists for all $t\in G$, where $G \subset (0,\infty)$ is a set of full measure. Since
\[
(\partial \stcurr{\mu})\mres\bigl((0,\infty)\times\R^d\bigr)=0,
\]
by the definition of $\young$, the slicing and cylinder formula \eqref{eqn:current-cylinder-formula} gives, for $s,t\in G$ with $s<t$,
\[
T_\mu(t)-T_\mu(s)
=
\partial\,\xproj{}_{*}
\left(\stcurr{\mu}\mres(s,t)\right).
\]
Therefore, combining  Lemma \ref{lem:flatdist-leq-var} and Lemma \ref{lem:variation-ineq}, we have
\begin{align*}
\mathbf F\bigl(T_\mu(t)-T_\mu(s)\bigr)
\leq \variation{\stcurr{\mu};(s,t)}
\leq \variation{\mu;(s,t)}.
\end{align*}
Recalling also Lemma \ref{lem:p-continuity-in-time}, we obtain the H\"older continuity estimate
\begin{align*}
\mathbf F\bigl(T_\mu(t)-T_\mu(s)\bigr)\leq
\|v\|_{L^p((s,t);L^p(\mu(\cdot)))}
\bigl(\gamma(t-s)\bigr)^{1-1/p}, \quad t,s\in G, s<t.
\end{align*}
Thus, on every compact interval $J\subset[0,\infty)$, the map
$t\mapsto T_\mu(t)$ is uniformly $(1-1/p)$-H\"older continuous on $G\cap J$
with respect to the flat norm. Since $G$ is dense in $[0,\infty)$ and the space $\statespace$ is compact, there exists a unique $\mathbf F$-continuous extension
\[
T_\mu\colon[0,\infty)\to\intcur{k}(\overline\Omega),
\]
satisfying the same H\"older continuity estimate. The cylinder formula \eqref{eqn:current-cylinder-formula} now extends to arbitrary $0\leq s<t<\infty$ by the H\"older continuity of $T_\mu$.
Further, directly from $\mu\in\young$ we have $$\mass{T_\mu(t)}\leq\mass{\slice{\mu}{t}}\leq\gamma,\quad \partial T_\mu(t)=0,\quad \text{for a.e.\ }t\in(0,\infty).$$By the flat continuity, lower semicontinuity of mass, and compactness of integral currents with uniformly bounded mass, these properties extend to all $t$, showing that $T_\mu(t)\in\statespace$ for all $t\in[0,\infty)$. This finishes the proof.
\end{proof}
The mass-energy bound gives local energy integrability for every $\mu\in\young$: $$\int_J\mathcal W(\stcurr{\mu}(t))\d t\leq E_\gamma|J|.$$

\begin{remark}[Slicing the current vs.\ varicurrent]
Lemma~\ref{lem:all-time-slices-Y} gives a uniformly $\mathbf F$-continuous representative of the barycenter current, so its state is defined at \textit{every} time. This is the representative used to prescribe initial data and to evaluate the energy. In contrast, the measure-valued slices of a varicurrent, as defined in Definition~\ref{def:velocity-attached-slices}, are instead defined only for \textit{almost every} time. No additional time regularity of those varicurrent slices is asserted in general. The two notions of slice therefore play different roles, even though their barycenters agree at almost every time.
\end{remark}

\begin{definition}
For $\mu\in\young$, we shall henceforth write
\[
\stcurr{\mu}(t):=T_\mu(t)
\qquad\text{for every }t\geq0.
\]
For $T_0\in\statespace$ we define the set of admissible trajectories
\[
\young(T_0)
:=
\left\{\mu\in\young:\stcurr{\mu}(0)=T_0\right\}.
\]
Equivalently,
\[
\young(T_0)
=
\left\{\mu\in\young:
\partial\stcurr{\mu}
=-\delta_0\times T_0
\right\}.
\]
\end{definition}
\begin{remark}
For any $T_0\in\statespace$, it holds that $\young(T_0)\neq\emptyset$, since it contains the stationary elementary varicurrent $$
\eym{\llbracket(0,\infty)\rrbracket\times T_0} = \d t\otimes \|T_0\|(\mathrm dx)\otimes \delta_{e_0\wedge\vec T_0(x)}(\mathrm d\eta)\in \young(T_0).
$$
\end{remark}

The next proposition combines stability of admissibility with locally uniform convergence of the current states. In addition to obtaining a limiting varicurrent, we need its continuous current representative to retain the prescribed initial value. The locally uniform convergence supplies this link between compactness of trajectories and the initial condition.

\begin{proposition}[Locally uniform convergence in time]\label{prop:uniform-convergence-state-curves}
Let $(\mu_n)_{n\in\N}\subset\young$ satisfy
\[
\mu_n\weakstarto\mu,
\]
for some $\mu\in\varicur{1+k}([0,\infty)\times\overline\Omega)$ and that
\begin{equation}\label{eqn:uniform-local-dissipation-bound}
\sup_{n\in\N}\int_J\mathcal R_p(\slice{\mu_n}{t})\d t<\infty \quad \text{for every compact interval $J\subset [0,\infty)$}.
\end{equation}
Then it holds that $\mu\in \young$
\begin{equation}\label{eqn:uniform-flat-convergence}
\sup_{t\in J}
\mathbf F\bigl(\stcurr{\mu_n}(t)-\stcurr{\mu}(t)\bigr)
\longrightarrow0
\qquad\text{for every compact interval }J\subset [0,\infty).
\end{equation}
In particular,
\[
\stcurr{\mu_n}(t)\weakstarto \stcurr{\mu}(t)
\qquad\text{in }\intcur{k}(\overline\Omega)\text{ for every }t\in[0,\infty).
\]

\end{proposition}

\begin{proof}
By Proposition \ref{prop:compactness-dissipation-potential},
the limit $\mu$ belongs to $\young$. Select by Lemma \ref{lem:compactness-varicurrents} a (non-relabeled) subsequence such that 
\begin{equation}\label{eqn:ae-time-converge-current}
\stcurr{\mu_n}(t)\weakstarto\stcurr{\mu}(t)\quad \text{in }\intcur{k}(\overline\Omega) \quad\text{for a.e. }t\in(0,\infty),
\end{equation}
equivalently, 
\begin{equation*}
\mathbf F\bigl(\stcurr{\mu_n}(t)-\stcurr{\mu}(t)\bigr)\longrightarrow0
\qquad\text{for a.e. }t\in (0,\infty).
\end{equation*}
Fix a compact interval $J\subset[0,\infty)$ and set
\[
C_J
:=
\gamma^{1-1/p}\sup_{n\in\N}\|v\|_{L^p(J;L^p(\mu_n(\cdot)))}
.
\]
By Lemma~\ref{lem:all-time-slices-Y},
\begin{equation}\label{eq:uniform-holder-state-curves}
\mathbf F\bigl(\stcurr{\mu_n}(t)-\stcurr{\mu_n}(s)\bigr)
\leq C_J|t-s|^{1-1/p}
\end{equation}
for every $n\in\N$ and all $s,t\in J$. This shows that the family of maps
$$
t\mapsto\stcurr{\mu_n}(t), n\in\N
$$
is $\mathbf F$-equicontinuous on $J$. By the uniform bound on slice mass and empty boundary, namely
$$
\mass{\stcurr{\mu_n}(t)}\leq \gamma, \quad \partial \stcurr{\mu_n}(t)=0,
$$
the family is also pointwise $\mathbf F$-compact. Therefore, by the Arzel\`a-Ascoli theorem, there is a subsequence converging $\mathbf F$-uniformly in $J$. By a diagonal argument we obtain a subsequence converging $\mathbf F$-uniformly on all compact sets.

Every uniformly convergent subsequence has limit
$\stcurr{\mu}(\cdot)$, by the almost-everywhere convergence
and continuity. Applying the same argument to every
subsequence proves convergence of the full sequence. This concludes the proof.
\end{proof}

In particular, this proposition entails the stability of initial values.
\begin{corollary}[Stability of initial values]\label{cor:stability-initial-values}
Under the assumptions of Proposition~\ref{prop:uniform-convergence-state-curves}, suppose that
\[
\stcurr{\mu_n}(0)\weakstarto T_0 \quad \text{in }\intcur{k}(\overline\Omega).
\]
Then
\[
\stcurr{\mu}(0)=T_0.
\]
In particular, if $\mu_n\in\young(T_0)$, $n\in\N$, then $\mu\in\young(T_0)$.
\end{corollary}

\subsection{The Solution Concept}
In a smooth gradient system, the energy-dissipation formulation involves both the potential evaluated on the velocity and the dual potential $\mathcal R_p^*(-D\mathcal W(\cdot))$ evaluated on the energy derivative. The latter measures the rate at which the energy can decrease relative to the dissipation cost. Here the state space consists of integral currents, and we do not introduce an energy derivative or a separate convex conjugate as separate objects. Instead, we define a \emph{dual dissipation rate} directly through admissible trajectories, comparing their energy decrease with their dissipation cost. This gives the needed quantity without assuming a differential structure on the space of currents.

\begin{definition}[Dual dissipation]\label{def:dual-dissipation}
We define the upper dual dissipation slope at a current $T \in \statespace$ as
$$\dualR{T} := \sup_{ \mu\in\young(T)} \limsup_{\delta\to 0} \left[ \frac{\mathcal W(T)-\mathcal W(\stcurr{\mu}(\delta))}{\delta} - \frac{1}{\delta}\int_0^\delta \mathcal R_p(\slice{\mu}{t}) \d t \right].$$

We define the lower dual dissipation slope $\dualRinf{\cdot}$ as the sequential weak${}^*$ lower semicontinuous relaxation (envelope) of $\dualR{\cdot}$, given by
$$\dualRinf{T} := \inf \left\{ \liminf_{n\to\infty} \dualR{T_n} : T_n\in\statespace,\ T_n \weakstarto T \text{ in } \intcur{k}(\overline\Omega) \right\}.$$
\end{definition}

\begin{remark}
In the smooth Banach-space setting with power-law dissipation, the trajectory-based definition agrees with $\mathcal R_p^*(-D\mathcal W(x))$.

Specifically, let $X$ be a Banach space satisfying the Radon-Nikod\'ym property\footnote{This property guarantees that every $u\in\AC([0,1];X)$ is strongly differentiable a.e.} (e.g. separable and reflexive), let $\mathcal W\in C^1( X; \R)$, and let $\mathcal R_p$ be defined
$$
\mathcal R_p(v):=\frac1p\|v\|_X^p,
\qquad v\in X.
$$
Its convex conjugate is, by the standard definition
$$
\mathcal R_p^*(\xi) = \sup_{v\in X}\bigl(\langle \xi, v\rangle - \mathcal R_p(v)\bigr)=\frac1{p'}\|\xi\|_{X^*}^{p'},\quad \xi\in X^*.
$$
Then
$$\mathcal R_p^*(-D\mathcal W(x))
=
\sup_{\substack{
u\in\AC([0,1];X)\\
u(0)=x}}
\limsup_{\delta\to 0}
\left[
\frac{\mathcal W(x)-\mathcal W(u(\delta))}{\delta}
-\frac1\delta\int_0^\delta
\mathcal R_p(\dot u(t))\,\d t
\right], \quad x\in X.$$
Indeed, for every such $u$ with finite $\mathcal R_p$-integral near 0, we have by the chain rule
\begin{align*}
\frac{\mathcal W(x)-\mathcal W(u(\delta)) }{\delta}-\frac1\delta \int_0^\delta \mathcal R_p(\dot u (t))\d t
&=\frac1\delta \int_0^\delta\langle -D\mathcal W(u(t)),\dot u(t)\rangle-\mathcal R_p(\dot u(t))\d t
\\
&\leq\frac1\delta \int_0^\delta\mathcal R_p^*(-D\mathcal W( u(t)))\d t \to \mathcal R_p^*(-D\mathcal W(x)),
\end{align*}
proving one inequality. The converse inequality follows from taking affine trajectories, as for any direction $v\in X$ it holds 
$$
\lim_{\delta \to 0 }\left[\frac{\mathcal W(x)-\mathcal W(x+\delta v)}{\delta}-\frac 1\delta \int_0^\delta \mathcal R_p(v) \d s\right] = \langle-D\mathcal W(x),v\rangle-\mathcal R_p(v).
$$
\end{remark}

Having a notion of dual dissipation rate established, we define our solution in the framework of the Energy-Dissipation Inequality (EDI), following \cite{mielkeVariationalConvergenceGradient2012} and others.

\begin{definition}[EDI Solution]\label{def:EDI-solution}
Given an initial configuration $T_0 \in\statespace$, we say that a varicurrent $\mu \in \young(T_0)$ is an \emph{EDI solution} to the gradient flow of $\mathcal W$ with dissipation $\mathcal R_p$, if the map  $t \mapsto \dualRinf{\stcurr{\mu}(t)}$ is locally integrable, and it satisfies the Energy-Dissipation Inequality
$$\mathcal W(\stcurr{\mu}(t)) + \int_0^t \mathcal R_p(\slice{\mu}{s}) \d s + \int_0^t \dualRinf{\stcurr{\mu}(s)} \d s \leq \mathcal W(T_0)$$
for all $t \in [0, \infty)$.
\end{definition}

\subsection{WED functional and value function}

To obtain solutions of the EDI, we construct approximate evolutions by the WED principle. Our approach adapts the metric-space construction of \cite{rossiWeightedEnergyDissipationPrinciple2019}, whose terminology we mostly follow.

\begin{definition}[WED functional]
Given $\varepsilon>0$, we define the \emph{WED functional}
\[
\mathcal J_\varepsilon\colon\young\to[0,\infty]
\]
by
$$
\mathcal J_\varepsilon(\mu):=\int_0^\infty\frac{e^{-t/\varepsilon}}{\varepsilon}
\left(\mathcal W(\stcurr{\mu}(t))+\varepsilon\mathcal R_p(\slice{\mu}{t})\right)\d t.
$$
\end{definition}

The first application of the preceding framework is the existence of a WED minimizer.

\begin{lemma} \label{lem:WEDminimizer}
Let $\varepsilon>0$ and $T_0\in\statespace$. Then, there exists a $\mu_\varepsilon\in\young(T_0)$ solving the minimization problem
$$
\mu_\varepsilon\in\operatorname{argmin}\left\{\mathcal J_\varepsilon(\mu):\mu\in\young(T_0)\right\}.
$$
\end{lemma}
\begin{proof}
The functional $\mathcal J_\varepsilon$ is non-negative, and the admissible class $\young(T_0)$ is non-empty because it contains the constant trajectory. Take a sequence $\mu^i\in\young(T_0)$ that realizes the infimum, that is,
$$
\lim_{i\to\infty} \mathcal J_\varepsilon(\mu^i)=\inf_{\mu\in\young(T_0)}\mathcal J_\varepsilon(\mu)=:m.
$$
Because the definition of $\mu^i\in\young(T_0)$ includes the bound $\|\mu^i\|_{\LinftySlices}\leq \gamma$ on the slices $\slice{\mu^i}{t}$, Lemmas \ref{lem:spacetime-mass-bdd-slice-variation} and \ref{lem:p-continuity-in-time} give
$$
\mass{\mu^i \mres J} \leq |J|\|\mu^i\|_{\LinftySlices(J)} + \variation{\mu^i;J}\leq |J|\gamma + \|v\|_{L^p(J;L^p(\mu^i(\cdot)))}\gamma^{1-1/p}|J|^{1-1/p}.
$$
Since $\mu^i$ is an infimizing sequence for $\mathcal J_\varepsilon$, we see that 
$$
\|v\|_{L^p(J;L^p(\mu^i(\cdot)))}=\left(p\int_J\mathcal R_p(\slice{\mu^i}{t})\d t\right)^{1/p}\leq \left( p e^{\frac{\sup J}{\varepsilon}}\mathcal J_\varepsilon(\mu^i)\right)^{1/p},
$$
and hence the mass bound for $\mu^i$ over a compact $J\subset[0,\infty)$ is uniform. Thus, the assumption of Compactness Lemma \ref{lem:compactness-varicurrents} is satisfied and we get
$$
\mu^i\weakstarto\mu_\varepsilon\quad\text{ for some }\mu_\varepsilon\in\varicur{1+k}([0,\infty)\times\overline\Omega)
.$$
Moreover, by Corollary \ref{cor:stability-initial-values} it holds that $\mu_\varepsilon\in\young(T_0)$. It remains to check that 
 $$
 \mathcal J_\varepsilon(\mu_\varepsilon) \leq \liminf_{i\to\infty} \mathcal  J_\varepsilon(\mu^i).
 $$

For the energy term, by the weak${}^*$  lower semicontinuity of $\mathcal W$ in $\intcur{k}(\overline\Omega)$, the convergence $\stcurr{\mu^i}(t)\weakstarto\stcurr{\mu_\varepsilon}(t)$ in $\intcur{k}(\overline\Omega)$ for a.e.\ $t\in(0,\infty)$ and the Fatou lemma, we obtain
$$
\int_0^\infty \frac{e^{-\frac t\varepsilon}}{\varepsilon}\mathcal W(\stcurr{\mu_\varepsilon}(t))\d t\leq
\int_0^\infty \liminf_{i\to\infty}\frac{e^{-\frac t\varepsilon}}{\varepsilon}\mathcal W(\stcurr{\mu^i}(t))\d t\leq
\liminf_{i\to\infty} \int_0^\infty \frac{e^{-\frac t\varepsilon}}{\varepsilon}\mathcal W(\stcurr{\mu^i}(t))\d t.
$$
 For the dissipation term, notice that Lemma \ref{lem:wlsc-dissipation} continues to hold if the integrand is multiplied by a positive weight, and consequently
$$
 \int_0^\infty  e^{-\frac t\varepsilon}\mathcal R_p(\slice{\mu_\varepsilon}{t})\d t\leq \liminf_{i\to\infty}\int_0^\infty  e^{-\frac t\varepsilon}\mathcal R_p(\slice{\mu^i}{t})\d t.
$$
 Together, this finishes the proof.
\end{proof}

For fixed $\varepsilon>0$ and $T_0\in\statespace$, we write $\mu_\varepsilon$ for any chosen minimizer supplied by Lemma \ref{lem:WEDminimizer} and call it an \emph{$\varepsilon$-minimizer}. No uniqueness is asserted.

\subsection{Dynamic programming and the value function}

Having obtained WED minimizers, we introduce the \emph{value function} $V_\varepsilon(T_0)$, which records the minimum WED cost from a prescribed initial state $T_0$. It allows us to compare the remaining cost along a minimizing trajectory with the optimal cost of starting afresh from its current state. This will lead to a dynamic programming principle and a description of the minimizers' evolution in terms of $V_\varepsilon$.

\begin{definition}[Value function]
For a fixed $\varepsilon>0$ and $T_0\in\statespace$, we define the \emph{$\varepsilon$-value function} as
$$
V_\varepsilon(T_0):=\min \left\{\mathcal J_\varepsilon(\mu): \mu\in \young(T_0) \right\}.
$$
Note that then $$V_\varepsilon(T_0)=\mathcal J_\varepsilon(\mu_\varepsilon).$$
\end{definition}

\begin{lemma}[Alternative representation of $V_\varepsilon$]\label{lem:value-func-alt-representation}
For any $\varepsilon > 0$ the value function can be rewritten as
    $$V_\varepsilon(T_0)=\min_{\mu\in\young(T_0)} \int_0^\infty e^{-t}\left(\varepsilon \int_{\mathcal Z_k}\tfrac1p \left|\frac {v}{\varepsilon} \right|^p\d\slice{\mu}{t}+\mathcal W(\stcurr{\mu}(t))\right)\d t,\quad T_0\in\statespace.$$

\end{lemma}

\begin{proof}
Let $\mu\in\young(T_0)$ and put $\tilde\mu:= \timepush{\left( t\mapsto t/\varepsilon\right)}{\mu}$. By Proposition \ref{prop:time-pushforward} and Lemma \ref{lem:dissipation-pushforward}, we have $\tilde\mu\in\young(T_0)$. Then
\begin{align*}
\mathcal J_\varepsilon(\mu)
&=\int_0^\infty\frac{e^{-\frac t\varepsilon}}{\varepsilon}\left(\varepsilon\int\tfrac1p|v|^p\d\slice{\mu}{t}+ \mathcal W(\stcurr{\mu}(t))\right)\d t
\\
&=\int_0^\infty  e^{-s}\left(\varepsilon \int\tfrac1p \left|\frac {v}{\varepsilon} \right|^p\d\slice{\tilde \mu}{s} +\mathcal W(\stcurr{\tilde \mu}(s))\right)\d s
\end{align*} As $\mu\mapsto\tilde \mu$ is a bijection, we can get $V_\varepsilon(T_0)$ by minimizing over $\tilde \mu$ instead of over $\mu$, finishing the proof.
\end{proof}

The exponential weight is preserved, up to a constant factor, under time translation. This gives a dynamic programming principle, from which we obtain the dissipation identity for $V_\varepsilon$ along minimizers.

\begin{proposition}[Dynamic programming principle]\label{prop:dynamic-programming-principle}
It holds for all $t\in(0,\infty)$ that $$V_\varepsilon (T_0)=\min_{\mu\in\young(T_0)}\left\{\int_0^t e^{-\frac s\varepsilon}\left(\mathcal R_p(\slice{\mu}{s}) + \frac1\varepsilon \mathcal W(\stcurr{\mu}(s))\right)\d s+e^{-\frac t\varepsilon} V_\varepsilon(\stcurr{\mu}(t))\right\}.$$
Moreover, any minimizer of $\mathcal J_\varepsilon$ is a minimizer of the functional on the right.
\end{proposition}

\begin{proof}
Let us by Lemma \ref{lem:WEDminimizer} pick a minimizer $$\mu_\varepsilon\in\operatorname*{argmin}_{\mu\in\young(T_0)}\mathcal J_\varepsilon(\mu).$$
Denote the expression on the right hand side as $W_\varepsilon(T_0)$.
First, we will show $W_\varepsilon(T_0)\leq V_\varepsilon(T_0)$.
Let $\mu\in\young(T_0)$ be arbitrary and fix $t>0$. Then, as $\timepush{(s\mapsto s-t)} {\mu}\mres [0,\infty)$ is an admissible competitor in $V_\varepsilon(\stcurr{\mu}(t))$, the following holds for every $t$ by Lemma \ref{lem:all-time-slices-Y}:
\begin{align*}
V_\varepsilon(\stcurr{\mu}(t))&\leq \int_0^\infty e^{-\frac r\varepsilon}\left(\mathcal R_p(\slice{\mu}{t+r}) + \frac1\varepsilon \mathcal W(\stcurr{\mu}(t+r))\right)\d r\\
&=
e^{\frac t\varepsilon}\int_t^\infty e^{-\frac s\varepsilon}\left(\mathcal R_p(\slice{\mu}{s}) + \frac1\varepsilon \mathcal W(\stcurr{\mu}(s))\right)\d s.
\end{align*}
Thus, we see that 
$$
W_\varepsilon(T_0)\leq \int_0^t e^{-\frac s\varepsilon}\left(\mathcal R_p(\slice{\mu}{s}) + \frac1\varepsilon \mathcal W(\stcurr{\mu}(s))\right)\d s+e^{-\frac t\varepsilon} V_\varepsilon(\stcurr{\mu}(t)) \leq \mathcal J_\varepsilon(\mu)
$$
and taking the infimum over all such $\mu$ gives $W_\varepsilon(T_0)\leq V_\varepsilon(T_0)$.

To get the other inequality $V_\varepsilon(T_0)\leq W_\varepsilon(T_0)$, let $\mu\in\young(T_0)$ be arbitrary and let 
$$\tilde \mu \in\operatorname*{argmin}_{\young(\stcurr{\mu}(t))}\mathcal J_\varepsilon.$$
Using the competitor
$$
\nu:=\mu\mres [0,t)
+\timepush{(s\mapsto s+t)}{\tilde\mu}.
$$
The boundary $\delta_t\times\stcurr{\mu}(t)$ of the first current cancels the initial boundary of the second, because $\stcurr{\tilde\mu}(0)=\stcurr{\mu}(t)$. Using $\nu$ in the minimization for $V_\varepsilon(T_0)$, we obtain
$$
V_\varepsilon(T_0)\leq \int_0^t e^{-\frac s\varepsilon}\left(\mathcal R_p(\slice{\mu}{s}) + \frac1\varepsilon \mathcal W(\stcurr{\mu}(s))\right)\d s + e^{-\frac t\varepsilon}\mathcal J_\varepsilon(\tilde \mu).
$$
Upon realizing that $\mathcal J_\varepsilon(\tilde \mu)=V_\varepsilon(\stcurr{\mu}(t))$ and taking the infimum over all $\mu$, the inequality  $V_\varepsilon(T_0)\leq W_\varepsilon(T_0)$ follows.
\end{proof}

We introduce the dual dissipation rate $G_\varepsilon$, corresponding to the $\varepsilon$-value function instead of the energy functional. This should be seen as the $\varepsilon$-level analogue of $\dualR{}$ from Definition~\ref{def:dual-dissipation}, as explained below.

The quantity $G_\varepsilon$ plays the role of a dual
dissipation rate: In a smooth gradient system with energy
$\mathcal W$ and dissipation potential $\mathcal R$, this
rate is $\mathcal R^*(-D\mathcal W(\cdot))$, where
$\mathcal R^*$ is the convex conjugate of $\mathcal R$.
Here the corresponding energy is $V_\varepsilon$, so
$G_\varepsilon$ is the analogue of
$\mathcal R^*(-D V_\varepsilon)$.
Its associated slope factor $(p'G_\varepsilon)^{1/p'}$
gives an upper-gradient estimate for $V_\varepsilon$
along admissible trajectories, as shown in
Proposition \ref{prop:chain-rule-epsilon}.

\begin{definition}[$\varepsilon$-dual dissipation rate]
We denote the \emph{$\varepsilon$-slope} by $$G_\varepsilon(T_0):=\frac{\mathcal W(T_0)-V_\varepsilon(T_0)}{\varepsilon},\quad T_0\in\statespace.$$
\end{definition}

\begin{proposition}[Dissipation identity for the value function]\label{prop:value-funct-gradient-flow}
The function
$$
t\mapsto V_\varepsilon(\stcurr{\mu_\varepsilon}(t)) \quad \text{is locally absolutely continuous on }[0,\infty).
$$
Moreover, it holds
\begin{equation}\label{eqn:value-funct-gradient-flow}
-\frac{\mathrm d}{\mathrm dt}V_\varepsilon(\stcurr{\mu_\varepsilon}(t))=\mathcal R_p(\slice{\mu_\varepsilon}{t})+G_\varepsilon(\stcurr{\mu_\varepsilon}(t))\quad \text{for a.e.\ }t\in(0,\infty).
\end{equation}
\end{proposition}
\begin{proof}
Let us denote
$$
f_\varepsilon(t)=e^{-\frac t\varepsilon} V_\varepsilon(\stcurr{\mu_\varepsilon}(t)), \quad t\in [0,\infty).
$$
We have by the Dynamic Programming principle (Proposition \ref{prop:dynamic-programming-principle}) that it holds for all $0\leq t <\infty$,
$$
f_\varepsilon(0)- f_\varepsilon(t)=\int_0^t e^{-\frac s\varepsilon}\left(\mathcal R_p(\slice{\mu_\varepsilon}{s}) + \frac1\varepsilon \mathcal W(\stcurr{\mu_\varepsilon}(s))\right)\d s.
$$
Now, the integrand on the right hand side satisfies
$$
s\mapsto e^{-\frac s\varepsilon}\left(\mathcal R_p(\slice{\mu_\varepsilon}{s}) + \frac1\varepsilon \mathcal W(\stcurr{\mu_\varepsilon}(s))\right)
\quad\text{is non-negative and integrable on } (0,\infty),
$$
and consequently $f_\varepsilon$ is absolutely continuous and decreasing on $(0,\infty)$
with 
$$
\dot f_\varepsilon(t)= -e^{-\frac t\varepsilon}\left(\mathcal R_p(\slice{\mu_\varepsilon}{t}) + \frac1\varepsilon \mathcal W(\stcurr{\mu_\varepsilon}(t))\right),\quad \text{ for a.e.\ }t\in(0,\infty).
$$
Now, by the product rule for $\AC$ functions we have that
$$
t\mapsto V_\varepsilon(\stcurr{\mu_\varepsilon}(t))=e^{\frac t\varepsilon} f_\varepsilon(t) \quad \text{ is locally absolutely continuous on }[0,\infty)
$$
with
$$
\begin{aligned}
\frac{\mathrm d}{\mathrm dt}V_\varepsilon(\stcurr{\mu_\varepsilon}(t))
&=\frac{\mathrm d}{\mathrm dt}\left( e^{\frac t\varepsilon} f_\varepsilon(t)\right)
= \frac{e^{\frac t \varepsilon}}{\varepsilon} f_\varepsilon(t) + e^{\frac t \varepsilon}\dot f_\varepsilon(t)
\\
&=\frac1\varepsilon V_\varepsilon(\stcurr{\mu_\varepsilon}(t)) -\left(\mathcal R_p(\slice{\mu_\varepsilon}{t}) + \frac1\varepsilon \mathcal W(\stcurr{\mu_\varepsilon}(t))\right),\quad \text{ for a.e.\ }t\in (0,\infty),
\end{aligned}
$$
which upon rearranging the terms is the desired equality.
\end{proof}

 We next establish an upper-gradient estimate for $V_\varepsilon$ along arbitrary admissible trajectories.

\begin{proposition}[Upper-gradient estimate for $V_\varepsilon$]\label{prop:chain-rule-epsilon}
Let $\varepsilon>0$. For every $\mu\in\young$, the map
$
t\mapsto V_\varepsilon(\stcurr{\mu}(t))
$
is locally absolutely continuous on $[0,\infty)$ and satisfies
\begin{equation}\label{eq:upper-gradient-epsilon}
\left|\frac{\mathrm d}{\mathrm d t}V_\varepsilon(\stcurr{\mu}(t))\right|
\leq \|v\|_{L^p(\mu(t))}
\left(p'G_\varepsilon(\stcurr{\mu}(t))\right)^{1/p'}
\quad\text{for a.e. }t\in(0,\infty).
\end{equation}
In particular, it holds
$$
\left|\frac{\mathrm d}{\mathrm dt}V_\varepsilon(\stcurr{\mu}(t))\right|\leq \mathcal R_p(\slice{\mu}{t}) +G_\varepsilon(\stcurr{\mu}(t)) \quad \text{for a.e. }t\in(0,\infty).
$$
\end{proposition}

\begin{proof}
Recall that the stationary competitor gives
$$
0\leq V_\varepsilon(T)\leq\mathcal W(T)
\quad\text{for every }T\in\statespace,
$$
and consequently $G_\varepsilon\geq0$.

Take any $\mu\in\young$ and $0\leq r<s<\infty$. We first establish the additive estimate. Using the Dynamic Programming Principle (Proposition \ref{prop:dynamic-programming-principle}) with initial state $\stcurr{\mu}(r)$, at time $s-r$, and the competitor
$$
\timepush{(t\mapsto t-r)}
{(\mu\mres (r,\infty))}
\in\young(\stcurr{\mu}(r)),
$$
we obtain
$$
V_\varepsilon(\stcurr{\mu}(r))
\leq
\int_0^{s-r}e^{-t/\varepsilon}
\left(
\mathcal R_p(\slice{\mu}{t+r})
+\frac1\varepsilon\mathcal W(\stcurr{\mu}(t+r))
\right)\d t
+e^{-(s-r)/\varepsilon}V_\varepsilon(\stcurr{\mu}(s)).
$$
Subtracting $V_\varepsilon(\stcurr{\mu}(s))$ and changing time variables gives
\begin{equation}\label{eqn:Vepsilon-ineq1}
\begin{aligned}
&V_\varepsilon(\stcurr{\mu}(r))-V_\varepsilon(\stcurr{\mu}(s))
\\
&\leq
\left(e^{-(s-r)/\varepsilon}-1\right)V_\varepsilon(\stcurr{\mu}(s))+
\int_r^s e^{-(t-r)/\varepsilon}\left(\mathcal R_p(\slice{\mu}{t})+\frac1\varepsilon\mathcal W(\stcurr{\mu}(t))\right)\d t.
\end{aligned}
\end{equation}

Let $\nu$ be the reversal (Lemma \ref{lem:time-reversal}) of $\mu$ on $(r,s)$ extended stationarily, that is, 
$$
\nu=
\left(\timepush{(t\mapsto t-r)}
{\mu\mres(r,s)}\right)^{\mathrm{rev}}
+\eym{\llbracket(s-r,\infty)\rrbracket\times\stcurr{\mu}(r)}.
$$
Consequently, we have $\nu\in\young(\stcurr{\mu}(s))$, with $\stcurr{\nu}(t)=\stcurr{\mu}(s-t)$ for $0\leq t\leq s-r$ and $\mathcal R_p(\slice{\nu}{t})=\mathcal R_p(\slice{\mu}{s-t})$ for a.e. $t\in(0,s-r)$.

Applying the Dynamic Programming Principle to $\nu$ yields
$$
V_\varepsilon(\stcurr{\mu}(s))
\leq
\int_0^{s-r}e^{-t/\varepsilon}
\left(
\mathcal R_p(\slice{\mu}{s-t})
+\frac1\varepsilon\mathcal W(\stcurr{\mu}(s-t))
\right)\d t
+e^{-(s-r)/\varepsilon}V_\varepsilon(\stcurr{\mu}(r)).
$$
Subtracting $V_\varepsilon(\stcurr{\mu}(r))$ and changing variables gives
\begin{equation}\label{eqn:Vepsilon-ineq2}
\begin{aligned}
&V_\varepsilon(\stcurr{\mu}(s))
-V_\varepsilon(\stcurr{\mu}(r))
\\
&\leq
\left(e^{-(s-r)/\varepsilon}-1\right)
V_\varepsilon(\stcurr{\mu}(r))+
\int_r^s e^{-(s-t)/\varepsilon}
\left(
\mathcal R_p(\slice{\mu}{t})
+\frac1\varepsilon\mathcal W(\stcurr{\mu}(t))
\right)\d t.
\end{aligned}
\end{equation}

Since $V_\varepsilon\geq0$ and the exponential weights are at most one,
\eqref{eqn:Vepsilon-ineq1} and \eqref{eqn:Vepsilon-ineq2} imply
$$
\left|
V_\varepsilon(\stcurr{\mu}(s))
-V_\varepsilon(\stcurr{\mu}(r))
\right|
\leq
\int_r^s
\left(
\mathcal R_p(\slice{\mu}{t})
+\frac1\varepsilon\mathcal W(\stcurr{\mu}(t))
\right)\d t.
$$
The integrand belongs to $L^1_{\loc}([0,\infty))$ by the definition
of $\young$. Hence
$t\mapsto V_\varepsilon(\stcurr{\mu}(t))$
is locally absolutely continuous on $[0,\infty)$.

Now divide \eqref{eqn:Vepsilon-ineq1} and
\eqref{eqn:Vepsilon-ineq2} by $s-r$ and let $s\to r$.
At differentiability points of
$t\mapsto V_\varepsilon(\stcurr{\mu}(t))$,
which are also right Lebesgue points of the integrand, the two
inequalities give
\begin{equation}\label{eq:upper-gradient-additive}
\begin{aligned}
\left|
\frac{\mathrm d}{\mathrm dt}
V_\varepsilon(\stcurr{\mu}(t))
\right|
&\leq
\mathcal R_p(\slice{\mu}{t})
+\frac1\varepsilon\mathcal W(\stcurr{\mu}(t))
-\frac1\varepsilon V_\varepsilon(\stcurr{\mu}(t))
\\
&=
\mathcal R_p(\slice{\mu}{t})
+G_\varepsilon(\stcurr{\mu}(t))
\quad\text{for a.e. }t>0.
\end{aligned}
\end{equation}

To obtain the multiplicative estimate, apply this inequality to the
increasing time rescaling
$$
\mu_\lambda:=\timepush{(t\mapsto t/\lambda)}{\mu},
\qquad \lambda>0.
$$
By Proposition \ref{prop:time-pushforward} and
Lemma \ref{lem:dissipation-pushforward}, we have
$\mu_\lambda\in\young$ and
$$
\stcurr{\mu_\lambda}(t/\lambda)=\stcurr{\mu}(t) \quad \text{for every }t\in[0,\infty),
\qquad
\mathcal R_p(\slice{\mu_\lambda}{t/\lambda})
=\lambda^p\mathcal R_p(\slice{\mu}{t}) \quad \text{for a.e. }t\in(0,\infty).
$$
Thus \eqref{eq:upper-gradient-additive} applied to $\mu_\lambda$ gives for a fixed $\lambda>0$:
$$
\lambda
\left|
\frac{\mathrm d}{\mathrm dt}
V_\varepsilon(\stcurr{\mu}(t))
\right|
\leq
\lambda^p\mathcal R_p(\slice{\mu}{t})
+G_\varepsilon(\stcurr{\mu}(t))\quad \text{for a.e. }t.
$$
Taking $\lambda$ from a countable dense set $D\subset(0,\infty)$, we obtain a common
full-measure set of times on which this inequality holds for all $\lambda\in D$. Continuity in $\lambda$ then extends this inequality from $\lambda \in D$ to all $\lambda>0$.

Therefore, writing
$\mathcal R_p(\slice{\mu}{t})=\frac1p\|v\|_{L^p(\mu(t))}^p$,
dividing by $\lambda$, we obtain for a.e. $t\in(0,\infty)$
$$
\left|
\frac{\mathrm d}{\mathrm dt}
V_\varepsilon(\stcurr{\mu}(t))
\right|
\leq
\inf_{\lambda>0}
\left(
\frac{\lambda^{p-1}}p\|v\|_{L^p(\mu(t))}^p
+\frac{G_\varepsilon(\stcurr{\mu}(t))}{\lambda}
\right)=
\|v\|_{L^p(\mu(t))}
\left(p'G_\varepsilon(\stcurr{\mu}(t))\right)^{1/p'},
$$
proving \eqref{eq:upper-gradient-epsilon}.

Finally, we remark that Young's inequality gives
$$
\|v\|_{L^p(\mu(t))}
\left(p'G_\varepsilon(\stcurr{\mu}(t))\right)^{1/p'}
\leq
\frac1p\|v\|_{L^p(\mu(t))}^p
+G_\varepsilon(\stcurr{\mu}(t))=
\mathcal R_p(\slice{\mu}{t})
+G_\varepsilon(\stcurr{\mu}(t)),
$$
which additionally shows that the additive estimate follows from the multiplicative one.
\end{proof}

Combining the upper-gradient estimate with the dissipation
identity along minimizers forces equality in Young's
inequality and identifies $G_\varepsilon$ along minimizers.

\begin{proposition}\label{prop:eps-slope-is-R*DR} 
It holds for a.e.\ $t\in(0,\infty)$ that $$G_\varepsilon(\stcurr{\mu_\varepsilon}(t))=(p-1)\mathcal R_p(\slice{\mu_\varepsilon}{t}).$$
\end{proposition}
\begin{proof}
The stationary competitor gives $0\leq V_\varepsilon\leq\mathcal W$, hence $G_\varepsilon\geq0$. Since also $\mathcal R_p\geq0$, the identity \eqref{eqn:value-funct-gradient-flow}, the upper-gradient estimate of Proposition \ref{prop:chain-rule-epsilon}, and Young's inequality give, for a.e.\ $t>0$,
$$
\begin{aligned}
&\mathcal R_p(\slice{\mu_\varepsilon}{t})
+G_\varepsilon(\stcurr{\mu_\varepsilon}(t))
\\
&=-\frac{\mathrm d}{\mathrm dt}V_\varepsilon(\stcurr{\mu_\varepsilon}(t))
\\
&\leq\|v\|_{L^p(\mu_\varepsilon(t))}
\left(p'G_\varepsilon(\stcurr{\mu_\varepsilon}(t))\right)^{1/p'}
\\
&\leq\frac1p\|v\|_{L^p(\mu_\varepsilon(t))}^p
+G_\varepsilon(\stcurr{\mu_\varepsilon}(t))
\\
&=\mathcal R_p(\slice{\mu_\varepsilon}{t})
+G_\varepsilon(\stcurr{\mu_\varepsilon}(t)).
\end{aligned}
$$
Both inequalities are therefore equalities. Consequently, having equality in Young's inequality gives
$$
p\mathcal R_p(\slice{\mu_\varepsilon}{t})
=p'G_\varepsilon(\stcurr{\mu_\varepsilon}(t)),
$$
which is the claimed identity. The equality criterion also includes the case in which either factor vanishes.
\end{proof}

\subsection{Energy Balance for \texorpdfstring{$\mu_\varepsilon$}{μ\_ε}}
The characterization of $G_\varepsilon$ along minimizers now gives an exact energy balance for $V_\varepsilon$. This balance relates the decrease of the value function to the dissipation along the minimizing trajectory. It will supply the uniform estimates used to extract a limit as $\varepsilon\to0$.

\begin{proposition}[Energy balance for $V_\varepsilon$]\label{prop:eps-energy-balance} 
It holds for all $t\in(0,\infty)$
$$V_\varepsilon(\stcurr{\mu_\varepsilon}(t))+\int_0^t p\mathcal R_p(\slice{\mu_\varepsilon}{s})\d s = V_\varepsilon(T_0).$$
\end{proposition}

\begin{proof}
The proof follows by using the representation of Proposition \ref{prop:eps-slope-is-R*DR} in Proposition \ref{prop:value-funct-gradient-flow} and integrating over $(0,t)$.
\end{proof}


\begin{proposition}[Energy estimates]\label{prop:energy-estimates} The minimizer $\mu_\varepsilon$ satisfies the following energy estimates:
$$\begin{aligned}
V_\varepsilon(\stcurr{\mu_\varepsilon}(t))&\leq \mathcal W(T_0)\quad\text{for all } t\in(0,\infty),
\\
p\int_0^t\mathcal R_p(\slice{\mu_\varepsilon}{s})\d s&\leq \mathcal W(T_0),
\\
\int_0^t \mathcal W(\stcurr{\mu_\varepsilon}(s))\d s&\leq t\mathcal W(T_0)+\varepsilon V_\varepsilon(T_0),
\end{aligned}$$

\end{proposition}
\begin{proof}
Integrating \eqref{eqn:value-funct-gradient-flow} in time, we get $$V_\varepsilon(\stcurr{\mu_\varepsilon}(t))+\int_0^tp\mathcal R_p(\slice{\mu_\varepsilon}{s})\d s = V_\varepsilon(T_0).$$
Further, we have $V_\varepsilon(T_0)\leq \mathcal W(T_0)$, since the constant trajectory $\eym{\llbracket [0,\infty)\rrbracket \times T_0}$ is an admissible competitor in the minimization defining $V_\varepsilon(T_0)$. This yields the first two desired estimates
$$V_\varepsilon(\stcurr{\mu_\varepsilon}(t))\leq \mathcal W(T_0),
\qquad p\int_0^t\mathcal R_p(\slice{\mu_\varepsilon}{s})\d s\leq \mathcal W(T_0).$$

Further, we have by Proposition~\ref{prop:eps-slope-is-R*DR}  that $$\mathcal W(\stcurr{\mu_\varepsilon}(s))= V_\varepsilon(\stcurr{\mu_\varepsilon}(s))+\varepsilon(p-1)\mathcal R_p(\slice{\mu_\varepsilon}{s}).$$
Combining with the above, we have 
$$\mathcal W(\stcurr{\mu_\varepsilon}(s))\leq \mathcal W(T_0)+\varepsilon (p-1)\mathcal R_p(\slice{\mu_\varepsilon}{s}).$$
Integrating this over $(0,t)$ yields the third estimate \begin{align*}\int_0^t \mathcal W(\stcurr{\mu_\varepsilon}(s))\d s&\leq t\mathcal W(T_0)+\varepsilon V_\varepsilon(T_0).\qedhere \end{align*}
\end{proof}

We also record the identity obtained by inner variations in time. This is an independent structural result and is not needed for the existence argument below.

\begin{proposition}[Inner variation equation]\label{prop:inner-variation-eqn}
The \emph{inner variation equation}
$$
\frac{\mathrm d}{\mathrm dt}\left(\mathcal W(\stcurr{\mu_\varepsilon}(t))-\varepsilon(p-1)\mathcal R_p(\slice{\mu_\varepsilon}{t}) \right)=-p\mathcal R_p(\slice{\mu_\varepsilon}{t})  \quad \text{in } \mathcal D'((0,\infty))
$$
holds, that is,
\[
\int_0^\infty \dot\psi(t)\left(\mathcal W(\stcurr{\mu_\varepsilon}(t))-\varepsilon(p-1)\mathcal R_p(\slice{\mu_\varepsilon}{t}) \right) \d t
=
\int_0^\infty\psi(t) p\mathcal R_p(\slice{\mu_\varepsilon}{t}) \d t
\]
for all $\psi\in C^1_c((0,\infty))$.
\end{proposition}

\begin{proof}
By rewriting the above weak formulation with $\tilde\psi(t):=\frac{e^{-\frac t \varepsilon}}{\varepsilon}\psi(t)$ in place of $\psi$, 
we see that it is equivalent to show that
$$
\begin{aligned}
0=&\int_0^\infty e^{-\frac {t}\varepsilon}\frac{-\psi(t)}{\varepsilon}p\mathcal R_p(\slice{\mu_\varepsilon}{t})\d t
\\
&+\int_0^\infty e^{-\frac {t}\varepsilon}\left(\frac{\psi(t)}{\varepsilon}-\dot\psi(t)\right)\left((p-1)\mathcal R_p(\slice{\mu_\varepsilon}{t}) -\frac1\varepsilon\mathcal W(\stcurr{\mu_\varepsilon}(t))\right)\d t
\quad\text{for all }\psi\in C^1_c((0,\infty)).
\end{aligned}
$$
To prove this, we will take the inner variation in time of the functional $\mathcal J_\varepsilon$.

Fix $\psi\in C^1_c((0,\infty))$. Then for all $0<\sigma<1/\|\dot\psi\|_{\infty}$ we have that the function $$\varphi_\sigma(t)=t+\sigma\psi(t),\quad t\in(0,\infty)$$ is monotone increasing and thus invertible, with $$(\varphi_\sigma^{-1})'(s)=\frac{1}{1+\sigma\dot\psi(\varphi_\sigma^{-1}(s))},\quad s\in(0,\infty).$$
Consider perturbations of $\mu_\varepsilon$ given by time pushforward (recall Definition \ref{def:time-pushforward}) $$\mu_\varepsilon^\sigma:=\timepush{(\varphi_\sigma)}{\mu_\varepsilon}.$$
We will now compute $$\left.\frac{\mathrm d}{\mathrm d\sigma}\right|_{\sigma=0}\mathcal J_\varepsilon(\mu_\varepsilon^\sigma).$$
If this derivative exists, it has to be equal to $0$ due to $\mu_\varepsilon=\mu_\varepsilon^0$ being the minimizer of $\mathcal J_\varepsilon$. 

Let us first rewrite 
$$\begin{aligned}
&\mathcal J_\varepsilon(\mu_\varepsilon^\sigma)=\int \ e^{-\frac s\varepsilon}\tfrac1p|v(\eta)|^p|\eta \mres \mathrm dt|\d\mu_\varepsilon^\sigma(s,x,\eta)+\int_0^\infty\frac{e^{-\frac s\varepsilon}}{\varepsilon} \mathcal W(\stcurr{\mu_\varepsilon}(\varphi_\sigma^{-1}(s)))\d s
\\
&=\int  e^{-\frac {t+\sigma\psi(t)}\varepsilon}\tfrac1p\left|\tfrac{v(\eta)}{1+\sigma\dot\psi(t)}\right|^p(1+\sigma\dot\psi(t))|\eta \mres \mathrm dt|\d\mu_\varepsilon+\int_0^\infty \frac{e^{-\frac {t+\sigma\psi(t)}\varepsilon}}{\varepsilon}\mathcal W(\stcurr{\mu_\varepsilon}(t))(1+\sigma\dot\psi(t))\d t,
\end{aligned}$$
where for the dissipation we used the reparametrization Lemma \ref{lem:dissipation-pushforward} as well as $\stcurr{\mu_\varepsilon^\sigma}(s)=\stcurr{\mu_\varepsilon}(\varphi_\sigma^{-1}(s))$ with the change of variables $s=\varphi_\sigma(t)$, $\d s=(1+\sigma\dot\psi(t))\d t$ in the $\mathcal W$ term.

A uniform upper bound of the integrand is provided by a multiple of $\mathcal W(\stcurr{\mu_\varepsilon}(t))+\mathcal R_p(\slice{\mu_\varepsilon}{t})$, which is integrable on the compact support of $\psi$.  Using the Lebesgue dominated convergence theorem to justify differentiation under the integral sign we can thus compute
$$\begin{aligned}
\frac{\mathrm d}{\mathrm d\sigma}\mathcal J_\varepsilon(\mu_\varepsilon^\sigma)
=&\int\left[ e^{-\frac {t+\sigma\psi(t)}\varepsilon}\frac{-\psi(t)}{\varepsilon}\tfrac1p\left|\frac{v(\eta)}{1+\sigma\dot\psi(t)}\right|^p(1+\sigma\dot\psi(t))\right.
\\
&+e^{-\frac {t+\sigma\psi(t)}\varepsilon}\frac{-|v(\eta)|^p}{|1+\sigma\dot\psi(t)|^{p+1}}\dot\psi(t)(1+\sigma\dot\psi(t))
\\
&\left.+e^{-\frac {t+\sigma\psi(t)}\varepsilon}\tfrac1p\left|\frac{v(\eta)}{1+\sigma\dot\psi(t)}\right|^p\dot\psi(t)\right]|\eta \mres \mathrm dt|\d\mu_\varepsilon
\\
&+\int_0^\infty \frac{e^{-\frac {t+\sigma\psi(t)}\varepsilon}}{\varepsilon}\frac{-\psi(t)}{\varepsilon}\mathcal W(\stcurr{\mu_\varepsilon}(t))+\frac{e^{-\frac {t+\sigma\psi(t)}\varepsilon}}{\varepsilon}\mathcal W(\stcurr{\mu_\varepsilon}(t))\dot\psi(t)\d t.
\end{aligned}$$
Thus we obtain, since $\mu_\varepsilon=\mu_\varepsilon^0$ minimizes $\mathcal J_\varepsilon$,
$$\begin{aligned}
&0=\left.\frac{\mathrm d}{\mathrm d\sigma}\right|_{\sigma=0}\mathcal J_\varepsilon(\mu_\varepsilon^\sigma )\\
&\quad=\int \Big[  e^{-\frac {t}\varepsilon}\frac{-\psi(t)}{\varepsilon}\tfrac1p|v(\eta)|^p+e^{-\frac {t}\varepsilon}(-|v(\eta)|^p\dot\psi(t))
+e^{-\frac {t}\varepsilon}\tfrac1p|v(\eta)|^p\dot\psi(t)
\Big] |\eta \mres \mathrm dt|\d\mu_\varepsilon
\\
&\quad+\int_0^\infty \frac{e^{-\frac {t}\varepsilon}}{\varepsilon}\frac{-\psi(t)}{\varepsilon}\mathcal W(\stcurr{\mu_\varepsilon}(t))+\frac{e^{-\frac {t}\varepsilon}}{\varepsilon}\mathcal W(\stcurr{\mu_\varepsilon}(t))\dot\psi(t)\d t.
\end{aligned}$$

Rearranging the terms, adding and subtracting, we obtain $$\begin{aligned}
0=\left.\frac{\mathrm d}{\mathrm d\sigma}\right|_{\sigma=0}\mathcal J_\varepsilon(\mu_\varepsilon^\sigma)
=&\int  e^{-\frac {t}\varepsilon}\frac{-\psi(t)}{\varepsilon}|v(\eta)|^p
|\eta \mres \mathrm dt|\d\mu_\varepsilon
\\
&+\int e^{-\frac {t}\varepsilon}\left(\frac{\psi(t)}{\varepsilon}-\dot\psi(t)\right)\tfrac{1}{p'}|v(\eta)|^p|\eta \mres \mathrm dt|
\d\mu_\varepsilon
\\
&+\int_0^\infty e^{-\frac {t}\varepsilon}\frac1\varepsilon
\left(\frac{-\psi(t)}{\varepsilon}+\dot\psi(t)\right)\mathcal W(\stcurr{\mu_\varepsilon}(t))\d t
\\
=&\int_0^\infty  e^{-\frac {t}\varepsilon}\frac{-\psi(t)}{\varepsilon} p\mathcal R_p(\slice{\mu_\varepsilon}{t})\d t
\\
&+\int_0^\infty  e^{-\frac {t}\varepsilon}\left(\frac{\psi(t)}{\varepsilon}-\dot\psi(t)\right)\left((p-1)\mathcal R_p(\slice{\mu_\varepsilon}{t}) -\frac1\varepsilon \mathcal W(\stcurr{\mu_\varepsilon}(t))\right)\d t.
\end{aligned}$$
which concludes the proof.
\end{proof}

\subsection{Value function convergence}

The energy balance of Proposition~\ref{prop:eps-energy-balance} is expressed in terms of the value function $V_\varepsilon$, whereas the desired EDI involves the original energy $\mathcal W$. We first show that $V_\varepsilon$ recovers $\mathcal W$ as $\varepsilon\to0$, including the lower bound needed when the current state also varies with $\varepsilon$. This handles the state-dependent energy terms in the limiting argument. The lower bound for the dual dissipation rate $G_\varepsilon$ requires a separate construction, developed in the next subsection.

\begin{proposition}[Convergence of $V_\varepsilon$]\label{prop:value-func-convergence}
It holds for any $T_0\in\statespace$ that $$V_\varepsilon(T_0)\nearrow \mathcal W(T_0),\quad \text{as }\varepsilon\to 0.$$
Moreover, also for $T_\varepsilon\in\statespace$ we have
$$\liminf_{\varepsilon\to 0}V_\varepsilon(T_\varepsilon)\geq \mathcal W(T_0),\quad \text{whenever }T_\varepsilon\weakstarto T_0 \text{ in }\intcur{k}(\R^d).$$
\end{proposition}
\begin{proof}
From the representation of $V_\varepsilon$ in Lemma \ref{lem:value-func-alt-representation} we readily see that  $\varepsilon \mapsto V_\varepsilon(T_0)$ is decreasing on $(0,\infty)$: $$0<\varepsilon_1\leq\varepsilon_2 \implies V_{\varepsilon_2}(T_0)\leq V_{\varepsilon_1}(T_0).$$
Moreover, since the constant trajectory $\eym{\llbracket [0,\infty)\rrbracket \times T_0}$ is a valid competitor, we have 
$$V_\varepsilon(T_0)\leq \mathcal W(T_0).$$

To complete the proof, we will show the lower estimate, namely $$\liminf_{\varepsilon\to 0}V_\varepsilon (T_\varepsilon)\geq \mathcal W(T_0).$$
Let $s>0$. Denoting $\mu_\varepsilon$ the $\varepsilon$-minimizer starting at $T_\varepsilon$, it follows from  the energy estimates (Proposition \ref{prop:energy-estimates}) and the H\"older continuity (Lemma \ref{lem:p-continuity-in-time}) that
\begin{align*}
\mathbf F(\stcurr{\mu_\varepsilon}(\varepsilon s)-T_\varepsilon)
&\leq\variation{\stcurr{\mu_\varepsilon};[0,\varepsilon s]}\\
&\leq\gamma^{1-1/p}\|v\|_{L^p((0,\varepsilon s);L^p(\mu_\varepsilon(\cdot)))}(\varepsilon s)^{1-1/p}\\
&\leq\gamma^{1-1/p}(pE_\gamma)^{1/p}(\varepsilon s)^{1-1/p}\to 0
\quad\text{as }\varepsilon\to0.
\end{align*}
Here \ref{item:E-bdd-by-mass} gives $\mathcal W(T_\varepsilon)\leq E_\gamma$. As weak${}^*$ convergence is equivalent to convergence in the flat distance on $\statespace$, we see that $$\mathbf F(\stcurr{\mu_\varepsilon}(\varepsilon s)-T_0)\leq \mathbf F(\stcurr{\mu_\varepsilon}(\varepsilon s)-T_\varepsilon)+\mathbf F(T_\varepsilon-T_0) \to 0, \quad \text{as }\varepsilon\to 0$$ and consequently $\stcurr{\mu_\varepsilon}(\varepsilon s)\weakstarto T_0$.

Hence, due to the lower semicontinuity of $\mathcal W$, it holds for all $s>0$ that $$\liminf_{\varepsilon\to 0}\mathcal W(\stcurr{\mu_\varepsilon}(\varepsilon s))\geq \mathcal W(T_0).$$
Now, we use the Fatou lemma to obtain 
$$\begin{aligned}
\liminf_{\varepsilon\to 0}V_\varepsilon (T_\varepsilon)& = \liminf_{\varepsilon\to 0}\int_0^\infty \frac{e^{-\frac t\varepsilon}}{\varepsilon}\left( \underbrace {\varepsilon\mathcal R_p(\slice{\mu_\varepsilon}{t})}_{\geq 0}+\mathcal W(\stcurr{\mu_\varepsilon}(t))\right)\d t\\
&\geq \liminf_{\varepsilon \to 0} \int_0^\infty e^{-s} \mathcal W(\stcurr{\mu_\varepsilon}(\varepsilon s))\d s
\\
&
\geq \int_0^\infty e^{-s}\liminf_{\varepsilon \to 0} \mathcal W(\stcurr{\mu_\varepsilon}(\varepsilon s))\d s\geq \int_0^\infty e^{-s}\mathcal W(T_0)\d s
=\mathcal W(T_0).
\end{aligned}$$
Altogether, we have proven that 
\begin{align*}V_\varepsilon(T_0)\nearrow \mathcal W(T_0),\quad &\text{as }\varepsilon\to 0. \qedhere\end{align*}
\end{proof}

\subsection{Moreau-Yosida Regularization and Slope Bounds}

Value-function convergence identifies the limiting energy, but it does not by itself give the required lower bound for the dual rate $G_\varepsilon$. We must relate that gap to the relaxed dual dissipation rate as $\varepsilon\to0$, allowing the initial current to vary as well. Our target is the following estimate.

\begin{proposition}[$\liminf$ estimate for $G_\varepsilon$]\label{prop:liminf-Gepsilon}
For any $T_\varepsilon,T_0\in\statespace$ with
$$
T_\varepsilon\weakstarto T_0 
$$
it holds
$$
\liminf_{\varepsilon\to 0} G_\varepsilon(T_\varepsilon)\geq \dualRinf{T_0}.
$$
\end{proposition}

The argument uses two auxiliary constructions. First, we construct a Moreau--Yosida-type regularization $\mathcal W_h$, analogous to the metric-space construction in \cite[Lemma 3.1.5]{ambrosioGradientFlowsMetric2008}. The limit $h\to0$ yields a bound involving the relaxed dual rate.

Second, at fixed $\varepsilon$, we introduce an auxiliary discrete WED functional with time step $\tau$. Short-time minimizers provide competitors for this discrete problem. Comparing it with the continuous WED functional as $\tau\to0$ gives a lower bound for $G_\varepsilon$ at fixed $\varepsilon$. We then use that bound in the passage $\varepsilon\to0$. The short-time regularization and the discrete approximation thus serve different purposes, and their roles will be kept separate below.

\begin{definition}[Moreau-Yosida regularization]\label{def:moreau-yosida-reg}
We define the regularization of $\mathcal W$ over trajectories as follows. For $h>0$ define 
$$\mathcal W_h(T_0):=\inf_{\mu\in\young(T_0)}\left\{\mathcal W(\stcurr{\mu}(h))+\int_0^h\mathcal R_p(\slice{\mu}{t})\d t\right\}, \quad T_0\in\statespace.$$
\end{definition}
\begin{lemma}[Moreau-Yosida minimizers]\label{lem:moreau-yosida-minimizer}
For every $h>0$ and $T_0\in\statespace$, there exists a minimizer $\mu_h\in\young(T_0)$, so that
\[
\mathcal W_h(T_0)=\mathcal W(\stcurr{\mu_h}(h))+\int_0^h\mathcal R_p(\slice{\mu_h}{t})\d t.
\]
The minimizer can be chosen stationary after $h$, meaning that its restriction to $(h,\infty)$ is the stationary elementary varicurrent
$$
\mu_h\mres (h,\infty)
=\eym{\llbracket(h,\infty)\rrbracket\times \stcurr{\mu_h}(h)}.
$$ Moreover, any such minimizer $\mu_h$ satisfies
\begin{equation}\label{eq:moreau-yosida-minimizer-bounds}
0\leq\mathcal W_h(T_0)\leq\mathcal W(T_0)\leq E_\gamma,
\qquad
\int_0^h\mathcal R_p(\slice{\mu_h}{t})\d t\leq\mathcal W(T_0).
\end{equation}
\end{lemma}
\begin{proof}
Choose an infimizing sequence $(\mu_i)\subset\young(T_0)$. We replace $\mu_i$ with a stationary path after time $h$, namely we define
$$
\tilde \mu_i:=\mu_i\mres [0,h]
+\eym{\llbracket(h,\infty)\rrbracket\times \stcurr{\mu_i}(h)}.
$$
Then we have
$$
\mathcal W(\stcurr{\tilde\mu_i}(h))+\int_0^h\mathcal R_p(\slice{\tilde\mu_i}{t})\d t
=\mathcal W(\stcurr{\mu_i}(h))+\int_0^h\mathcal R_p(\slice{\mu_i}{t})\d t
$$
and hence the modified sequence $\tilde\mu_i$ is also infimizing. As $\mathcal R_p(\slice{\tilde\mu_i}{t})=0$ for $t>h$, we may use Proposition \ref{prop:compactness-dissipation-potential} and the energy bound (Proposition \ref{prop:energy-estimates}) to obtain a subsequence with
$$\tilde\mu_i\weakstarto\mu_h$$ for some $\mu_h\in\young(T_0)$, and, in particular, by Lemma \ref{lem:all-time-slices-Y} also 
$$
\stcurr{\mu_i}(h) \weakstarto \stcurr{\mu_h}(h)\quad\text{in }\intcur{k}(\overline\Omega).
$$
Invoking the lower semicontinuity of $\mathcal W$ in \ref{item:E-wlsc} and of the $\mathcal R_p$-time integral in Lemma \ref{lem:wlsc-dissipation} finishes the proof.
\end{proof}
We next control the dependence of $\mathcal W_h$ on the time horizon $h$, with the current state held fixed.
\begin{lemma}[Lipschitz continuity of $\mathcal W_h$]\label{lem:yosida-is-lipschitz}
For any $T_0\in\statespace$ it holds that 
$$
h\mapsto \mathcal W_h(T_0)\quad \text{ is nonincreasing and locally Lipschitz on }(0,\infty).
$$
Moreover,
$$
\lim_{h\to0}\mathcal W_h(T_0)=\mathcal W(T_0),
$$
and consequently, 
$$
\mathcal W_h(T_0)-\mathcal W (T_0) = \int_0^h \frac{\mathrm d}{\mathrm dr}\mathcal W_r(T_0)\d r.
$$
\end{lemma}
\begin{proof}
Let $0<r<s$. By Lemma \ref{lem:moreau-yosida-minimizer}, choose a minimizer $\mu_r\in\young(T_0)$ for $\mathcal W_r(T_0)$ that is stationary after $r$. Using this same trajectory as a competitor for $\mathcal W_s(T_0)$, we see that
$$
\mathcal W_s(T_0)\leq  \mathcal W(\stcurr{\mu_r}(r))+\int_0^r\mathcal R_p(\slice{\mu_r}{t})\d t=\mathcal W_r(T_0)
$$
which proves the monotonicity.
Now, conversely, choose an exact minimizer $\mu_s\in\young(T_0)$ for $\mathcal W_s(T_0)$ by Lemma \ref{lem:moreau-yosida-minimizer} and construct
$$
\tilde \mu_r=\timepush{\left(t\mapsto \frac rs t\right)}{\mu_s}.
$$
Proposition \ref{prop:time-pushforward}, Lemma \ref{lem:dissipation-pushforward}, and the energy bound show that $\tilde\mu_r\in\young(T_0)$. By Lemma \ref{lem:dissipation-pushforward}, we have
$$
\int_0^r\mathcal R_p(\slice{\tilde\mu_r}{t})\d t=\left(\frac s r\right)^{p-1}\int_0^s\mathcal R_p(\slice{\mu_s}{t})\d t.
$$
Using this and $\tilde\mu_r$ as a competitor for $\mathcal W_r(T_0)$, we obtain

\[
\mathcal W_r(T_0)
 \leq \mathcal W(\stcurr{\mu_s}(s))+\left(\frac sr\right)^{p-1}\int_0^s\mathcal R_p(\slice{\mu_s}{t})\d t
 =\mathcal W_s(T_0)+ \left( \left(\frac sr\right)^{p-1}-1\right)\int_0^s\mathcal R_p(\slice{\mu_s}{t})\d t.
\]
It holds that 
$$
s\mapsto \left(\frac sr\right)^{p-1}-1= \frac{s^{p-1}-r^{p-1}}{r^{p-1}} \text{ is locally Lipschitz on a neighborhood of }r
$$
and, by \eqref{eq:moreau-yosida-minimizer-bounds},
$$
\int_0^s\mathcal R_p(\slice{\mu_s}{t})\d t\leq \mathcal W(T_0).
$$
So, altogether, it holds
$$
0\leq\mathcal W_r(T_0)- \mathcal W_s(T_0)\leq  C(r_0,s_0)\mathcal W(T_0)|s-r| \quad \text{for } 0<r_0\leq r <s \leq s_0<\infty.
$$
This proves the local Lipschitz continuity in $(0,\infty)$.

It remains to identify the limit at zero. Choose a minimizer \(\mu_h\) for
\(\mathcal W_h(T_0)\) by Lemma \ref{lem:moreau-yosida-minimizer}. By \eqref{eq:moreau-yosida-minimizer-bounds},
\[
\mathcal W(\stcurr{\mu_h}(h))+\int_0^h\mathcal R_p(\slice{\mu_h}{t})\d t\leq\mathcal W(T_0).
\]
Using Lemma \ref{lem:all-time-slices-Y}, we have
\[
\mathbf F(\stcurr{\mu_h}(h)-T_0)
\leq
\bigl(p\mathcal W(T_0)\bigr)^{1/p}(\gamma h)^{1-1/p}
\leq
(pE_\gamma)^{1/p}
(\gamma h)^{1-1/p}
\longrightarrow0.
\]
By the lower semicontinuity of \(\mathcal W\),
\[
\mathcal W(T_0)
\leq
\liminf_{h\downarrow0}\mathcal W(\stcurr{\mu_h}(h))
\leq
\liminf_{h\downarrow0}\mathcal W_h(T_0).
\]
The reverse inequality follows from \(\mathcal W_h(T_0)\leq\mathcal W(T_0)\),
which proves $\lim_{h\to 0}\mathcal W_h(T_0)=\mathcal W(T_0)$.

Finally, local Lipschitz continuity gives
\[
\mathcal W_h(T_0)-\mathcal W_a(T_0)
=
\int_a^h\frac{\d}{\d r}\mathcal W_r(T_0)\,\d r
\qquad (0<a<h).
\]
Utilizing the non-positivity of the derivative and the continuity at 0, we conclude the proof by passing to \(a\to 0\).
\end{proof}
The next estimate connects variation of the time horizon to the dual dissipation rate.
\begin{lemma}\label{lem:yosida-derivative-geq-slope}
Let $T_0\in\statespace$. At every $r>0$ where $r\mapsto\mathcal W_r(T_0)$ is differentiable, and for every minimizer $\mu_r\in\young(T_0)$ supplied by Lemma \ref{lem:moreau-yosida-minimizer}, it holds
$$
-\frac{\mathrm d}{\mathrm dr}\mathcal W_r(T_0)\geq \dualR{\stcurr{\mu_r}(r)}.
$$
\end{lemma}

\begin{proof}
Fix such an $r>0$ and any minimizer $\mu_r$. For $\mu\in\young(\stcurr{\mu_r}(r))$, form the space-time varicurrent by continuing $\mu_r$ by $\mu$ after time $r$, that is,
$$
\tilde\mu=\mu_r\mres [0,r)
+\timepush{(t\mapsto t+r)}{\mu}.
$$
It satisfies $\tilde\mu\in\young(T_0)$. Using it as a competitor in $\mathcal W_{r+s}(T_0)$ for $s>0$ gives
\begin{align*}
\mathcal W_{r+s}(T_0)
&\leq\mathcal W(\stcurr{\tilde\mu}(r+s))+\int_0^{r+s}\mathcal R_p(\slice{\tilde\mu}{t})\d t
\\
&=\mathcal W(\stcurr{\mu}(s))+\int_0^r\mathcal R_p(\slice{\mu_r}{t})\d t+\int_0^s\mathcal R_p(\slice{\mu}{t})\d t.
\end{align*}
We observe that 
$$
\int_0^r\mathcal R_p(\slice{\mu_r}{t})\d t=\mathcal W_r(T_0)-\mathcal W(\stcurr{\mu_r}(r)).
$$
Plugging this above and rearranging the terms yields
\begin{align*}
\mathcal W_r(T_0)-\mathcal W_{r+s}(T_0)
&\geq\mathcal W(\stcurr{\mu_r}(r))-\mathcal W(\stcurr{\mu}(s))-\int_0^s\mathcal R_p(\slice{\mu}{t})\d t.
\end{align*}
Dividing by $s$ and taking the $\limsup$ as $s\to 0$, differentiability at $r$ gives
\begin{align*}
-\frac{\mathrm d}{\mathrm dr}\mathcal W_r(T_0)
\geq\limsup_{s\to 0}\left[\frac{\mathcal W(\stcurr{\mu_r}(r))-\mathcal W(\stcurr{\mu}(s))}{s}-\frac1s\int_0^s\mathcal R_p(\slice{\mu}{t})\d t\right].
\end{align*}
Taking the supremum over all $\mu\in\young(\stcurr{\mu_r}(r))$ proves the claim for the arbitrary minimizer $\mu_r$.
\end{proof}

To connect the preceding estimate with the limiting state, we now let the horizon tend to zero. The resulting bound involves the relaxed dual rate, which is the quantity used in the EDI.

\begin{lemma}[$\Gamma$-$\liminf$ for the relaxed slope]\label{lem:yosida-relaxed-slope}
For all $T_h,T_0\in\statespace$ with
$$
T_h\weakstarto T_0\quad \text{in }\intcur{k}(\overline\Omega)
$$
it holds that
$$\liminf_{h \to 0}\frac{\mathcal W(T_h)-\mathcal W_h(T_h)}{h}\geq \dualRinf{T_0}.$$
\end{lemma}
\begin{proof}
For fixed $h>0$, define the infimum
\[
c_h:=\inf\left\{\dualR{\stcurr{\nu}(s)}:
s\in(0,h),\ \nu\in\young(T_h),
\displaystyle\mathcal W_s(T_h)=\mathcal W(\stcurr{\nu}(s))+\int_0^s\mathcal R_p(\slice{\nu}{t})\d t
\right\}.
\]
Lemma \ref{lem:moreau-yosida-minimizer} ensures that the set is nonempty. The stationary competitor in the definition of $\dualR{\cdot}$ gives $c_h\geq0$. At every differentiability point $s\in(0,h)$, Lemma \ref{lem:yosida-derivative-geq-slope} applies to every minimizer, and hence
\[
-\frac{\mathrm d}{\mathrm ds}\mathcal W_s(T_h)\geq c_h\quad \text{for a.e. }s\in(0,h).
\]
The derivative exists and is finite almost everywhere by Lemma \ref{lem:yosida-is-lipschitz}, so $c_h<\infty$. Integrating this and using the same lemma yields
\[
\frac{\mathcal W(T_h)-\mathcal W_h(T_h)}{h}
=\frac1h\int_0^h-\frac{\mathrm d}{\mathrm ds}\mathcal W_s(T_h)\d s\geq c_h.
\]

For every endpoint occurring in the infimum, Lemma \ref{lem:all-time-slices-Y} and \eqref{eq:moreau-yosida-minimizer-bounds} give the uniform bound
\[
\mathbf F\bigl(\stcurr{\nu}(s)-T_h\bigr)
\leq(p\mathcal W(T_h))^{1/p}(\gamma s)^{1-1/p}
\leq(pE_\gamma)^{1/p}(\gamma h)^{1-1/p}.
\]
It remains to prove $$\liminf_{h\downarrow0}c_h\geq\dualRinf{T_0}.$$ Choose $h_j\to0$ which realizes the $\liminf$, and choose $s_j\in(0,h_j)$ and a minimizer $\nu_j\in\young(T_{h_j})$ such that
\[
\dualR{\stcurr{\nu_j}(s_j)}\leq c_{h_j}+\frac1j.
\]
The endpoint bound and $T_{h_j}\weakstarto T_0$ imply $\stcurr{\nu_j}(s_j)\weakstarto T_0$, since flat and weak${}^*$ convergence agree on $\statespace$ (see Section~\ref{sec:admissible-trajectories}). By the definition of the relaxed slope,
\[
\dualRinf{T_0}\leq\liminf_{j\to\infty}\dualR{\stcurr{\nu_j}(s_j)}
\leq\lim_{j\to\infty}c_{h_j}
\]
which concludes the proof.
\end{proof}

We now fix $\varepsilon>0$ and introduce an auxiliary time step $\tau$. The next lemma compares the discrete WED functional with the continuous one as $\tau\to0$. This limit is taken before the subsequent passage $\varepsilon\to0$.

\begin{lemma}[$\Gamma$-convergence of discrete WED functional]\label{lem:gamma-conv-discr-WED}
Let $\varepsilon>0$ be fixed.
For a time-step parameter $\tau>0$, we define the discrete weights $w_n^\tau$ and the piecewise constant weight $\bar w^\tau\colon (0,\infty) \to (0,\infty)$ by
\begin{equation}
w_n^\tau := \int_{(n-1)\tau}^{n\tau} \frac{e^{-\frac t \varepsilon}}{\varepsilon}\d t = e^{-\tfrac{(n-1)\tau}{\varepsilon}}-e^{-\tfrac{n\tau}{\varepsilon}}, \quad \text{and} \quad \bar w^\tau(t) := \frac{1}{\tau} w_n^\tau \quad \text{for } t \in [(n-1)\tau, n\tau).
\end{equation}
We define the discrete WED functional $\mathcal J_\varepsilon^\tau\colon\young\to[0,\infty]$
\begin{equation}
\begin{aligned}
\mathcal J_\varepsilon^\tau(\mu) &:= \int_0^\infty \bar w^\tau(t) \left( \mathcal W(\bar S_\tau(t)) + \varepsilon \mathcal R_p(\slice{\mu}{t}) \right) \d t\\
&= \sum_{n=1}^\infty w_n^\tau \left( \mathcal W(\stcurr{\mu}(n\tau)) + \varepsilon \frac1\tau \int_{(n-1)\tau}^{n\tau}\mathcal R_p(\slice{\mu}{t})\d t \right),
\end{aligned}
\end{equation}
where for $\mu\in\young$ we denoted the piecewise-constant interpolation of $\stcurr{\mu}$ by $\bar S_\tau(t):=\stcurr{\mu}(n\tau)$ for $t\in[(n-1)\tau,n\tau)$, $n\in\N$.

Then it holds 
\begin{equation}
    \mathcal J_\varepsilon^\tau \overset{\Gamma}{\to} \mathcal J_\varepsilon\quad\text{on $\young(T_0)$ as } \tau\to 0.
\end{equation}
Specifically, this means that for every $\mu\in\young$ and $\mu_\tau\in\young(\stcurr{\mu}(0))$ with  $\mu_\tau\weakstarto\mu$ it holds
$$
\liminf_{\tau\to 0}\mathcal J_\varepsilon^\tau(\mu_\tau)\geq \mathcal J_\varepsilon(\mu),
$$
and for every $\mu\in\young$ there exists a recovery sequence $\mu_\tau\in\young(\stcurr{\mu}(0))$ satisfying $\mu_\tau\weakstarto\mu$ with
$$
\limsup_{\tau\to 0}\mathcal J_\varepsilon^\tau(\mu_\tau)\leq \mathcal J_\varepsilon(\mu).
$$
\end{lemma}
\begin{proof}
\textit{Recovery sequence ($\Gamma$-$\limsup$):}

Let $\mu\in\young$. As the energy is evaluated only at discrete times
$n\tau$, the constant recovery sequence may fail, as it may persistently
sample points where $\mathcal W(\stcurr{\mu}(n\tau))$ is too large.
We overcome this by defining the recovery sequence by shifting forward
in time by some $0<\delta_\tau<\tau$, extending it by the initial condition
on $(0,\delta_\tau)$. Namely, we define $\mu_\tau\in\young(\stcurr{\mu}(0))$ by
$$\mu_\tau
:=
\eym{\llbracket(0,\delta_\tau)\rrbracket
      \times \stcurr{\mu}(0)}
+
\timepush{(t\mapsto t+\delta_\tau)}{\mu}.
$$
Moreover, since $\delta_\tau\to0$, the translated measures converge
locally weakly$^*$ to $\mu$, while the mass of the prefix is
$\delta_\tau\mass{\stcurr{\mu}(0)}\to0$.
Consequently, $\mu_\tau\weakstarto\mu$.

To choose $\delta_\tau$, compute the value of the energy part of
$\mathcal J_\varepsilon^\tau$ averaged over all possible shifts
$\delta\in(0,\tau)$:
$$
A_\tau:= \frac 1\tau\int_0^\tau\sum_{n=1}^\infty w_n^\tau \mathcal W(\stcurr{\mu}(n\tau-\delta)) \d\delta = \frac 1\tau\sum_{n=1}^\infty w_n^\tau \int_{(n-1)\tau}^{n\tau} \mathcal W(\stcurr{\mu}(s)) \d s.
$$
For any $s \in ((n-1)\tau,n\tau)$, we have that
$$
w_n^\tau = e^{-\frac{(n-1)\tau}{\varepsilon}}(1-e^{-\frac\tau\varepsilon})\leq \frac \tau\varepsilon e^{-\frac{(n-1)\tau}{\varepsilon}} \leq \frac\tau\varepsilon e^{\frac{\tau}{\varepsilon}} e^{-\frac{s}{\varepsilon}}.
$$
Substituting this upper bound into the expression for $A_\tau$ yields
$$
A_\tau \leq e^{\frac{\tau}{\varepsilon}} \sum_{n=1}^\infty \int_{(n-1)\tau}^{n\tau} \frac{e^{-\frac s\varepsilon}}{\varepsilon}\mathcal W(\stcurr{\mu}(s)) \d s
= e^{\frac{\tau}{\varepsilon}} \int_0^\infty \frac{e^{-\frac s\varepsilon}}{\varepsilon}\mathcal W(\stcurr{\mu}(s)) \d s.
$$
Hence there exists a $\delta_\tau\in(0,\tau)$ such that
$$
\sum_{n=1}^\infty w_n^\tau \mathcal W(\stcurr{\mu}(n\tau-\delta_\tau)) \leq e^{\frac{\tau}{\varepsilon}} \int_0^\infty \frac{e^{-\frac s \varepsilon}}{\varepsilon}\mathcal W(\stcurr{\mu}(s)) \d s
$$
and consequently,
$$
\limsup_{\tau\to 0}\sum_{n=1}^\infty w_n^\tau
\mathcal W(\stcurr{\mu_\tau}(n\tau))
\leq \int_0^\infty \frac{e^{-\frac s \varepsilon}}{\varepsilon}
\mathcal W(\stcurr{\mu}(s)) \d s.
$$

We now verify the convergence of the $\mathcal R_p$ term.
Denoting the piecewise constant weight
$$
\bar w^\tau(t) = \frac 1\tau w_n^\tau,\quad t\in [(n-1)\tau,n\tau),
$$
we write
$$
\sum_{n=1}^\infty w_n^\tau \frac \varepsilon\tau\int_{(n-1)\tau}^{n\tau}\mathcal R_p(\slice{\mu_\tau}{t})\d t
= \varepsilon\int_0^\infty \bar w^\tau(s+\delta_\tau)\mathcal R_p(\slice{\mu}{s})\d s.
$$
As it holds
$$
\lim_{\tau \to 0} \bar w^\tau(s+\delta_\tau)
= \frac{e^{-\frac s\varepsilon}}{\varepsilon},
\qquad
0\leq \bar w^\tau(s+\delta_\tau)
\leq e^{\frac\tau\varepsilon}\frac{e^{-\frac s\varepsilon}}{\varepsilon},
$$
and $\mathcal J_\varepsilon(\mu)<\infty$ gives
$\int_0^\infty e^{-s/\varepsilon}\mathcal R_p(\slice{\mu}{s})\d s<\infty$,
we get by the Lebesgue dominated convergence theorem that
$$
\sum_{n=1}^\infty w_n^\tau \frac \varepsilon\tau
\int_{(n-1)\tau}^{n\tau}\mathcal R_p(\slice{\mu_\tau}{t})\d t
\to \varepsilon\int_0^\infty
\frac{e^{-\frac t\varepsilon}}{\varepsilon}
\mathcal R_p(\slice{\mu}{t})\d t,
$$
altogether showing the desired $\limsup$-inequality
$$
\limsup_{\tau\to 0}\mathcal J_\varepsilon^\tau(\mu_\tau)
\leq \mathcal J_\varepsilon(\mu).
$$

\textit{Lower semicontinuity ($\Gamma$-$\liminf$):}
Let $(\mu_\tau)_{\tau>0}\subset \young$, $\mu_\tau\weakstarto\mu$ and assume that 
$$
\liminf_{\tau\to0}\mathcal J_\varepsilon^\tau(\mu_\tau)=:m<\infty,
$$
as otherwise there is nothing to prove.

Let us consider the piecewise constant interpolation $\bar S_\tau\in L^\infty_{w^*}((0,\infty);\intcur{k}(\overline\Omega))$ by setting
$$
\bar S_\tau(t)=\stcurr{\mu_\tau}(n\tau), \quad t\in[(n-1)\tau,n\tau), n\in\N.
$$
We will now show that 
$$
\bar S_\tau(t)\weakstarto \stcurr{\mu}(t),\quad \text{for a.e.\ }t\in(0,\infty).
$$
We compute for $t\in((n-1)\tau,n\tau)$ by means of Lemma \ref{lem:variation-ineq} and Lemma \ref{lem:p-continuity-in-time}
\begin{align*}
\mathbf F(\stcurr{\mu_\tau}(t)-\stcurr{\mu_\tau}(n\tau))\leq \variation{\stcurr{\mu_\tau};[t,n\tau]}\leq\variation{\mu_\tau;[t,n\tau]}
\\
\leq \|v\|_{L^p(((n-1)\tau,n\tau);L^p(\mu_\tau(\cdot)))} \gamma^{1-1/p} |n\tau-t|^{1-1/p}\leq \left(p e^{\frac{\tau n}{\varepsilon}}\mathcal J_\varepsilon^\tau(\mu_\tau)\right)^{1/p}\gamma^{1-1/p}\tau^{1-1/p}.
\end{align*}
As the exponential term is bounded by $e^{\frac{\tau n}{\varepsilon}}\leq e^{\frac{t+\tau}{\varepsilon}}$ and thus bounded uniformly in $\tau$, this shows that 
$$
\mathbf F(\stcurr{\mu_\tau}(t)-\overline S_\tau(t))=
\mathbf F(\stcurr{\mu_\tau}(t)-\stcurr{\mu_\tau}(n\tau)) \to 0 \quad\text{ for all }t\in[0,\infty).
$$
Together with $\stcurr{\mu_\tau}(t)\weakstarto\stcurr{\mu}(t)$ for every $t$,
this implies 
$$
\overline S_\tau(t) \weakstarto\stcurr{\mu}(t) \quad \text{in $\intcur{k}(\overline\Omega)$ for a.e.\ } t\in(0,\infty).
$$
We use this, combined with the Fatou lemma and weak${}^*$ lower semicontinuity of $\mathcal W$ to obtain 
$$
\liminf_{\tau\to0}\sum_{n=1}^\infty w_n^\tau\mathcal W(\stcurr{\mu_\tau}(n\tau))= \liminf_{\tau\to 0}\int_0^\infty\bar w^\tau (t) \mathcal W(\bar S_\tau (t))\d t \geq \int_0^\infty \frac{e^{-\frac t\varepsilon}}{\varepsilon}\mathcal W(\stcurr{\mu}(t))\d t.
$$
For the dissipation term, we remark that $\bar w^\tau(t)\geq e^{-\frac \tau\varepsilon}\frac{e^{-\frac t\varepsilon}}{\varepsilon}$ and hence by the weighted conclusion of Proposition \ref{prop:compactness-dissipation-potential}, it holds
\begin{align*}
\liminf_{\tau\to 0}\sum_{n=1}^\infty w_n^\tau\varepsilon \frac1\tau \int_{(n-1)\tau}^{n\tau}\mathcal R_p(\slice{\mu_\tau}{t})\d t\geq
\liminf_{\tau\to 0} e^{-\frac \tau\varepsilon}\int_0^\infty e^{-\frac t\varepsilon}\mathcal R_p(\slice{\mu_\tau}{t})\d t 
\geq
\int_0^\infty e^{-\frac t\varepsilon}\mathcal R_p(\slice{\mu}{t})\d t,
\end{align*}
which completes the proof of the $\liminf$ inequality.
\end{proof}

We now use the short-time minimizers to build a minimizing-movement competitor $\tilde\mu_\varepsilon^\tau$ for the discrete functional $\mathcal J_\varepsilon^\tau$. For fixed $\varepsilon$, these competitors relate the regularized energy decrease to the WED value function. The $\Gamma$-convergence of Lemma~\ref{lem:gamma-conv-discr-WED} allows us to pass to the continuous WED functional as $\tau\to0$.

\begin{lemma}[Fixed-\texorpdfstring{$\varepsilon$}{epsilon} lower bound]
\label{lem:fixed-epsilon-G-lower-bound}
Let \(\varepsilon>0\) be fixed  and $T_\varepsilon\in\statespace$ with $\mathcal W(T_\varepsilon)<\infty$. Then
there exists a trajectory \(\mu_{\varepsilon,T_\varepsilon}\in\young(T_\varepsilon)\) such that
\begin{equation}\label{eq:fixed-epsilon-G-lower-bound}
G_\varepsilon(T_\varepsilon)
\geq
\int_0^\infty
\frac{e^{-t/\varepsilon}}{\varepsilon}
\dualRinf{\stcurr{\mu_{\varepsilon,T_\varepsilon}}(t)}\,\d t,
\end{equation}
and
\begin{equation}\label{eq:fixed-epsilon-trajectory-action}
\int_0^\infty
\mathcal R_p(\slice{\mu_{\varepsilon,T_\varepsilon}}{t})\,\d t
\leq
\mathcal W(T_\varepsilon).
\end{equation}
\end{lemma}

\begin{proof}
Let us consider the minimizing movements-like scheme for a time step $\tau>0$, that is, we define
$$
\mu_{\varepsilon,n}^\tau\in\operatorname*{argmin}\left\{ \mathcal W(\stcurr{\mu}(\tau))+\int_0^\tau\mathcal R_p(\slice{\mu}{t})\d t: \mu\in\young\left(\stcurr{\mu_{\varepsilon,n-1}^\tau}(\tau)\right)\right\}
$$
recursively over $n\in\N$, with initial state $T_\varepsilon$ when $n=1$. Each terminal current belongs to $\statespace$ by Lemma \ref{lem:all-time-slices-Y}, so Lemma \ref{lem:moreau-yosida-minimizer} supplies a minimizer at every step. Define the space-time varicurrent
$$
\tilde\mu_\varepsilon^\tau
:=\sum_{n=1}^\infty\timepush{(t\mapsto t+(n-1)\tau)}
{\left(\mu_{\varepsilon,n}^\tau\mres [0,\tau)\right)}.
$$
Due to the condition $\stcurr{\mu_{\varepsilon,n}^\tau}(0)=\stcurr{\mu_{\varepsilon,n-1}^\tau}(\tau)$, we see that $\tilde\mu_\varepsilon^\tau\in\young(T_\varepsilon)$.
In what follows, we will show that 
$$
\tilde \mu_\varepsilon^\tau\weakstarto\mu_\varepsilon \quad \text{as }\tau\to0
$$
for some $\mu_\varepsilon\in\young(T_\varepsilon)$ with the desired properties.

Lemma \ref{lem:moreau-yosida-minimizer} gives
\[
\mathcal W(\stcurr{\tilde\mu_\varepsilon^\tau}(n\tau))+\int_{(n-1)\tau}^{n\tau}\mathcal R_p(\slice{\tilde\mu_\varepsilon^\tau}{t})\d t
=\mathcal W_\tau(\stcurr{\tilde\mu_\varepsilon^\tau}((n-1)\tau)) \leq\mathcal W(\stcurr{\tilde\mu_\varepsilon^\tau}((n-1)\tau)).
\]
Summing over $n=1,\ldots,N$ gives
$$
\mathcal W(\stcurr{\tilde\mu_\varepsilon^\tau}(N\tau))
+\int_0^{N\tau}\mathcal R_p(\slice{\tilde\mu_\varepsilon^\tau}{t})\d t\leq\mathcal W(T_\varepsilon).
$$
Nonnegativity and monotone convergence as $N\to\infty$ imply
$$
\int_0^\infty \mathcal R_p(\slice{\tilde\mu_\varepsilon^{\tau}}{t})\d t \leq \mathcal W(T_\varepsilon).
$$
Lemmas \ref{lem:variation-ineq} and \ref{lem:p-continuity-in-time} give
$$
\variation{\stcurr{\tilde\mu_\varepsilon^\tau};[s,t]}
\leq(p\mathcal W(T_\varepsilon))^{1/p}\gamma^{1-1/p}|t-s|^{1-1/p}.
$$
Let us denote
$$E_\varepsilon^{\tau,n}:=\mathcal W(\stcurr{\tilde\mu_\varepsilon^\tau}(n\tau)),\quad R_\varepsilon^{\tau,n}:=\int_{(n-1)\tau}^{n\tau}\int\tfrac1p|v|^p \d\slice{\tilde\mu_\varepsilon^\tau}{t}\d t.$$
Here $E_\varepsilon^{\tau,0}=\mathcal W(T_\varepsilon)$.
We rewrite
$$
\sum_{n=1}^\infty\left(e^{-\tfrac{(n-1)\tau}{\varepsilon}}-e^{-\frac{n\tau}{\varepsilon}}\right)E_\varepsilon^{\tau,n}
= E_\varepsilon^{\tau,0}+\sum_{n=1}^\infty e^{-\frac{(n-1)\tau}{\varepsilon}}\left(E_\varepsilon^{\tau,n}-E_\varepsilon^{\tau,n-1}\right)
$$
and
$$
e^{-\tfrac{(n-1)\tau}{\varepsilon}}-e^{-\frac{n\tau}{\varepsilon}}= e^{-\tfrac{(n-1)\tau}{\varepsilon}}\left(1-e^{-\frac{\tau}{\varepsilon}}\right),
$$
to obtain
\begin{align*}
\mathcal J_\varepsilon^\tau(\tilde\mu_\varepsilon^\tau)
=&E_\varepsilon^{\tau,0}+\sum_{n=1}^\infty e^{-\frac{(n-1)\tau}{\varepsilon}}\biggl[E_\varepsilon^{\tau,n}-E_\varepsilon^{\tau,n-1}
+\varepsilon \frac{\left(1-e^{-\frac{\tau}{\varepsilon}}\right)}{\tau} R_\varepsilon^{\tau,n}\biggr].
\end{align*}
By Lemma \ref{lem:moreau-yosida-minimizer},
$$
\mathcal W_\tau(\stcurr{\tilde\mu_\varepsilon^\tau}((n-1)\tau))=E_\varepsilon^{\tau,n}+R_\varepsilon^{\tau,n},
$$
we obtain
$$
\mathcal J_\varepsilon^\tau(\tilde\mu_\varepsilon^\tau)
=E_\varepsilon^{\tau,0}+\sum_{n=1}^\infty e^{-\frac{(n-1)\tau}{\varepsilon}}\biggl[\mathcal W_\tau(\stcurr{\tilde\mu_\varepsilon^\tau}((n-1)\tau))-E_\varepsilon^{\tau,n-1}
+\left(\varepsilon \tfrac{\left(1-e^{-\frac{\tau}{\varepsilon}}\right)}{\tau} -1\right)R_\varepsilon^{\tau,n}\biggr].
$$
Hence
\begin{align*}
&\frac{\mathcal W(T_\varepsilon)-\inf_{\mu\in\young(T_\varepsilon)}\mathcal J_\varepsilon^\tau(\mu)}{\varepsilon}\geq \frac{E_\varepsilon^{\tau,0}-\mathcal J_\varepsilon^\tau(\tilde \mu_\varepsilon^\tau)}{\varepsilon}\\
&\quad= \sum_{n=1}^\infty \frac{e^{-\frac{(n-1)\tau}{\varepsilon}}}{\varepsilon}\biggl[
E_\varepsilon^{\tau,n-1}-\mathcal W_\tau(\stcurr{\tilde\mu_\varepsilon^\tau}((n-1)\tau))
-\left(\varepsilon \tfrac{\left(1-e^{-\frac{\tau}{\varepsilon}}\right)}{\tau} -1\right)R_\varepsilon^{\tau,n}
\biggr].
\end{align*}
Using the identity
$$
e^{-\frac{(n-1)\tau}{\varepsilon}}=\frac{e^{-\frac{(n-1)\tau}{\varepsilon}}-e^{-\frac{n\tau}{\varepsilon}}}{1-e^{-\frac\tau\varepsilon}}
$$
we arrive at
\begin{align*}
\frac{\mathcal W(T_\varepsilon)-\inf_{\mu\in\young(T_\varepsilon)}\mathcal J_\varepsilon^\tau(\mu)}{\varepsilon}
=&
\frac{\tau}{\varepsilon(1-e^{-\frac\tau\varepsilon})}\sum_{n=1}^\infty \left(e^{-\frac{(n-1)\tau}{\varepsilon}}-e^{-\frac{n\tau}{\varepsilon}}\right)\frac{E_\varepsilon^{\tau,n-1}-\mathcal W_\tau(\stcurr{\tilde\mu_\varepsilon^\tau}((n-1)\tau))}{\tau}
\\
&+\underbrace{\left(1-\varepsilon \tfrac{\left(1-e^{-\frac{\tau}{\varepsilon}}\right)}{\tau} \right)}_{\geq 0} \sum_{n=1}^\infty \frac{e^{-\frac{(n-1)\tau}{\varepsilon}}}{\varepsilon} R_\varepsilon^{\tau,n}.
\end{align*}
Now by $\Gamma$-convergence of the discrete WED functional (Lemma \ref{lem:gamma-conv-discr-WED}) we have 
$$
\lim_{\tau\to0} \inf_{\mu\in\young(T_\varepsilon)}\mathcal J_\varepsilon^\tau(\mu) = V_\varepsilon(T_\varepsilon),
$$
and the limit
$$
\lim_{\tau\to 0} \frac{1-e^{-\frac\tau\varepsilon}}{\tau}=\frac1\varepsilon
$$
holds, so that we obtain
\begin{align*}
\frac{\mathcal W(T_\varepsilon)-V_\varepsilon(T_\varepsilon)}{\varepsilon}
\geq&
\liminf_{\tau\to0}\frac{\tau}{\varepsilon(1-e^{-\frac\tau\varepsilon})}\sum_{n=1}^\infty \left(e^{-\frac{(n-1)\tau}{\varepsilon}}-e^{-\frac{n\tau}{\varepsilon}}\right)\frac{E_\varepsilon^{\tau,n-1}-\mathcal W_\tau(\stcurr{\tilde\mu_\varepsilon^\tau}((n-1)\tau))}{\tau}
\\
&+
\liminf_{\tau\to0}\underbrace{\left(1-\varepsilon \frac{\left(1-e^{-\frac{\tau}{\varepsilon}}\right)}{\tau} \right)}_{\geq 0} \sum_{n=1}^\infty \frac{e^{-\frac{(n-1)\tau}{\varepsilon}}}{\varepsilon} R_\varepsilon^{\tau,n}.
\end{align*}
Denoting the piecewise-constant current $\bar S_\varepsilon^\tau\in L^\infty_{w^*}((0,\infty);\intcur{k}(\overline\Omega))$
$$
\bar S_\varepsilon^\tau(t)=\stcurr{\tilde\mu_\varepsilon^\tau}((n-1)\tau),\quad t\in((n-1)\tau,n\tau)
$$
we can rewrite 
$$
\frac{\mathcal W(T_\varepsilon)-V_\varepsilon(T_\varepsilon)}{\varepsilon}
\geq\liminf_{\tau\to0}\frac{\tau}{\varepsilon(1-e^{-\frac\tau\varepsilon})}\int_0^\infty \bar w^\tau(t)\frac{\mathcal W(\bar S_\varepsilon^{\tau}(t))-\mathcal W_\tau(\bar S_\varepsilon^\tau(t))}{\tau}\d t.
$$
We take a sequence of time steps realizing the lower limit on the right-hand side of the preceding inequality when it is finite; if it is infinite, any sequence tending to zero suffices. All limits below are along this sequence or a subsequence of it. Using the limits
$$
\lim_{\tau\to 0}\frac{\tau}{1-e^{-\frac\tau\varepsilon}}=\varepsilon,\quad \text{and}\quad \lim_{\tau\to0}\bar w^\tau(t)= \lim_{\tau\to0}  e^{-\frac{\lfloor \frac{t}{\tau} \rfloor \tau}{\varepsilon}} \frac{1 - e^{-\frac{\tau}{\varepsilon}}}{\tau}= \frac{e^{-\frac t\varepsilon}}{\varepsilon}, \quad t\in(0,\infty)$$
and the Fatou lemma, we get
$$
\frac{\mathcal W(T_\varepsilon)-V_\varepsilon(T_\varepsilon)}{\varepsilon}
\geq\int_0^\infty \frac{e^{-\frac t\varepsilon}}{\varepsilon}\liminf_{\tau\to0}\frac{\mathcal W(\bar S_\varepsilon^{\tau}(t))-\mathcal W_\tau(\bar S_\varepsilon^\tau(t))}{\tau}\d t.
$$
Now Proposition \ref{prop:compactness-dissipation-potential} gives, along a subsequence,
$$
\tilde\mu_\varepsilon^\tau\weakstarto\mu_\varepsilon
$$
for some $\mu_\varepsilon\in\young(T_\varepsilon)$, so in particular
$$
\stcurr{\tilde\mu_\varepsilon^\tau}(t)\weakstarto \stcurr{\mu_\varepsilon}(t)\quad \text{ for every }t\in[0,\infty).
$$
The lower semicontinuity in Proposition \ref{prop:compactness-dissipation-potential} proves \eqref{eq:fixed-epsilon-trajectory-action}.
Similarly to the previous proof, by the uniform H\"older continuity of the variation of $\tilde\mu_\varepsilon^\tau$  we have
$$
\mathbf F(\bar S_\varepsilon^\tau(t)-\stcurr{\tilde \mu_\varepsilon^\tau}(t))\leq C\tau ^{1-1/p}\to 0 \text{ as }\tau \to 0,
$$
which altogether shows that 
$$
\bar S_\varepsilon^\tau(t)\weakstarto \stcurr{\mu_\varepsilon}(t)\quad \text{in }\intcur{k}(\overline\Omega),\quad \text{for all } t>0
$$
Using this, we continue our inequality, by  Lemma \ref{lem:yosida-relaxed-slope} we obtain
\begin{align*}
\frac{\mathcal W(T_\varepsilon)-V_\varepsilon(T_\varepsilon)}{\varepsilon}
&\geq\int_0^\infty \frac{e^{-\frac t\varepsilon}}{\varepsilon}\liminf_{\tau\to0}\frac{\mathcal W(\bar S_\varepsilon^{\tau}(t))-\mathcal W_\tau(\bar S_\varepsilon^\tau(t))}{\tau}\d t
\\&\geq
\int_0^\infty \frac{e^{-\frac t\varepsilon}}{\varepsilon}\dualRinf{\stcurr{\mu_\varepsilon}(t)}\d t.
\end{align*}
Setting $\mu_{\varepsilon,T_\varepsilon}:=\mu_\varepsilon$ concludes the proof.
\end{proof}

Now we can proceed with the proof of the main result of this section, Proposition \ref{prop:liminf-Gepsilon}.

\begin{proof}[Proof of Proposition \ref{prop:liminf-Gepsilon}]
For every $\varepsilon>0$, let us take $\mu_\varepsilon\in\young(T_\varepsilon)$ from Lemma \ref{lem:fixed-epsilon-G-lower-bound}.
As $\int_0^\infty\mathcal R_p(\slice{\mu_\varepsilon}{t})\d t\leq\mathcal W(T_\varepsilon)\leq E_\gamma$, Lemma \ref{lem:all-time-slices-Y} gives
$$
\mathbf F(\stcurr{\mu_\varepsilon}(\varepsilon s)- T_\varepsilon)\leq\gamma^{1-1/p}(pE_\gamma)^{1/p}(\varepsilon s)^{1-1/p}\to 0 \text{ as } \varepsilon\to 0 \text{ for all } s>0.
$$
Together with $T_\varepsilon\weakstarto T_0$, we get that 
$$
\stcurr{\mu_\varepsilon}(\varepsilon s)\overset{\ast}{\rightharpoonup} T_0,\quad \text{as }\varepsilon\to 0 \text{ for all }s>0.$$
Hence, by the Fatou lemma and  weak${}^*$ lower semicontinuity of $\dualRinf{\cdot}$ we obtain

\begin{align*}
\liminf_{\varepsilon\to 0} G_\varepsilon(T_\varepsilon)&\geq \liminf_{\varepsilon\to 0} \int_0^\infty \frac{e^{-t/\varepsilon}}{\varepsilon}\dualRinf{\stcurr{\mu_\varepsilon}(t)}\d t\\
&=
 \liminf_{\varepsilon\to 0} \int_0^\infty e^{-s}\dualRinf{\stcurr{\mu_\varepsilon}(\varepsilon s)}\d s\geq \int_0^\infty e^{-s}\dualRinf{T_0}\d s= \dualRinf{T_0}. \qedhere
\end{align*}
\end{proof}

\subsection{The limit \texorpdfstring{$\varepsilon\to 0$}{ε→0}}

We now pass to the limit in the energy balance for the $\varepsilon$-minimizers. In the $\varepsilon$-regime, our minimizers satisfy an exact energy equality (Proposition \ref{prop:eps-energy-balance}). However, as is standard in the calculus of variations, the weak topologies required to extract convergent subsequences only guarantee lower semicontinuity for our functionals, both for  $V_\varepsilon$  (Proposition \ref{prop:value-func-convergence})  and $G_\varepsilon$ (Proposition \ref{prop:liminf-Gepsilon}). So the conclusion is an energy-dissipation inequality. The following theorem gives the global evolution announced in the introduction.

\begin{theorem}[Existence of EDI solutions]\label{thm:energy-dissipation-inequality}
For any initial condition $T_0\in\statespace$, there exists a $\mu\in\young(T_0)$ which is an EDI solution in the sense of Definition \ref{def:EDI-solution}.
\end{theorem}

\begin{proof}
We will pass to the limit in the gradient flow for $V_\varepsilon$, see Proposition \ref{prop:value-funct-gradient-flow}:
$$-\frac{\mathrm d}{\mathrm dt}V_\varepsilon(\stcurr{\mu_\varepsilon}(t))=\mathcal R_p(\slice{\mu_\varepsilon}{t})+\frac {\mathcal W(\stcurr{\mu_\varepsilon}(t))- V_\varepsilon(\stcurr{\mu_\varepsilon}(t))}{\varepsilon}$$
in order to obtain the limit EDI solution after the limit $\varepsilon\to0$.

The bounds of Proposition \ref{prop:energy-estimates} allow us to apply Proposition \ref{prop:compactness-dissipation-potential}. This gives $\mu\in\young(T_0)$ and, for a subsequence of $\varepsilon\to0$,
$$\mu_\varepsilon\weakstarto \mu.$$
Integrating the gradient flow for $\mu_\varepsilon$ over $(0,t)$ we have 
$$\int_0^t \mathcal R_p(\slice{\mu_\varepsilon}{s})+\frac {\mathcal W(\stcurr{\mu_\varepsilon}(s))- V_\varepsilon(\stcurr{\mu_\varepsilon}(s))}{\varepsilon}\d s+ V_\varepsilon(\stcurr{\mu_\varepsilon}(t))= V_\varepsilon(T_0).$$
Let us take $\liminf_{\varepsilon\to0}$. On the right hand side, we have by Proposition~\ref{prop:value-func-convergence} $$\lim_{\varepsilon\to0}V_\varepsilon(T_0)=\mathcal W(T_0).$$
On the left hand side, Proposition \ref{prop:uniform-convergence-state-curves} gives
$$
\stcurr{\mu_\varepsilon}(t)\weakstarto\stcurr{\mu}(t),\quad \text{for all }t\in[0,\infty),
$$
so that we obtain by Proposition~\ref{prop:value-func-convergence} the inequality
$$\liminf_{\varepsilon\to 0} V_\varepsilon(\stcurr{\mu_\varepsilon}(t)) \geq \mathcal W(\stcurr{\mu}(t)),\quad \text{for all }t\in(0,\infty).$$
Further, Proposition \ref{prop:compactness-dissipation-potential} gives lower semicontinuity on every $[0,t]$, since the absolute continuity in time of the limit has now been proved:
$$\liminf_{\varepsilon\to0}\int_0^t\mathcal R_p(\slice{\mu_\varepsilon}{s})\d s\geq \int_0^t\mathcal R_p(\slice{\mu}{s} )\d s\quad \text{for all }t\in(0,\infty).$$
Finally, by the Fatou lemma and Proposition \ref{prop:liminf-Gepsilon} we obtain
\begin{align*}
&\liminf_{\varepsilon\to0}\int_0^t\frac {\mathcal W(\stcurr{\mu_\varepsilon}(s))- V_\varepsilon(\stcurr{\mu_\varepsilon}(s))}{\varepsilon}\d s\\
&\quad\geq \int_0^t\liminf_{\varepsilon\to0}\frac {\mathcal W(\stcurr{\mu_\varepsilon}(s))- V_\varepsilon(\stcurr{\mu_\varepsilon}(s))}{\varepsilon}\d s\geq \int_0^t\dualRinf{\stcurr{\mu}(s)}\d s.
\end{align*}
Altogether, this gives the desired EDI inequality
$$\mathcal W(\stcurr{\mu}(t))+\int_0^t\mathcal R_p(\slice{\mu}{s})\d s+\int_0^t \dualRinf{\stcurr{\mu}(s)}\d s\leq \mathcal W(T_0),$$
finishing the proof.
\end{proof}

\section{Energy-Dissipation balance} \label{sec:balance}

In this section we give conditions under which the energy-dissipation inequality becomes an equality under an assumption of absolute continuity in time.

First, we present an upper-gradient estimate for $\mathcal W$, analogous to Proposition~\ref{prop:chain-rule-epsilon}, with the unrelaxed dual dissipation rate $\dualR{}$ in place of $G_\varepsilon$.

\begin{lemma}[Upper-gradient estimate]\label{lem:fenchel-inequality}
Take any $\mu\in\young$ such that $t\mapsto \mathcal W(\stcurr{\mu}(t))$ is locally absolutely continuous on $[0,\infty)$. Then it holds
$$
\left|\frac{\mathrm d}{\mathrm dt}\mathcal W(\stcurr{\mu}(t))\right|
\leq
\|v\|_{L^p(\mu(t))}
\left(p'\dualR{\stcurr{\mu}(t)}\right)^{1/p'}
\quad\text{for a.e. }t\in(0,\infty).
$$
In particular, we have
$$
\left|\frac{\mathrm d}{\mathrm dt}\mathcal W(\stcurr{\mu}(t))\right|
\leq
\mathcal R_p(\slice{\mu}{t})+\dualR{\stcurr{\mu}(t)}
\quad\text{for a.e. }t\in(0,\infty).
$$
\end{lemma}
\begin{proof}
We proceed analogously as in the proof of Proposition \ref{prop:chain-rule-epsilon}.
Fix a time $t\in (0,\infty)$ where $s\mapsto \mathcal W(\stcurr{\mu}(s))$ is differentiable and which is both left and right Lebesgue point of $s\mapsto\mathcal R_p(\slice{\mu}{s})$. By using the restriction of $\mu$ to $(t,\infty)$, translated back
in time by $t$, as a competitor in $\dualR{\stcurr{\mu}(t)}$, we obtain
$$-\frac{\mathrm d}{\mathrm dt}\mathcal W(\stcurr{\mu}(t))- \mathcal R_p(\slice{\mu}{t})\leq\dualR{\stcurr{\mu}(t)}.$$
For the opposite sign, we use as competitor the reversal $(\mu\mres[0,t])^{\mathrm{rev}}$ (see Lemma \ref{lem:time-reversal})  extended arbitrarily to $(0,\infty)$, to obtain
$$\frac{\mathrm d}{\mathrm dt}\mathcal W(\stcurr{\mu}(t))- \mathcal R_p(\slice{\mu}{t})=-\frac{\mathrm d}{\mathrm dt}\mathcal W(\stcurr{(\mu\mres(0,t))^{\mathrm{rev}}}(t))-\mathcal R_p(\slice{(\mu\mres(0,t))^{\mathrm{rev}}}{0})\leq\dualR{\stcurr{\mu}(t)}.$$
Consequently, we have shown the additive estimate
$$\abs{\frac{\mathrm d}{\mathrm dt}\mathcal W(\stcurr{\mu}(t))}\leq\dualR{\stcurr{\mu}(t)} + \mathcal R_p(\slice{\mu}{t})\quad \text{for a.e. }t\in(0,\infty).$$

To prove the multiplicative estimate, we now rescale time as in the proof of Proposition \ref{prop:chain-rule-epsilon}. We define 
$$
\mu_\lambda:=\timepush{(s\mapsto s/\lambda)}{\mu},
\qquad \lambda>0, 
$$
so that we have (see Lemma \ref{lem:dissipation-pushforward})
\begin{align*}
\stcurr{\mu_\lambda}(s/\lambda)&=\stcurr{\mu}(s) \quad \text{for every }s\in[0,\infty),\\
\mathcal R_p(\slice{\mu_\lambda}{s/\lambda})
&=\lambda^p\mathcal R_p(\slice{\mu}{s}) \quad \text{for a.e. }s\in(0,\infty).
\end{align*}
Thus, the additive estimate applied to $\mu_\lambda\in\young$ for some fixed $\lambda >0$ gives
$$\lambda\abs{\frac{\mathrm d}{\mathrm dt}\mathcal W(\stcurr{\mu}(t))}\leq\dualR{\stcurr{\mu}(t)}+\lambda^p\mathcal R_p(\slice{\mu}{t})\quad\text{ for a.e. }t\in(0,\infty).$$
Taking $\lambda$ from a countable dense set $D\subset(0,\infty)$, we obtain a common
full-measure set of times on which this inequality holds for all $\lambda\in D$. Continuity in $\lambda$ then extends this inequality from $\lambda \in D$ to all $\lambda>0$.

Writing
$\mathcal R_p(\slice{\mu}{t})=\frac1p\|v\|_{L^p(\mu(t))}^p$,
dividing by $\lambda$, and minimizing over $\lambda$, we obtain
\begin{align*}
\abs{\frac{\mathrm d}{\mathrm dt}\mathcal W(\stcurr{\mu}(t))}
&\leq
\inf_{\lambda>0}
\left(
\frac{\lambda^{p-1}}p\|v\|_{L^p(\mu(t))}^p
+\frac{\dualR{\stcurr{\mu}(t)}}{\lambda}
\right)
\\
&=
\|v\|_{L^p(\mu(t))}
\left(p'\dualR{\stcurr{\mu}(t)}\right)^{1/p'},\quad \text{for a.e.\ }t\in(0,\infty),
\end{align*}
which completes the proof.

We remark that the additive estimate follows from the multiplicative one by Young's inequality 
\begin{align*}
\|v\|_{L^p(\mu(t))}
\left(p'\dualR{\stcurr{\mu}(t)}\right)^{1/p'}&\leq
\frac1p\|v\|_{L^p(\mu(t))}^p
+\dualR{\stcurr{\mu}(t)}=
\mathcal R_p(\slice{\mu}{t})
+\dualR{\stcurr{\mu}(t)}. \qedhere
\end{align*}
\end{proof}

Consequently, if the dual dissipation $\dualR{}$ is lower semicontinuous and the energy is absolutely continuous along the solution, the upper-gradient estimate gives the reverse inequality and hence an energy-dissipation balance.

\begin{corollary}[Energy-Dissipation Balance]\label{cor:energy-dissipation-balance}
Let $\mu\in\young(T_0)$ be a solution of the Energy-Dissipation Inequality of Theorem \ref{thm:energy-dissipation-inequality}, and assume that $\dualR{}$ is weak$^*$ lower semicontinuous on $\statespace$. If the map $t \mapsto \mathcal W(\stcurr{\mu}(t))$ is locally absolutely continuous on $[0,\infty)$, then $\mu$ satisfies the Energy-Dissipation Balance for all $t\in(0,\infty)$:
$$
\mathcal W(\stcurr{\mu}(t))+\int_0^t\mathcal R_p(\slice{\mu}{s})\d s+\int_0^t \dualRinf{\stcurr{\mu}(s)}\d s = \mathcal W(T_0).
$$
\end{corollary}
\begin{proof}
The singular part of the derivative of $t \mapsto \mathcal W(\stcurr{\mu}(t))$ is zero, and by assumption $\dualRinf{}=\dualR{}$. The two inequalities in Theorem \ref{thm:energy-dissipation-inequality} and Lemma \ref{lem:fenchel-inequality} give the desired equality.
\end{proof}

The following simple example illustrates the occurrence of a jump down of the energy in the $\varepsilon\to 0$ limit. The example shows two points of opposite signs attracting, and annihilating in finite time. It is noteworthy that the energy decreases by the amount equal to the extra dissipation of the WED $\varepsilon$-minimizers.

\begin{example}[Energy drop in the $\varepsilon\to 0$ limit]\label{ex:energy-drop}
We take $p=2$, $k=0$, $d=1$, $\Omega=(-1,1)$, $\gamma=2$.
We define the energy functional
$$\mathcal W(T):=\mass{T}+|T(\operatorname{id})|,$$
where $\operatorname{id}: \R \to \R$ is treated as a zero form. More explicitly this means
$$
\mathcal W(T)=\sum_{j=1}^m|\theta_j|
+\left|\sum_{j=1}^m\theta_jx_j\right|\quad \text{for }T=\sum_{j=1}^m\theta_j\llbracket x_j\rrbracket\quad \theta_j\in\mathbb Z\setminus\{0\},\quad x_j\in[-1,1] \text{ distinct}.
$$
The energy satisfies
\ref{item:E-bdd-below}--\ref{item:E-bdd-by-mass}. Notably the mass is bounded by $\mathcal{W}(T)$ directly. We fix $x_0\in(0,1)$ and set the initial state as
\[
T_0=\llbracket x_0\rrbracket-\llbracket-x_0\rrbracket\in\statespace.
\]
We will show that for the purposes of minimization of the WED functional $\mathcal J_\varepsilon$, it is sufficient to consider trajectories which are two symmetric points around zero and possibly annihilate in finite time. Precisely, $\mu_\varepsilon \in\young(T_0)$ with slices
\begin{equation}\label{eqn:symmteric-trajectory}
S[\mu_\varepsilon](t)=\begin{cases}
\llbracket x_\varepsilon(t)\rrbracket-\llbracket-x_\varepsilon(t)\rrbracket, &t\in[0,t_\varepsilon)\\
0 & t\in [t_\varepsilon,\infty)
\end{cases}
\end{equation}
for some $t_\varepsilon\in(0,\infty]$ and some $x_\varepsilon\in W^{1,2}_{\loc}([0,\infty))$ with $x_\varepsilon(0)=x_0$, $x_\varepsilon>0$ on $[0,t_\varepsilon)$, and $x_\varepsilon(t_\varepsilon)=0$ if $t_\varepsilon$ is finite.

To see this, take any $\mu\in\young(T_0)$. By the cylinder formula
\eqref{eq:all-time-cylinder-formula} we have 
\[
\stcurr{\mu}(t)-\stcurr{\mu}(s)=\partial\xproj{}_*(\stcurr{\mu}\mres(s,t)),
\quad 0\leq s<t,
\]
so testing by $\omega\equiv1 \in \mathcal D^0([s,t]\times [0,1])$ gives conservation of signed multiplicity: $\langle S[\mu](t),1\rangle=\langle T_0,1\rangle=1-1=0$ for every $t\geq0$. Thus absolute values of multiplicities have to add up to an even number. Since $\mass{S[\mu](t)}\leq\gamma=2$ by Lemma~\ref{lem:all-time-slices-Y}, this means either
\begin{equation}\label{eqn:two-point-slice}
S[\mu](t)=\llbracket a(t)\rrbracket-\llbracket b(t)\rrbracket,
\qquad a(t)\neq b(t)
\end{equation}
or $S(t)=0$. Furthermore, as defining $t_0:=\inf\{t\geq0:S(t)=0\}\in(0,\infty]$ and replacing $S$ with constant zero on $(t_0,\infty)$ if $t_0<\infty$, does not increase the value of $\mathcal J_\varepsilon(\mu)$, any minimizer will start in the first form and may then vanish at some point.

Using also that $\orient$ consists of unit space-time vectors, we get by convexity arguments that $\mu$ can be chosen elementary and thus $\mu=\eym{S}$.

Applying the definition of the dissipation potential and energy we get for $t< t_0$
$$
\mathcal R_2(\mu(t))=\frac{|\dot a(t)|^2+|\dot b(t)|^2}{2} \text{ and } \mathcal{W}(S[\mu](t)) = 2 +\abs{a(t)-b(t)}.
$$
From this it is clear by Jensen's inequality that minimizers must be symmetrical with $a(t) = -b(t)$ and thus satisfy \eqref{eqn:symmteric-trajectory}.

We can thus instead write the WED functional for these trajectories in the form 
\begin{equation*}
\tilde {\mathcal{J}}_\varepsilon(x,t_0)
= \int_0^{t_0} \frac{e^{-t/\varepsilon}}{\varepsilon}\left( 2 + 2x(t) + \varepsilon |\dot{x}(t)|^2 \right)\d t,
\end{equation*}
over pairs $t_0\in(0,\infty]$ and $x\in W^{1,2}_{\loc}([0,\infty))$ satisfying $x(0)=x_0$, $0<x\leq1$ on $[0,t_0)$, and $x(t_0)=0$ when $t_0<\infty$. Denote a minimizing pair by $(x_\varepsilon,t_\varepsilon)$.

Variation in $x$ gives the Euler-Lagrange equation (omitting the common factor 2):
$$\frac{e^{-t/\varepsilon}}{\varepsilon}-\frac{\d}{\d t}\left( e^{-t/\varepsilon}\dot{x}_\varepsilon \right)=0,
\qquad\text{equivalently}\qquad
\frac1\varepsilon+\frac1\varepsilon\dot x_\varepsilon-\ddot x_\varepsilon=0
\quad\text{on }(0,t_\varepsilon).
$$
Its solution with $x_\varepsilon(0)=x_0$ is, with some constant $A_\varepsilon$,
$$x_\varepsilon(t)=x_0+A_\varepsilon\varepsilon\left(e^{t/\varepsilon}-1\right)-t,\qquad 0\leq t<t_\varepsilon.$$
No such solution stays in $(0,x_0]$ for all $t\geq0$, so $t_\varepsilon<\infty$.
Variation in $t_0$ with $x(t_0)=0$ gives
$$0=e^{-t_\varepsilon/\varepsilon}\left(\frac2\varepsilon-|\dot x_\varepsilon(t_\varepsilon^-)|^2\right)
\quad
\text{and thus}
\quad
\dot x_\varepsilon(t_\varepsilon^-)=-\sqrt{\frac2\varepsilon}.$$
Hence $A_\varepsilon$ and $t_\varepsilon$ are given by the system (the second equation is $x_\varepsilon(t_\varepsilon)=0$)
$$
\begin{aligned}
    A_\varepsilon &=e^{-t_\varepsilon/\varepsilon}\left(1-\sqrt{\frac2\varepsilon}\right), \\
    t_\varepsilon &= x_0 + \varepsilon \left(1-\sqrt{\frac2\varepsilon}\right)\left(1-e^{-t_\varepsilon/\varepsilon}\right) \to x_0, \quad \text{as }\varepsilon\to 0.
\end{aligned}
$$

Limit dissipation:
By Proposition~\ref{prop:eps-slope-is-R*DR}, for a.e.\ $t\in(0,t_\varepsilon)$ we have
$$\begin{aligned}
G_\varepsilon(\stcurr{\mu_\varepsilon}(t))
&= \frac{\mathcal W(\stcurr{\mu_\varepsilon}(t))-V_\varepsilon(\stcurr{\mu_\varepsilon}(t))}{\varepsilon}
=\mathcal R_2(\slice{\mu_\varepsilon}{t})=|\dot x_\varepsilon(t)|^2.
\end{aligned}$$
So, the total dissipation of $x_\varepsilon$ is given by
$$ I_\varepsilon:=\int_0^{t_\varepsilon}\left(|\dot x_\varepsilon(t)|^2+G_\varepsilon(\stcurr{\mu_\varepsilon}(t))\right)\d t = \int _0^{t_\varepsilon} 2|\dot x_\varepsilon(t)|^2 \d t.$$
Substituting the velocity gives
\begin{align*}
I_\varepsilon&=\int _0^{t_\varepsilon}2 |\dot x_\varepsilon(t)|^2 \d t
 \\&= \varepsilon \left( 1 - \sqrt{\frac2\varepsilon} \right)^{2} \left( 1 - e^{-2t_\varepsilon/\varepsilon} \right)
- 4\varepsilon \left( 1 - \sqrt{\frac2\varepsilon} \right) \left( 1 - e^{-t_\varepsilon/\varepsilon} \right)
+ 2t_\varepsilon\to 2+2x_0,\quad\text{as } \varepsilon\to0.
\end{align*}

The trajectories converge locally uniformly to a solution of the (formal) evolution equation for such trajectories
$$D\tilde{\mathcal W}(x(t))+D\tilde{\mathcal R}_2(\dot{x}(t))=0,\qquad 0<t<x_0,$$
with $\tilde{\mathcal W}(x)=2+2|x|$ for $x\neq0$, $\tilde{\mathcal W}(0)=0$, and $\tilde{\mathcal R}_2(\dot{x})=|\dot{x}|^2$, which under the same assumptions gives $$x(t)=\begin{cases}x_0-t, &0\leq t<x_0 \\ 0, & x_0\leq t <\infty. \end{cases}$$

Notably, however, the total dissipation in that case is
$$\int_0^{x_0} \dot{x}(t) D\tilde{\mathcal R}_2(\dot{x}(t))\d t = \int_0^{x_0} 2\d t = 2x_0.$$

Thus the limiting dissipation of the WED minimizers exceeds the dissipation of the limit trajectory by $2$. This is precisely the energy jump from $2$ to $0$ when the two particles annihilate at $t=x_0$.
\end{example}

\section{Appendix: Calculations for the examples}\label{appendix}

Here we provide the full calculations for the examples of Sections \ref{ssec:examples-currents} and \ref{ssec:examples-revisited}.
For the calculations below, we recall the orientation formula for the pushforward of a $2$-current in a plane.
Let $Q\subset \R^2$ be a rectangle and let $F\in C^{1}(\overline Q;\R^{1+d})$ be injective with $\operatorname{rank}\nabla F=2$. Then 
    $F_*\llbracket Q \rrbracket\in\intcur{2}(\R^{1+d})$ with the orienting $2$-vector
\begin{equation}\label{eqn:rectangle-pushforward}
    \overrightarrow{F_*\llbracket Q\rrbracket}(F(t,x)) = \frac{\partial_0 F(t,x)\wedge \partial _1 F(t,x)}{|\partial_0 F(t,x)\wedge \partial _1 F(t,x)|},\quad \text{a.e. }(t,x)\in Q.
\end{equation}
This follows from the definition of current pushforward and the area formula.

\subsection{Calculations for Example~\ref{ex:zigzag-current}: Zig-zag}
\label{sec:zigzag-appendix}

Recall the full current and its decomposition $$S_{\Delta t,\alpha,N}=S_{\Delta t,\alpha,N}^{\mathrm m}+S_{\Delta t,\alpha,N}^{\mathrm s}$$ from Example \ref{ex:zigzag-current}, with zigzag height $h_{\alpha,N}=\cot\alpha/(2N)$. We write \begin{equation}\label{eq:zigzag-slice-graph-appendix}
y_t^{\Delta t,\alpha,N}(x):=\frac{t}{\Delta t}+h_{\alpha,N}\zigzag_{1/N}(x), \quad t\in[0,\Delta  t] \quad \text{with }h_{\alpha,N}:=\frac{\cot\alpha}{2N},
\end{equation}
and for fixed $\Delta t,\alpha,N$, we abbreviate $y_t=y_t^{\Delta t,\alpha,N}$ and set $y_t=y_{\Delta t}$ for $t\in[\Delta t,1]$. For $t\in[0,1]$ and almost every $x\in(0,1)$, we define the sign
$$\sigma_N(x):=\operatorname{sgn}\zigzag'_{1/N}(x)\in\{-1,+1\},\qquad \text{so that}\quad y'_t(x)=\sigma_N(x)\cot\alpha,$$
and in the following calculations we abbreviate $\sigma=\sigma_N(x)$.

The slice $S_{\Delta t, \alpha, N}(t)$ is the graph of $x\mapsto y_t(x)$ in \eqref{eq:zigzag-slice-graph}, so we have
$$
\vec S_{\Delta t, \alpha, N}(t)(x,y_t(x))= \frac{e_1+ y_t'e_2}{\sqrt{1+|y'_t|^2}}= \sin\alpha\,e_1+\sigma\cos\alpha\,e_2=: \tau_\sigma,
$$
the latter expression verifies that this in fact corresponds to Figure \ref{fig:zigzag-slice}.

The normal velocity vector is obtained by projecting $(1/\Delta t)e_2$ to the direction normal to $\tau_\sigma$, namely
$$
\frac{D}{Dt}S_{\Delta t,\alpha,N}=\frac1{\Delta t}\left(e_2-\sigma\cos\alpha\,\tau_\sigma\right)
=\frac{-\sigma\sin\alpha\cos\alpha\,e_1+\sin^2\alpha\,e_2}{\Delta t}=:v_\sigma.
$$

Applying the pushforward formula \eqref{eqn:rectangle-pushforward} with
\[
F(t,x)=(t,x,y_t(x)),
\]
so that, for $t\in(0,\Delta t)$,
\[
\partial_0F(t,x)=e_0+e_2/\Delta t,
\qquad
\partial_1F(t,x)=e_1+\sigma\cot\alpha\,e_2,
\]
gives
$$
\partial _0 F \wedge \partial_1F = e_0\wedge e_1+ \sigma \cot \alpha\, e_0\wedge e_2 + \frac{1}{\Delta t}e_2\wedge e_1,\quad \text{and}\quad |\partial _0 F \wedge \partial_1F|= \sqrt{1+\sigma^2\cot^2 \alpha + 1/(\Delta t)^2}.
$$
Normalizing this $2$-vector and multiplying numerator and denominator by $\sin\alpha$ gives the space-time orienting $2$-vector
\begin{align*}
\vec S_{\Delta t,\alpha,N}\left(t,x,y_t(x)\right)
&=\frac{\partial _0 F \wedge \partial_1F}{\abs{\partial _0 F \wedge \partial_1F}}\\
&=\frac{\sin\alpha\,e_0\wedge e_1+\sigma\cos\alpha\,e_0\wedge e_2+(\sin\alpha/\Delta t)e_2\wedge e_1}{\sqrt{1+(\sin\alpha/\Delta t)^2}}
=:\eta_{\Delta t,\alpha}^{\sigma}.
\end{align*}
One can readily verify that
$$
\eta_{\Delta t,\alpha}^\sigma
=\frac{(e_0+v_\sigma)\wedge\tau_\sigma}{\sqrt{1+|v_\sigma|^2}}
=|\eta_{\Delta t,\alpha}^\sigma \mres \mathrm dt|(e_0+v_\sigma)\wedge\tau_\sigma
$$
which shows that, in the sense of Lemma \ref{lem:space-time-decomposition},
$$
v_\sigma = v(\eta_{\Delta t,\alpha}^\sigma)\quad \text{and } \tau_\sigma=\tau(\eta_{\Delta t,\alpha}^\sigma).
$$
Thus the slice orientation and normal velocity obtained from the geometric description agree with those recovered from the space-time $2$-vector computed by the pushforward formula.

The normal speed then satisfies
$$
\abs{\frac{D}{Dt}S_{\Delta t,\alpha,N}}=|v_\sigma|=\frac{\sin\alpha}{\Delta t}.
$$

Thus, the moving slices and their dissipation potential satisfy
\begin{align*}
\mass{S_{\Delta t,\alpha,N}(t)}&=\frac1{\sin\alpha},
&|v_\sigma|&=\frac{\sin\alpha}{\Delta t},\\
\mathcal R_p(\slice{\mu_{\Delta t,\alpha,N}^{\mathrm m}}{t})
&=\frac{\sin^{p-1}\alpha}{p(\Delta t)^p},
&&0<t<\Delta t.
\end{align*}
The stationary tail $S_{\Delta t,\alpha,N}^{\mathrm s}(t)$ has zero velocity and therefore does not contribute to $\Psi_p$. Consequently,
\begin{align*}
\Psi_p(S_{\Delta t,\alpha,N};[0,1])
&= \frac1p\int_0^{\Delta t}\int|v_\sigma|^p\d\|S_{\Delta t,\alpha,N}(t)\| \d t\\
&=\frac1p\int_0^{\Delta t}\left(\frac{\sin\alpha}{\Delta t}\right)^p\frac1{\sin\alpha}\d t
=\frac1p\left(\frac{\sin\alpha}{\Delta t}\right)^{p-1}
\end{align*}
which proves \eqref{eqn:zigzag-finite-cost}.

\textbf{Current convergence.}
In the convergence arguments, we write $\sigma_N=\sigma_N(x)$ and $\sigma_{N_i}=\sigma_{N_i}(x)$, suppressing the $x$ argument.
Take a test form $\omega\in\mathcal D^2([0,1]\times\overline\Omega)$.
Using the time-disintegration formula in Proposition \ref{prop:time-disintegr-NC-current}, we have
\begin{align*}
\langle S_{\Delta t,\alpha,N}^{\mathrm m},\omega\rangle 
&= \int_0^{\Delta t} \int \langle (e_0+v_{\sigma_N})\wedge \tau _{\sigma_N},\omega\rangle \d\|S_{\Delta t,\alpha,N}^{\mathrm m}(t)\|\d t
\\
&=\int_0^{\Delta t} \int \left\langle (e_0+v_{\sigma_N})\wedge \tau _{\sigma_N},\omega\right\rangle \d\|S_{\Delta t,\alpha,N}^{\mathrm m}(t)\|\d t.
\end{align*}
Substituting the expressions for $v_{\sigma_N}$ and $\tau_{\sigma_N}$ and simplifying the wedge product yields
\begin{align*}
\langle S_{\Delta t,\alpha,N}^{\mathrm m},\omega\rangle 
&=
\int_0^{\Delta t} \int \sin\alpha\left\langle (e_0+e_2/\Delta t)\wedge e_1+ \sigma_N\Delta t\cot \alpha / \Delta t\, e_0\wedge e_2 ,\omega\right\rangle \d\|S_{\Delta t,\alpha,N}^{\mathrm m}(t)\|\d t.
\end{align*}
Parametrizing the slice $S_{\Delta t,\alpha,N}^{\mathrm m}(t)$ by $x\in(0,1)$ gives the area factor $1/(\sin\alpha)$
\begin{align*}
\langle S_{\Delta t,\alpha,N}^{\mathrm m},\omega\rangle 
&=
\int_0^{\Delta t} \int_0^1 \left\langle (e_0+e_2/\Delta t)\wedge e_1 +\sigma_N\Delta t\cot \alpha / \Delta t e_0\wedge e_2 ,\omega(t,x,y_t(x))\right\rangle \d x\d t.
\end{align*}
Changing the variables\footnote{This change of variables is necessary only for the case $\Delta t_i\to 0$, where the limit moving part becomes purely spatial and can no longer be parametrized by $t$.} $y=t/\Delta t$ gives
\begin{align*}
&\langle S_{\Delta t,\alpha,N}^{\mathrm m},\omega\rangle\\
&\quad=\int_0^{1} \int_0^1 \left\langle (\Delta t e_0+e_2)\wedge e_1 +\sigma_N\Delta t\cot \alpha \, e_0\wedge e_2 ,\omega(\Delta t \, y,x,y+h_{\alpha,N}\zigzag_{1/N}(x))\right\rangle \d x\d y.
\end{align*}
Now take a sequence $\Delta t_i,\alpha_i,N_i$ satisfying \eqref{eqn:limits-zigzag}, recall that this reads
\begin{equation*}\tag{\ref*{eqn:limits-zigzag}}
\lim_{i\to\infty} N_i=\infty, \quad \lim_{i\to\infty} \alpha_i=\alpha\in[0,\pi/2),\quad \lim_{i\to\infty} \Delta t_i =\Delta t\in[0,1],\quad \lim_{i\to\infty} h_{\alpha_i,N_i}=0.
\end{equation*}
We have 
$$\sigma_{N_i}\weakstarto0\quad \text{} \quad \text{in }L^\infty((0,1)),\quad h_{\alpha_i,N_i}\to 0, \quad |\zigzag_{1/N_i}|\leq 1,\quad  \Delta t_i \cot \alpha_i \to \beta\cos\alpha$$
and thus we can pass to the limit to obtain
\begin{align*}
&\langle S_{\Delta t_i,\alpha_i,N_i}^{\mathrm m},\omega\rangle\\
&\quad=\int_0^{1} \int_0^1 \left\langle (\Delta t_i e_0+e_2)\wedge e_1 +\sigma_{N_i}\Delta t_i\cot \alpha_i \, e_0\wedge e_2 ,\omega(\Delta t_i \, y,x,y+h_{\alpha_i,N_i}\zigzag_{1/N_i}(x))\right\rangle \d x\d y
\\
&\quad\to
\int_0^{1} \int_0^1 \left\langle (\Delta t e_0+e_2)\wedge e_1 ,\omega(\Delta t\,y,x,y)\right\rangle \d x\d y=:\langle S_{\Delta t}^{\mathrm m},\omega\rangle.
\end{align*}
For $\Delta t=0$, this formula gives
$$
S_{\Delta t}^{\mathrm m}= -\llbracket\{0\}\times (0,1)\times(0,1)\rrbracket,
$$
where the minus is due to the limit having orientation $e_2\wedge e_1=-e_1\wedge e_2$.

For $\Delta t>0$, we change the variables back to $t=\Delta t \, y$ to obtain 
$$
\langle S_{\Delta t}^{\mathrm m},\omega\rangle =\int_0^{\Delta t} \int_0^1 \left\langle (e_0+e_2/\Delta t)\wedge e_1 ,\omega(t,x,t/\Delta t)\right\rangle \d x\d t.
$$
By the time-disintegration formula in Proposition \ref{prop:time-disintegr-NC-current}, the slices have
$$
S_{\Delta t}^{\mathrm m}(t) = \llbracket (0,1)\times \{t/\Delta t\} \rrbracket, \quad  \frac{D}{Dt}S_{\Delta t}^{\mathrm m} = e_2/\Delta t ,\quad t\in(0,\Delta t). 
$$
Alternatively, the pushforward formula \eqref{eqn:rectangle-pushforward} gives the space-time representation of the same moving current directly:
$$S_{\Delta t}^{\mathrm m}=
\left((t,x)\mapsto(t,x,t/\Delta t)\right)_*
\llbracket(0,\Delta t)\times(0,1)\rrbracket.
$$
For the stationary parts, we have
\begin{align*}
\lim_{i\to\infty}\langle S_{\Delta t_i,\alpha_i,N_i}^{\mathrm s},\omega\rangle
&=\int_{\Delta t}^1\int_0^1
\langle e_0\wedge e_1,\omega(t,x,1)\rangle\d x\d t=\langle S_{\Delta t}^{\mathrm s},\omega\rangle.
\end{align*}
Thus we have 
$$
S_{\Delta t_i,\alpha_i,N_i}^{\mathrm s} \weakstarto S_{\Delta t}^{\mathrm s} = \llbracket (\Delta t,1) \times (0,1)\times\{1\}\rrbracket.
$$
Combining the two parts gives the convergence
$$
 S_{\Delta t_i,\alpha_i,N_i} \weakstarto S_{\Delta t} = S_{\Delta t}^{\mathrm m}+S_{\Delta t}^{\mathrm s}.
$$

For $\Delta t>0$, the straight moving line has unit length and normal speed $1/\Delta t$, while its tail is stationary. Hence
$$
\lim_{i\to\infty}\Psi_p(S_{\Delta t_i,\alpha_i,N_i};[0,1])
=\frac1p\left(\frac{\sin\alpha}{\Delta t}\right)^{p-1}
<\frac1p\left(\frac1{\Delta t}\right)^{p-1}
=\Psi_p(S_{\Delta t};[0,1]),
$$
where the inequality is strict for all $0\leq\alpha<\pi/2$. At the limiting angle $\alpha=0$, the left-hand side is zero and the right-hand side is positive. The restriction $\alpha>0$ is only needed for the finite varicurrent limit below.

For $\Delta t=0$, the limit is instantaneous. In addition to \eqref{eqn:limits-zigzag}, recall the scaling condition \eqref{eqn:beta-zigzag}:
\begin{equation*}\tag{\ref*{eqn:beta-zigzag}}
\lim_{i\to\infty}\frac{\Delta t_i}{\sin\alpha_i}
=\beta\in(0,\infty).
\end{equation*}
The spatial current $S_0^{\mathrm m}$ has positive mass and fails \eqref{eqn:NC-current}, so we use the extended definition of $\Psi_p$, which gives $\Psi_p(S_0;[0,1])=\infty$. Under this scaling,
$$
\lim_{i\to\infty}\Psi_p(S_{\Delta t_i,\alpha_i,N_i};[0,1])
=\frac1p\beta^{1-p}<\infty=\Psi_p(S_0;[0,1]).
$$

\textbf{Varicurrent convergence.}

Take a test function $g\in C_c([0,1]\times\overline\Omega\times\orient)$. For $\Delta t>0$ and $\alpha>0$, put
$$
\beta_{\Delta t,\alpha}:=\frac{\Delta t}{\sin\alpha},
\qquad
\eta_{\Delta t,\alpha}^{\pm}
:=
\frac{
(\Delta t e_0+e_2)\wedge e_1
\pm\Delta t\cot\alpha\,e_0\wedge e_2
}{
\sqrt{1+\beta_{\Delta t,\alpha}^2}
}.
$$
These are the two orientations of the moving part, namely $\vec S_{\Delta t, \alpha,N}(t,x,y_t(x))=\eta_{\Delta t,\alpha}^{\sigma_N(x)}$.

We proceed analogously as in the calculation for currents. The difference here is that the two orientations do not average, but are both retained in the limit.

Using the moving orientations $\vec S_{\Delta t,\alpha,N}^{\mathrm m}=\eta_{\Delta t,\alpha}^{\sigma_N}$ and the time-disintegration formula in Proposition~\ref{prop:time-disintegr-NC-current}, with $1/|\eta_{\Delta t,\alpha}^{\sigma_N}\mres \mathrm dt|=\sqrt{1+(\sin\alpha/\Delta t)^2}$, gives
\begin{align*}
\left\langle  g, \eym{S_{\Delta t,\alpha,N}^{\mathrm m}}\right\rangle
&=
\int g(t,x,y,\eta_{\Delta t,\alpha}^{\sigma_N})
\,\d\|S_{\Delta t,\alpha,N}^{\mathrm m}\|(t,x,y)
\\
&=
\int_0^{\Delta t}\int_{\overline\Omega}
\,
g(t,x,y,\eta_{\Delta t,\alpha}^{\sigma_N})
\sqrt{1+(\sin\alpha/\Delta t)^2}\d\|S_{\Delta t,\alpha,N}^{\mathrm m}(t)\|(x,y)\,\d t.
\end{align*}
Parametrizing the slice $S_{\Delta t,\alpha,N}^{\mathrm m}(t)$ by $x\in(0,1)$ gives the area factor $1/\sin \alpha$:
\begin{align*}
\left\langle  g, \eym{S_{\Delta t,\alpha,N}^{\mathrm m}}\right\rangle
&=
\sqrt{1/\sin^2\alpha+1/(\Delta t)^2}
\int_0^{\Delta t}\int_0^1
g\left(
t,x,y_t(x),\eta_{\Delta t,\alpha}^{\sigma_N}
\right)\d x\,\d t.
\end{align*}
Changing variables $y=t/\Delta t$ gives
\begin{align*}
\left\langle  g, \eym{S_{\Delta t,\alpha,N}^{\mathrm m}}\right\rangle
&=
\sqrt{(\Delta t/\sin\alpha)^2+1}
\int_0^1\int_0^1
g\left(
\Delta t\,y,x,
y+h_{\alpha,N}\zigzag_{1/N}(x),
\eta_{\Delta t,\alpha}^{\sigma_N}
\right)\d x\,\d y.
\end{align*}

Take a sequence satisfying \eqref{eqn:limits-zigzag} and the ratio limit in \eqref{eqn:beta-zigzag}, namely
\begin{equation*}\tag{\ref*{eqn:beta-zigzag}}
\lim_{i\to\infty}\frac{\Delta t_i}{\sin\alpha_i}
=\beta\in(0,\infty).
\end{equation*}
For $\Delta t>0$, this means $\alpha>0$ and $\beta=\Delta t/\sin\alpha$.
For $\Delta t=0$, this is an additional scaling assumption.

In both cases, $\Delta t=\beta\sin\alpha$. We use this relation throughout the following calculation.

For each of the two fixed signs,
$$
\eta_{\Delta t_i,\alpha_i}^{\pm}\longrightarrow
\eta_{\alpha,\beta}^{\pm}
:=
\frac{
(\Delta t e_0+e_2)\wedge e_1
\pm\beta\cos\alpha\,e_0\wedge e_2
}{
\sqrt{1+\beta^2}
},\quad \text{where }\Delta t =\beta\sin\alpha.
$$

As $\sigma_{N_i}$ oscillates evenly over $x\in(0,1)$ on intervals of length $1/N_i$, the term $h_{\alpha_i,N_i}\zigzag_{1/N_i}$ converges uniformly to $0$ and $g$ is uniformly continuous, we obtain
\begin{align*}
\left\langle  g, \eym{S_{\Delta t_i,\alpha_i,N_i}^{\mathrm m}}\right\rangle
&\to
\frac{\sqrt{\beta^2+1}}{2}
\int_0^1\int_0^1
g\left(
\Delta t\,y,x,
y,
\eta_{\alpha,\beta}^+
\right)+g\left(
\Delta t\,y,x,
y,
\eta_{\alpha,\beta}^-
\right)\d x\,\d y\\
&=:\langle g,\mu_{\alpha,\beta}^{\mathrm m}\rangle.
\end{align*}
For $\Delta t = 0$, we have $\alpha=0$ and 
$$
\eta_{0,\beta}^{\pm}
=\frac{
e_2 \wedge e_1
\pm\beta\,e_0\wedge e_2
}{
\sqrt{1+\beta^2}
}.
$$
Thus the limiting measure is
$$
\mu_{0,\beta}^{\mathrm m}=\delta_0(\mathrm dt)\otimes \mathcal H^2((0,1)\times(0,1))(\mathrm dx,\mathrm dy) \otimes \sqrt{1+\beta^2}\frac{\delta_{\eta_{0,\beta}^+}+\delta_{\eta_{0,\beta}^-}}{2}(\mathrm d\eta).
$$
For $\Delta t=\beta\sin\alpha>0$ we change the variables back to $t=\Delta t\,y$ to obtain
$$
\langle g,\mu_{\alpha,\beta}^{\mathrm m}\rangle=\frac{\sqrt{1+\beta^2}}{2\Delta t}
\int_0^{\Delta t}\int_0^1
g\left(
t,x,t/\Delta t,\eta_{\alpha,\beta}^+
\right)+g\left(
t,x,t/\Delta t,\eta_{\alpha,\beta}^-
\right)\d x\,\d t,
$$
and we read
$$
\mu_{\alpha,\beta}^{\mathrm m}
=\mathcal H^2((0,\Delta t)\times (0,1)\times\{t/\Delta t\})(\mathrm dt,\mathrm dx,\mathrm dy)
\otimes\frac{\sqrt{1+\beta^2}}{2\Delta t}
(\delta_{\eta_{\alpha,\beta}^+}+\delta_{\eta_{\alpha,\beta}^-})(\mathrm d\eta).
$$

For $\Delta t>0$, convergence of the stationary part follows by the same argument:
\begin{align*}
\left\langle  g, \eym{S_{\Delta t_i,\alpha_i,N_i}^{\mathrm s}}\right\rangle
&=
\frac1{\sin\alpha_i}
\int_{\Delta t_i}^1\int_0^1
g\left(
t,x,y_{\Delta t}^{\Delta t_i,\alpha_i,N_i}(x),
\eta_{\alpha_i}^{\sigma_{N_i}}
\right)\d x\,\d t
\\
&\to
\frac1{2\sin\alpha}
\int_{\Delta t}^1\int_0^1
g(t,x,1,\eta_{\alpha,\beta}^+)+g(t,x,1,\zeta_{\alpha,\beta}^-)\d x\,\d t.
\end{align*}
Thus the stationary parts converge
\begin{align*}
\eym{S_{\Delta t_i,\alpha_i,N_i}^{\mathrm s}}
\weakstarto\mu_{\Delta t,\alpha}^{\mathrm s}
:=\mathcal H^2((\Delta t,1)\times(0,1))(\mathrm dt,\mathrm dx)\otimes \frac{1}{\sin\alpha}\frac{\delta_{\eta_{\alpha,\beta}^+}+\delta_{\eta_{\alpha,\beta}^-}}{2}(\mathrm d\eta)
\end{align*}
Combining the moving and stationary parts gives
$$
\eym{S_{\Delta t_i,\alpha_i,N_i}}
\weakstarto
\mu_{\alpha,\beta}
:=
\mu_{\alpha,\beta}^{\mathrm m}
+
\mu_{\alpha,\beta}^{\mathrm s}.
$$
For $\Delta t=0$ we have $\alpha_i\to0$ and the stationary masses $\mass{S_{\Delta t_i,\alpha_i,N_i}^{\mathrm s}}=(1-\Delta t_i)/\sin\alpha_i$ diverge. There is no finite varicurrent limit on the stationary time region. This mass could be made finite by straightening $S_{\Delta t_i,\alpha_i,N_i}$ shortly after time $\Delta t_i$; we will refrain from giving the formal construction as this is not the point of this demonstration. Thus only the moving limit $\mu_{0,\beta}^{\mathrm m}=\mu_\beta^{\mathrm m}$ is asserted in case $\Delta t=0$, the currents still converge by the preceding current calculation.

The following identities give the extended costs in the two regimes. For $\Delta t>0$ and $\alpha>0$,
$$
\int_0^1\mathcal R_p(\slice{\mu_{\alpha,\beta}}{t})\d t
=
\frac1p\left(\frac{\sin\alpha}{\Delta t}\right)^{p-1}
=
\lim_{i\to\infty}
\Psi_p(S_{\Delta t_i,\alpha_i,N_i};[0,1]).
$$

For $\Delta t=0$ the extended dissipation potential is given by
$$
\Psi_p(\mu_\beta^{\mathrm m};[0,1])
=
\frac1p\beta^{1-p}
=
\lim_{i\to\infty}
\Psi_p(S_{\Delta t_i,\alpha_i,N_i}^{\mathrm m};[0,1]).
$$

\subsection{Calculations for Example~\ref{ex:thousand-cuts-current}: Thousand Cuts}
\label{sec:thousand-cuts-appendix}

Recall the space-time current $S_{\Delta t,N}\in\intcur{2}([0,1]\times\overline\Omega)$ from Example \ref{ex:thousand-cuts-current} and its decomposition
$$
S_{\Delta t,N}=S_{\Delta t,N}^{\mathrm m}+S_{\Delta t,N}^{\mathrm s}.
$$
Recall that the $y$-coordinate of the $j$-th moving line at time $t$ is given by
$$
y_{j,\Delta t,N}^\pm(t):=\frac{2j\pm t/\Delta t}{2N+1},\qquad 0\leq t\leq\Delta t,
$$
where we abbreviate $y_{j}^\pm:=y_{j,\Delta t,N}^\pm$ for fixed $\Delta t, N$, with $j=0,\ldots,N$ for $y_j^+$, and $j=1,\ldots,N$ for $y_j^-$.
Here $0<\Delta t\leq1$, and $S_{\Delta t,N}^{\mathrm s}$ is its stationary part with orientation $e_0\wedge e_1$. For $0<t<\Delta t$ and $0<x<1$, the normal velocity and spatial slice orientation on the moving lines are
\begin{align*}
v\bigl(\vec S_{\Delta t,N}(t,x,y_j^\pm(t))\bigr)
&=\pm\frac{e_2}{(2N+1)\Delta t},\\
\vec S_{\Delta t,N}|_t(t,x,y_j^\pm(t))
&=\pm e_1.
\end{align*}
The space-time orientations of $S_{\Delta t, N}$ on the moving parts are
\begin{align*}
\vec S_{\Delta t,N}(t,x,y_j^+(t))=\eta_N^+
&:=\frac{((2N+1)\Delta t\,e_0+e_2)\wedge e_1}{\ell_{\Delta t,N}}, \quad j=0,\dots,N\\
\vec S_{\Delta t,N}(t,x,y_j^-(t))=\eta_N^-
&:=\frac{((2N+1)\Delta t\,e_0-e_2)\wedge(-e_1)}{\ell_{\Delta t,N}}\quad j=1,\dots,N.
\end{align*}
Here $t\in(0,\Delta t)$ and $x\in(0,1)$. The quantity $\ell_{\Delta t,N}$ is defined by \eqref{eqn:ell-cuts}:
\begin{equation*}\tag{\ref*{eqn:ell-cuts}}
\ell_{\Delta t,N}:=\sqrt{1+(2N+1)^2\Delta t^2},
\qquad
\ell_\beta:=\sqrt{1+\beta^2}.
\end{equation*}

For $0<t<\Delta t$, it follows that
$$
\mass{S_{\Delta t,N}(t)}=2N+1,\qquad |v|=\frac1{(2N+1)\Delta t},\qquad
\mathcal R_p(\slice{\eym{S_{\Delta t,N}}}{t})=\frac{2N+1}{p((2N+1)\Delta t)^p}.
$$
The tail has zero velocity, so
\[
\Psi_p(S_{\Delta t,N};[0,1])
=\Psi_p(S_{\Delta t,N}^{\mathrm m};[0,1]) =\int_0^{\Delta t}\mathcal R_p(\slice{\eym{S_{\Delta t,N}}}{t})\d t
=\frac1p\bigl((2N+1)\Delta t\bigr)^{1-p}.
\]

For $\omega\in\mathcal D^2([0,1]\times\overline\Omega)$, the time-disintegration formula of Proposition \ref{prop:time-disintegr-NC-current} gives
\begin{align*}
\langle S_{\Delta t,N}^{\mathrm m},\omega\rangle 
=&\int_{0}^{\Delta t}\int_{\overline\Omega}\left\langle (e_0+v(\vec S_{\Delta t, N}(t,\cdot)))\wedge \vec S_{\Delta t,N}|_t, \omega (t,\cdot)\right\rangle \d\|S_{\Delta t,N}(t)\|\d t
\\
=&\int_{0}^{\Delta t} \int_0^1   \left\langle \left(e_0+\frac{e_2}{(2N+1)\Delta t}\right)
\wedge 
e_1,\sum_{j=0}^N\omega\left(t,x, y_j^+(t) \right)  \right\rangle \d x \d t
 \\
 &+\int_{0}^{\Delta t} \int_0^1   \left\langle \left(e_0-\frac{e_2}{(2N+1)\Delta t}\right)
\wedge 
(-e_1),\sum_{j=1}^N\omega\left(t,x, y_j^-(t) \right)  \right\rangle \d x \d t.
\end{align*}
We change to the vertical coordinate by $y=y_j^\pm(t)$ separately on the upward and downward lines:
\begin{align*}
t=(y_j^+)^{-1}(y)&=\Delta t\bigl((2N+1)y-2j\bigr),
&&\frac{2j}{2N+1}<y<\frac{2j+1}{2N+1},\quad j=0,\ldots,N,\\
t=(y_j^-)^{-1}(y)&=\Delta t\bigl(2j-(2N+1)y\bigr),
&&\frac{2j-1}{2N+1}<y<\frac{2j}{2N+1},\quad j=1,\ldots,N.
\end{align*}
On a downward line, $\d t=-(2N+1)\Delta t\d y$ and reversing the integration limits gives the positive factor $(2N+1)\Delta t$, while retaining the negative slice orientation. Consequently,
\begin{align*}
&\langle S_{\Delta t,N}^{\mathrm m},\omega\rangle\\
&\quad=\sum_{j=0}^N\int_{\frac{2j}{2N+1}}^{\frac{2j+1}{2N+1}}\int_0^1
\left\langle((2N+1)\Delta t\,e_0+e_2)\wedge e_1,
\omega\bigl((y_j^+)^{-1}(y),x,y\bigr)\right\rangle\d x\d y\\
&\qquad+\sum_{j=1}^N\int_{\frac{2j-1}{2N+1}}^{\frac{2j}{2N+1}}\int_0^1
\left\langle(-(2N+1)\Delta t\,e_0+e_2)\wedge e_1,
\omega\bigl((y_j^-)^{-1}(y),x,y\bigr)\right\rangle\d x\d y.
\end{align*}
Parametrizing this curve by height $y$, the inverse functions give its time coordinate on successive upward and downward branches, namely
$$
\Delta t\zigzag_{\frac{2}{2N+1}}(y)=
\begin{cases}
(y_j^+)^{-1}(y),
&\displaystyle\frac{2j}{2N+1}\leq y\leq\frac{2j+1}{2N+1},\quad j=0,\ldots,N,\\
(y_j^-)^{-1}(y),
&\displaystyle\frac{2j-1}{2N+1}\leq y\leq\frac{2j}{2N+1},\quad j=1,\ldots,N.
\end{cases}
$$
Define the oscillating sign by
\[
\sigma_N(y)
:=\operatorname{sgn}\zigzag_{\frac{2}{2N+1}}'(y)
\in\{-1,+1\}
\qquad\text{for a.e. }y\in(0,1).
\]
Thus $\sigma_N(y)=+1$ on the upward lines and $\sigma_N(y)=-1$ on the downward lines.
We subsequently abbreviate $\sigma_N=\sigma_N(y)$ for readability.
Then
$$
\Delta t\zigzag_{\frac{2}{2N+1}}'(y)
=\sigma_N(y)(2N+1)\Delta t
\quad\text{for a.e. }y\in(0,1).
$$
Combining the integrals over these intervals therefore gives the formula
\begin{align*}
&\langle S_{\Delta t,N}^{\mathrm m},\omega\rangle\\
&\quad=\int_0^1\int_0^1
\left\langle
\bigl((2N+1)\Delta t\,\sigma_N e_0+e_2\bigr)
\wedge e_1,
\omega\left(\Delta t\zigzag_{\frac{2}{2N+1}}(y),x,y\right)\right\rangle\d x\d y.
\end{align*}
This explains that the graph in the $(t,y)$-plane is a zig-zag, as illustrated in Figure \ref{fig:thousand-cuts-space-time}.

\textbf{Current convergence.}

Consider a sequence $\Delta t_i, N_i$ as in \eqref{eqn:limits-cuts}, namely
\begin{equation*}\tag{\ref*{eqn:limits-cuts}}
 \lim_{i\to\infty}N_i= \infty,
 \quad
 \lim_{i\to\infty}\Delta t_i=0,
 \quad
 \lim_{i\to\infty}(2N_i+1)\Delta t_i=\beta\in[0,\infty).
\end{equation*}

As $\sigma_{N_i}$ oscillates evenly between $\pm 1$ and $(2N_i+1)\Delta t_i$ is bounded, we see that
$$
\sigma_{N_i}(2N_i+1)\Delta t_i \weakstarto 0 \quad\text{ in }L^\infty((0,1)).
$$
Since $0\leq \zigzag_{\frac{2}{2N+1}}\leq1$, the functions $\Delta t_i \zigzag_{\frac{2}{2N_i+1}}$ converge uniformly to zero. Therefore we have 
$$
\langle S_{\Delta t_i,N_i}^{\mathrm m},\omega\rangle
\to
\int_0^1\int_0^1
\langle e_2\wedge e_1,\omega(0,x,y)\rangle\d x\d y
=\langle S_0^{\mathrm m},\omega\rangle.
$$
For the stationary tail we have separately
$$
\langle S_{\Delta t_i,N_i}^{\mathrm s},\omega\rangle
=\int_{\Delta t_i}^1\int_0^1
\langle e_0\wedge e_1,\omega(t,x,1)\rangle\d x\d t
\to\langle S_0^{\mathrm s},\omega\rangle.
$$
Combining the moving part with the stationary tail, we have shown the convergence of the full currents
$$S_{\Delta t_i, N_i}\weakstarto S_0\quad \text{in }\intcur{2}([0,1]\times\overline\Omega).$$ The current $S_0$ fails \eqref{eqn:NC-current} because of the purely spatial part $S_0^{\mathrm m}$ at time zero, and therefore
$$
\Psi_p(S_0;[0,1])=\infty.
$$
For $\beta>0$,
$$
\lim_{i\to\infty}\Psi_p(S_{\Delta t_i,N_i};[0,1])
=\frac1p\beta^{1-p}<\infty=\Psi_p(S_0;[0,1]).
$$
At $\beta=0$, we have $\infty$ on both sides.

\textbf{Varicurrent convergence.}

We next compute the varicurrent limit by testing against $g\in C([0,1]\times\overline\Omega\times\orient)$. Each moving line has space-time area density $\ell_{\Delta t,N}/((2N+1)\Delta t)$ in $(t,x)$ coordinates. As in the current calculation, we change to the height variable using $t=(y_j^\pm)^{-1}(y)$. This gives
\begin{align*}
&\langle\eym{S_{\Delta t,N}^{\mathrm m}},g\rangle\\
&\quad=\ell_{\Delta t,N}\sum_{j=0}^N
\int_{\frac{2j}{2N+1}}^{\frac{2j+1}{2N+1}}\int_0^1
g\bigl((y_j^+)^{-1}(y),x,y,\eta_N^+\bigr)\d x\d y\\
&\qquad+\ell_{\Delta t,N}\sum_{j=1}^N
\int_{\frac{2j-1}{2N+1}}^{\frac{2j}{2N+1}}\int_0^1
g\bigl((y_j^-)^{-1}(y),x,y,\eta_N^-\bigr)\d x\d y\\
&\quad=\ell_{\Delta t,N}\int_0^1\int_0^1
g\Biggl(\Delta t\zigzag_{\frac{2}{2N+1}}(y),x,y,\eta_N^{\sigma_N}\Biggr)\d x\d y.
\end{align*}

Along the sequence, $\eta_{N_i}^\pm\to\eta_\beta^\pm$ as defined in the continued Example \ref{ex:thousand-cuts-current} in Section \ref{ssec:examples-revisited}. Uniform continuity of $g$ and $\Delta t_i\to 0$ give
$$
g(\Delta t_i \zigzag_{\frac{2}{2N_i+1}}(y),x,y,\eta_{N_i}^\pm)
\to g(0,x,y,\eta_\beta^\pm)\quad\text{uniformly on }(x,y)\in[0,1]^2.
$$
As $\sigma_{N_i}$ oscillates between values $1$ and $-1$ on $N_i+1$ and $N_i$ intervals, respectively, of equal length $1/(2N_i+1)$, we obtain
\[
\lim_{i\to\infty}\langle\eym{S_{\Delta t_i,N_i}^{\mathrm m}},g\rangle
=\frac{\ell_\beta}2\int_0^1\int_0^1
\bigl(g(0,x,y,\eta_\beta^+)+g(0,x,y,\eta_\beta^-)\bigr)\d x\d y =\langle\mu_\beta^{\mathrm m},g\rangle.
\]
Taking $g\equiv1$ in the preceding identity gives the limiting mass $\ell_\beta$, since the spatial unit square has area one and the two orientation weights sum to $\ell_\beta$. Moreover,
$$
\frac{\ell_\beta}2(\eta_\beta^++\eta_\beta^-)=e_2\wedge e_1,
$$
which verifies $\stcurr{\mu_\beta^{\mathrm m}}=S_0^{\mathrm m}$.

The stationary varicurrents satisfy
$$
\langle\eym{S_{\Delta t_i,N_i}^{\mathrm s}},g\rangle
=\int_{\Delta t_i}^1\int_0^1g(t,x,1,e_0\wedge e_1)\d x\d t
\longrightarrow\langle\eym{S_0^{\mathrm s}},g\rangle.
$$
Combining the moving and stationary varicurrents, we obtain the full limit
$$
\eym{S_{\Delta t_i,N_i}}\weakstarto\mu_\beta^{\mathrm m}+\eym{S_0^{\mathrm s}}=\mu_\beta.
$$

For $\beta>0$, the orientations satisfy
$$
|\eta_\beta^\pm \mres\mathrm dt|=\frac{\beta}{\ell_\beta},\qquad
|v(\eta_\beta^\pm)|=\frac1{\beta}.
$$
Using the total moving mass in the extended cost gives
\[
\Psi_p(\mu_\beta^{\mathrm m};[0,1])
=\frac1p\left(\frac1{\beta}\right)^p
\frac{\beta}{\ell_\beta}\ell_\beta =\frac1p\beta^{1-p}
=\lim_{i\to\infty}\Psi_p(\eym{S_{\Delta t_i,N_i}^{\mathrm m}};[0,1]).
\]
The measure $\mu_\beta^{\mathrm m}$ satisfies \eqref{eqn:NC-varicurrent} but is concentrated at time zero, so its cost is not a Lebesgue-time slice integral. Adding the stationary part changes neither the approximating nor the limiting cost.

For $\beta=0$, the two orientation atoms coincide at $e_2\wedge e_1$ with combined weight one. Thus $\mu_0^{\mathrm m}=\eym{S_0^{\mathrm m}}$ fails \eqref{eqn:NC-varicurrent}, and
$$
\lim_{i\to\infty}\Psi_p(\eym{S_{\Delta t_i,N_i}^{\mathrm m}};[0,1])
=\lim_{i\to\infty}\frac1p\bigl((2N_i+1)\Delta t_i\bigr)^{1-p}
=\infty=\Psi_p(\mu_0^{\mathrm m};[0,1]).
$$

\bibliographystyle{alpha}
\bibliography{biblio}

\end{document}